\documentclass[a4paper]{article}

\usepackage[OT1]{fontenc}

\usepackage[english]{babel}

\usepackage[utf8]{inputenc}

\usepackage[sc]{mathpazo}

\usepackage{amsmath,amssymb,amsfonts}
\usepackage{supertabular}
\usepackage{amscd}
\usepackage{textcomp}
\usepackage{authblk}
\usepackage{accents}

\usepackage{latexsym}
\usepackage{mathtools}
\usepackage{mathdots}
\usepackage[only,llbracket,rrbracket]{stmaryrd}

\usepackage[amsmath,thmmarks]{ntheorem}

\usepackage{graphicx}

\usepackage{soul}

\usepackage{pdfpages}

\usepackage{longtable}

\usepackage{varioref}

\usepackage{datetime}

\usepackage[h]{esvect}

\usepackage{array}

\usepackage{booktabs}

\usepackage[mathscr]{euscript}

\usepackage{mathtools}

\usepackage{tikz-cd}

\numberwithin{equation}{section}

\newtheorem{theorem}{Theorem}[section]

\newtheorem{re}[theorem]{Remark}
\newtheorem{co}[theorem]{Corollary}
\newtheorem{defi}[theorem]{Definition}
\newtheorem{lem}[theorem]{Lemma}
\newtheorem{pr}[theorem]{Proposition}

\theoremstyle{nonumberplain}
\theorembodyfont{\normalfont}
\theoremsymbol{\ensuremath{\square}}
\newtheorem{proof}{Proof}

\def\o{\circ}

\def\Lip{{\mathcal{L}ip}}

\def\Nz{\mathbb N}

\def\Rz{\mathbb R}

\def\Rz{\mathbb R}

\def\Wz{\mathbb W}

\def\Zz{\mathbb Z}

\def\ep{\varepsilon}

\def\Cc{\mathcal C}

\def\Ec{\mathcal E}

\def\Hc{\mathcal H}

\def\Lc{\mathcal L}
\def\Mc{\mathcal M}

\def\Sc{\mathcal S}

\def\Wc{\mathcal W}
\def\Xc{\mathcal X}

\def\Ef{\mathfrak E}

\def\Pf{\mathcal{}\mathfrak P}

\def\pf{\mathfrak p}

\def\ol#1{\overline{#1}}
\def\0{{\mathbf 0}}
\def\1{{\mathbf 1}}

\def\d{\partial}
\newcommand{\der}{\operatorname{d}}

\def\o{\circ}

\def\supp{{\operatorname{supp}}}

\def\eps{\varepsilon}
\def\al{\alpha}
\def\be{\beta}
\def\ga{\gamma}
\def\de{\delta}

\def\la{\lambda}

\def\rh{\rho}

\def\ph{\varphi}

\def\Ph{\Phi}

\def\ol#1{\overline{#1}}

\def\aequiv#1{\;\raise-4pt\hbox{$\substack{{\displaystyle\sim}\\\raise4pt\hbox{$\scriptstyle #1$}}$}
\;}

\DeclareMathAccent{\wtilde}{\mathord}{largesymbols}{"65}

\makeatletter
\def\cn#1{\def\tmp@cn{#1}\ifx\tmp@cn\empty\nabla\else\nabla^{\tmp@cn}\fi}
\makeatother

\def\div{\operatorname{div}}

\def\o{\circ}

\def\Xct{\Xc^{2,\infty}}

\usepackage[linkcolor=black,colorlinks=true,citecolor=black,filecolor=black]{hyperref}

\usepackage[numbers,sort&compress]{natbib}

\title{A convergent Lagrangian scheme for the Navier–Stokes equations on the normed Wasserstein space}
\author{Arthur Schichl}
\affil{\small KOF Swiss Economic Institute\\ETH Zürich\\Leonhardstrasse 21\\8092 Zürich\\ Switzerland\\
aschichl@ethz.ch}
\date{February 3, 2026}

\begin{document}
\def\d{\partial}
\def\o{\circ}


\maketitle

\begin{abstract}
We extend the differential structure on Wasserstein spaces of \citet{Schichl2025} to the normed Wasserstein space $\Wz_1$ of finite signed measures of \citet{PiccoliRossiTournus2023},
the natural setting for measure differential equations with source terms. We generalize the averaged flow dispersion and expectational flow operators
to $\Wz_1$ and show that along divergence--free tangent vectors they act by measure--preserving flows, and further define source operators. As an application, we construct a purely Lagrangian splitting 
scheme for the incompressible Navier--Stokes equations on $\Rz^N$, inspired by the Averaged Flow Approximation Series (AFAS) in \cite{Schichl2025}. We recover, with explicit constants, the classical lifespan bound, a blow--up rate and the inviscid 
limit with rate $O(\eta^2)$, and a linearised scheme yields $L^{\infty}$, Hölder, $L^1$ and moment estimates for classical solutions.
\end{abstract}

\newpage
\tableofcontents

\section{Introduction}\label{ch:intro}
In \citet{Schichl2025}, the Wasserstein spaces $\Wc_p(M)$ of $p$--integrable probability measures on a Riemannian manifold $(M,g)$, $p>2$, were endowed with a new differential structure, which is finer than the classical one induced by the continuity equation. 
The central objects of that theory are \emph{probability vector--fields}, i.e., maps sending a measure $\mu \in \Wc_p(M)$ onto a probability measure $V[\mu]$ on the space $\Xc^{2,\infty}(M)$ of vector--fields, and the associated non--linear and possibly degenerate parabolic measure differential equation, which reads in its weak form as,
\begin{equation}
    \tag{MDE}
    \begin{aligned}
        \int_M \ph \der \mu(s) &- \int_M \ph \der \mu(0) \\&= 
        \frac{1}{2}\int_0^s \int_M \int_{\Xct(M)} \Lc_{X-\ol{V}[\mu(t)]}^2(\ph) \der V[\mu(t)] \der \mu(t) \der t \\ &\qquad+ \int_0^s \int_M \Lc_{\ol{V}[\mu(t)]}(\ph) \der \mu(t) \der t,
    \end{aligned}
\end{equation} 
where $\ol{V}$ is the average vector--field of $V$ defined in \cite[Definition 3.8]{Schichl2025}. Under weak assumptions on $V$, solutions of (MDE) were constructed as uniform limits of an explicit Euler--type scheme, the \emph{Average Flow Approximation Series} (AFAS), built from the \emph{averaged flow dispersion operators} and the \emph{expectational flow operators}. 

However, the theory of \citet{Schichl2025} is intrinsically restricted to probability measures. The aim of this work is to develop a similar theory for evolution equations carrying a source term, i.e., equations whose solutions do not conserve mass. To this end, instead of the Wasserstein spaces $\Wc_p(M)$, we work on a space of signed measures, namely the \emph{normed Wasserstein space} $\Wz_1(M)$ of Section~\ref{ch:diffprelim}, which is the space of finite signed Borel--measures on $M$ endowed with the $1$--Wasserstein--norm. This norm was first introduced by Piccoli and his coauthors \cite{PiccoliRossi2014, PiccoliRossiTournus2023}, precisely in order to study non--local transport equations with source terms, see also \cite{Piccoli2019} for the related notion of measure differential equations. 
While the $\Wz_1$--norm is not equivalent to the standard Wasserstein--$1$--distance and $\Wz_1(M)$ is in general not even complete, we shall see in Section~\ref{ch:diffprelim} that many topological properties of $\Wc_1(M)$ carry over to $\Wz_1(M)$ or hold on suitable subsets. This allows us to reproduce convergence arguments of \cite{Schichl2025} in the present setting.

In Section~\ref{ch:parstructwas}, we illustrate the scope of the theory of \citet{Schichl2025} on normed Wasserstein spaces by building a convergent approximation scheme for solutions of a system of particularly prominent non--linear parabolic PDEs: the incompressible Navier--Stokes equations. We introduce the analogues of the averaged flow dispersion operators and of the expectational flow operators of \cite{Schichl2025} on $\Wz_1(M)$, 
together with two new families of operators: the \emph{source operators}, which encode a general source term and the \emph{dievergence operators}, which encodes the specific source term for the incompressible Navier--Stokes equations. We then focus on divergence--free tangent vectors, for which the expectational flow operators act on densities by composition with a measure--preserving flow. This leads to particularly strong estimates on the operators and allows us to build Euler--type schemes---the \emph{Incompressible Flow Approximation Series} (IFAS) and their linearised version (LIFAS)---for the incompressible Navier--Stokes equations, which we write in the form \eqref{eq:NS} below. 
Our first main result, Theorem~\ref{th:LIFAS}, states that the LIFAS associated to a classical solution of \eqref{eq:NS} converge uniformly to that solution with respect to the $\Wz_1$--norm. As a by--product, we obtain a priori estimates for classical solutions of \eqref{eq:NS}. Finally, we discuss the convergence of the non--linear scheme (IFAS) in Theorem~\ref{th:IFAS}. In particular, we show that the (IFAS) uniformly converge to a solution of the incompressible Navier--Stokes equations in the $\Wz_1$-- and $\Cc^{2}$--norms on a given time interval, provided that the initial condition lies in a sufficently regular weighted Hölder--space. 

Two connections should be pointed out. 
First, the intended scope of the framework is the class of measure differential equations with source terms studied in \cite{PiccoliRossi2014, Piccoli2019, PiccoliRossiTournus2023}, such as transport--diffusion equations with birth and death terms.
Second, the non--linear scheme we build for the Navier--Stokes equations is a purely Lagrangian splitting scheme, in which the diffusion is realised by a symmetric average of translations---the expectation over one step of a random walk---the transport by composition with the flow of the current velocity, and the pressure by a Leray projection. 
It is thus the Euler--step analogue of the stochastic Lagrangian representation of the Navier--Stokes equations of \citet{ConstantinIyer2008}, and from the numerical point of view, a projection method of Chorin type with semi--Lagrangian advection and random--walk diffusion \cite{Chorin1968, Chorin1973}. Our convergence result is accordingly
uniform in the viscosity and includes the Euler equations.

\section{Preliminaries}\label{ch:diffprelim}

In this preliminary section, we introduce the objects and notions needed to provide a suitable setting for a generalization of the setting of \citet{Schichl2025} and the study of solutions of equations of the incompressible Navier--Stokes equations in this new framework.

First we define the normed Wasserstein space $\Wz_1(M)$, which is the space of finite signed Borel--measures on a Polish space $(M,d)$ endowed with the $1$--Wasserstein--norm. This norm was first introduced by \citet{PiccoliRossi2014, PiccoliRossiTournus2023}.
While it is not equivalent to the standard Wasserstein--$1$--distance, we shall see that many topological properties of $\Wc_1(M)$ carry over to $\Wz_1(M)$ or hold on suitable subsets.  

Then, we introduce the notion of tightness in the $p$--th moment, which is equivalent to relative compactness in $\Wc_p(M)$, and prove some further interesting properties of such families of measures in $\Mc_p(M)$. 

Finally, we include a short collection of known results on weighted Hölder--spaces, which are the natural functional spaces for the study of solutions of the incompressible Navier--Stokes equations. 
In particular, we establish relations between the $\Wz_1$--norm and the weighted Hölder--norms and apply the resulting estimates to standard objects of differential geometry, which shall be used 
to obtain bounds for the evolution operators defined in Section~\ref{sbsec:flowestims}.

\subsection{Notation}

The basic number sets will be denoted by $\Nz$ (natural numbers, $0 \in \Nz$), $\Zz$ (integers) and $\Rz$ (real numbers), and we will denote by $\Rz_+$ the set of all positive real numbers. If $a,b \in \Zz$, we will write $\llbracket a, b \rrbracket \coloneqq \{a, a+1, \dots, b\}$. For $N \in \Nz$, $(e_i)_{i \in \llbracket 0, N-1 \rrbracket}$ denotes the canonical basis of $\Rz^N$, $\|\cdot\|$ the Euclidean norm, $B(x,r)$ the open ball of radius $r$ around $x$, and $\delta_x$ the Dirac measure at $x$.

Let $\phi : \Rz^N \to \Rz^m$ be a differentiable function. We will denote by $D\phi(x)$ the Jacobian of $\phi$ at $x \in \Rz^N$, by $\nabla^k \phi$ its $k$--th derivative, and by $\|A\|$ the operator norm of a linear map $A$. For a measurable map $\phi: M \to M'$ between measurable spaces and a measure $\mu$ on $M$, $\phi_*\mu \coloneqq \mu \circ \phi^{-1}$ denotes the push--forward of $\mu$. The Fourier transform of a tempered distribution $f$ on $\Rz^N$ is denoted by $\hat f$. For $T > 0$, a \emph{partition} $\Pf$ of $[0,T]$ is a finite subset of $[0,T]$ containing $0$ and $T$; we write $|\Pf|$ for its cardinality and $\pf(\Pf)$ for its maximal step--length, i.e., the maximal distance between two consecutive elements of $\Pf$.

We will need several spaces of measures and real function spaces throughout this work. Let $(Z,d)$ be a Polish space with a fixed base point $x_0 \in Z$, $(M,g)$ a complete Riemannian manifold with volume form $\la$, and $A \subset \Rz^N$ be closed. We define (see \citet{Schichl2025} and the references therein for formal introductions):

\bigskip

\begin{longtable}{p{3.2cm}p{8.2cm}}
    $\Mc(Z)$ & The space of finite signed Borel--measures on $Z$. For $\mu \in \Mc(Z)$, $\mu = \mu^+ - \mu^-$ denotes its Jordan decomposition, $|\mu| \coloneqq \mu^+ + \mu^-$ its total variation measure and $|\mu|(Z)$ its mass.\\
    $\Mc_+(Z)$ & The set of finite positive Borel--measures on $Z$.\\
    $\Mc_p(\mu)$ & For $p > 0$ and $\mu \in \Mc(Z)$, the $p$--th moment $\int_Z d(x,x_0)^p \der |\mu|(x)$ of $\mu$. For $Z = \Rz^N$ we take $x_0 = 0$.\\
    $\Mc_p(Z)$ & For $p > 0$ the space of all $\mu \in \Mc(Z)$ with $\Mc_p(\mu) < \infty$.\\
    $\Wc_p(Z)$ & The Wasserstein space of $p$--integrable probability Borel--measures on $Z$ for $p \geq 1$.\\
    $\Wc_p(\mu,\nu)$ & The $p$--Wasserstein distance between $\mu, \nu \in \Wc_p(Z)$, for $p \geq 1$.\\
    $\Cc(A, \Rz^m)$ & The set of all continuous functions $f : A \to \Rz^m$. We write $\Cc_b(A,\Rz^m)$ for the bounded ones and $\Cc_0(A, \Rz^m)$ for those vanishing at infinity. For a normed space $E$ and $T > 0$, $\Cc([0,T], E)$ denotes the set of continuous curves $[0,T] \to E$.\\
    $\Cc^k(A, \Rz^m)$ & For $k \in \Nz \cup \{\infty\}$ the set of all $k$--times continuously differentiable functions $f : A \to \Rz^m$.\\
    $\Cc^k_c(A, \Rz^m)$ & For $k \in \Nz \cup \{\infty\}$ the set of all $f \in \Cc^k(A,\Rz^m)$ with $\supp f$ compact.\\
    $\Cc^k_b(A, \Rz^m)$ & For $k \in \Nz \cup \{\infty\}$ the set of all $f \in \Cc^k(A, \Rz^m)$ where $f$ and all of its derivatives up to order $k$ are bounded, endowed with the norm $\|f\|_{\Cc^k} \coloneqq \sum_{j \leq k}\|\nabla^j f\|_{\infty}$.\\
    $\Cc^{\al}(A, \Rz^m)$ & For $\al \in (0,1]$ the set of all bounded $\al$--Hölder continuous functions $f: A \to \Rz^m$. We write $\{f\}_{\Cc^{\al}} \coloneqq \sup_{x \neq y}|f(x) - f(y)|/\|x - y\|^{\al}$ for the Hölder seminorm and $\|f\|_{\Cc^{\al}} \coloneqq \|f\|_{\infty} + \{f\}_{\Cc^{\al}}$ for the Hölder norm. $\dot{\Cc}^{\al}(A, \Rz^m)$ denotes the homogeneous Hölder space, i.e., the set of all $f$ with $\{f\}_{\Cc^{\al}} < \infty$.\\
    $\Cc^{k,\al}(A, \Rz^m)$, $\Cc^{k,\al}_b(A,\Rz^m)$ & For $k \in \Nz$ and $\al \in (0,1]$ the set of all $f \in \Cc^k(A,\Rz^m)$, resp. $f \in \Cc^k_b(A,\Rz^m)$, such that $\nabla^k f \in \dot{\Cc}^{\al}$, the latter endowed with $\|f\|_{\Cc^{k,\al}} \coloneqq \|f\|_{\Cc^k} + \{\nabla^k f\}_{\Cc^{\al}}$.\\
    $\Lip(A, \Rz^m)$ & The set of all Lipschitz continuous functions $f : A \to \Rz^m$, with $\|f\|_{\Lip}$ the Lipschitz constant of $f$. For $f \in \Cc^1_b$, $\|f\|_{\Lip} = \|\nabla f\|_{\infty}$.\\
    $L^p(A, \Rz^m)$ & For $1 \leq p < \infty$ the space of all (equivalence classes of) measurable functions $f : A \to \Rz^m$ for which $\int_A |f|^p \der \la$ converges.\\
    $L^{\infty}(A, \Rz^m)$ & The set of all (equivalence classes of) measurable functions $f: A \to \Rz^m$ that are essentially bounded. We write $\|f\|_{\infty}$ for the essential supremum.\\
    $W^{k,p}(A, \Rz^m)$ & For $1 \leq p \leq \infty$ and $k \in \Nz$ the space of functions $f \in L^p(A, \Rz^m)$ which are weakly $k$--times differentiable, and for which all derivatives up to order $k$ are also in $L^p$.\\
    $H^k(A, \Rz^m)$ & $W^{k,2}(A, \Rz^m)$.\\
    $\Hc^1(\Rz^N, \Rz^m)$ & The real Hardy space on $\Rz^N$ (see \cite[Chapter III]{Stein1993}).\\
    $\Sc(\Rz^N, \Rz^m)$ & The Schwartz space of rapidly decreasing smooth functions $f : \Rz^N \to \Rz^m$.\\
    $\Xc(M)$ & The space of all smooth vector--fields on $M$.\\
    $\Xc^{k,\infty}(M)$ & For $k \in \Nz$, the space of vector--fields on $M$ which are essentially bounded, weakly $k$--times differentiable, and for which all derivatives up to order $k$ are also essentially bounded. We endow it with the classical $W^{k,\infty}$--norm.\\
    $\rh^X_t$ & For $X \in \Xc^{1,\infty}(M)$ and $t \in \Rz$, the flow of $X$ at time $t$, i.e., $\frac{\der}{\der t}\rh^X_t(x) = X(\rh^X_t(x))$ and $\rh^X_0 = \mathrm{id}_M$.\\
    $\Lc_X \ph$ & For $X \in \Xc^{1,\infty}(M)$ and $\ph \in \Cc^1(M)$, the Lie derivative $\der\ph(X)$ of $\ph$ along $X$.\\
    $\div X$ & The divergence of $X \in \Xc^{1,\infty}(M)$ with respect to $\la$. On $\Rz^N$, $\Delta$ denotes the Laplacian, $\Gamma$ the fundamental solution of $\Delta$, normalised by $\Delta\Gamma = \delta_0$, and $\Delta^{-1}f \coloneqq \Gamma * f$ whenever the convolution is defined.\\
    $(R_j)_{j \in \llbracket 0, N-1 \rrbracket}$ & The Riesz transforms on $\Rz^N$ (see \cite[Chapter III, \S 1]{Stein1970}).\\
\end{longtable}

\bigskip

When the target space $\Rz^m$ or the domain is clear from the context, we omit it from the notation, e.g., $\Cc^{\al}(\Rz^N) = \Cc^{\al}(\Rz^N, \Rz)$. Signed measures with a density with respect to $\la$ are identified with their densities. In particular, $L^1(\Rz^N) \subset \Mc(\Rz^N)$, and a vector--field $\nu \in L^1(\Rz^N, \Rz^N)$ is identified with the $N$--tuple $(\nu_i)_{i \in \llbracket 0, N-1 \rrbracket}$ of signed measures given by its components, on which all operators on $\Mc(\Rz^N)$ and all norms of measures are applied componentwise. Unless stated otherwise, $C$ denotes a positive constant depending only on the quantities indicated in brackets, which may change from line to line.

\subsection{The normed Wasserstein space $\Wz_1$}

The normed Wasserstein space $\Wz_1(M)$ was introduced by Piccoli and his coauthors \cite{PiccoliRossi2014, PiccoliRossiTournus2023}.

\begin{defi}[$1$--Wasserstein--norm]\label{de:Wassnorm}
   Let $(M,d)$ be a Polish space. Define for all $\mu \in \Mc(M)$,
   \begin{equation}
      \|\mu\|_{\Wz_1} \coloneqq \sup \biggl\{ \int f \der \mu \mid f \in L^{\infty}(M) \cap \Lip(M) : \max(\|f\|_{\infty}, \|f\|_{\Lip}) \leq 1 \biggr\}.
   \end{equation}
   Then, $\Wz_1(M) \coloneqq \bigl(\Mc(M), \|\|_{\Wz^1}\bigr)$ is a normed vector--space. 
\end{defi}
While $\Wz_1(M)$ is in general not a Banach space, convergence of Cauchy--sequences 
can be recovered under relatively mild assumptions.
\begin{pr}
   Let $(\mu_n)_{n \in \Nz}$ be a Cauchy--sequence in $\Wz_1(M)$. Suppose that $(\mu_n)_{n \in \Nz}$ is tight and that there exists $B>0$
   such that $\mu_n(M) \leq B$ for all $n \in \Nz$. Then, $(\mu_n)_{n \in \Nz}$ converges in $\Wz_1(M)$.
\end{pr}
\begin{proof}
   See \cite[Theorem 29]{PiccoliRossiTournus2023}, where the statement is proved for the equivalent norms $\|\cdot\|^{a,b}$ on $\Mc(\Rz^d)$. The proof carries over verbatim to Polish spaces.
\end{proof}
Despite its name, the $1$--Wasserstein--norm of Definition~\ref{de:Wassnorm} is not equivalent to the classical $1$--Wasserstein distance on 
on $\Wc_1(M)$. Convergence in $\Wc_1$ is a stronger property than convergence in $\Wz_1$.
\begin{lem}
   Let $\mu,\nu \in \Wc_1(M)$. Then, $\|\mu - \nu\|_{\Wz_1} \leq \Wc_1(\mu, \nu)$.
\end{lem}
While convergence in $\Wc_1$ implies narrow convergence, there exists sequences that do not converge narrowly, but converge in $\Wz_1(M)$, see \cite[Remarks 24 and 25]{PiccoliRossiTournus2023}. However, we recover the following property: 
\begin{lem}\label{le:weakerconv}
   Let $\mu \in \Wz_1(M)$ and let $(\mu_n) \in \Wz_1(M)^{\Nz}$ such that $\lim_{n \to \infty} \mu_n = \mu$ in $\Wz_1(M)$.
   Then, for all $\ph \in \Cc_b(M) \cap \Lip(M)$, 
   \begin{equation}
      \lim_{n \to \infty} \int_M \ph \der \mu_n = \int_M \ph \der \mu.
   \end{equation} 
\end{lem}
 Lemma~\ref{le:weakerconv} does not imply pointwise convergence on arbitrary test functions. However, the statement is strong enough to deduce that limits of sequences converging with respect to the $\Wz_1$--norm preserve
 important properties of the sequence. For the coming results, we take $N \in \Nz$ and work with $M = \Rz^N$.
\begin{pr}
   Let $(\mu_n)_{n \in \Nz}$ be a sequence in $\Wz_1(\Rz^N)$ such that $\lim_{n \to \infty} \mu_n = \mu$ in $\Wz_1(\Rz^N)$. Then, 
   \begin{itemize}
      \item if for all $n \in \Nz$, $\mu_n \in L^{\infty}(\Rz^N)$ and there exists $B > 0$ such that for all $n \in \Nz$, $\|\mu_n\|_{L^{\infty}} \leq B$, then, $\mu \in L^{\infty}(\Rz^N)$ and $\|\mu\|_{L^{\infty}} \leq B$.
      \item if there exist $\al \in (0,1]$ and $B >0$, such that for all $n \in \Nz$, $\mu_n \in \Cc^{\al}(\Rz^N)$ and $\|\mu_n\|_{\Cc^{\al}} \leq B$, then $\mu \in \Cc^{\al}(\Rz^N)$ and $\|\mu\|_{\Cc^{\al}} \leq B$.
   \end{itemize}
\end{pr}
\begin{proof}
   Let us prove the first point. For all $\ep > 0$, let $\ph_{\ep} \in \Cc^{\infty}(\Rz^N)$ be a standard mollifier with support in $B(0,\ep)$. In particular, $\ph_{\ep} \in \Cc_b(\Rz^N) \cap \Lip(\Rz^N)$. Hence, for all $x \in \Rz^N$,
   \begin{equation}
      \lim_{n \to \infty} \int_{\Rz^N} \ph_{\ep}(x-y) \der \mu_n(y) = \int_{\Rz^N} \ph_{\ep}(x-y) \der \mu(y).
   \end{equation}
   Now, for all $x \in \Rz^N$,
   \begin{equation}
      \int_{\Rz^N} \ph_{\ep}(x-y) \der \mu_n(y) \leq B. 
   \end{equation}
   Hence, for all $\ep >0$,  $\|\ph_{\ep} * \mu\|_{\infty} \leq B$.  Therefore, $\mu \in L^{\infty}(\Rz^N)$ and $\|\mu\|_{L^{\infty}} \leq B$. 

   For the second point, denote for all $x,y \in \Rz^N$ with $x \neq y$ by 
   \begin{equation}
      \ph_{\ep}^{x,y}: z \mapsto \ph_{\ep}(z - x) - \ph_{\ep}(z - y).
   \end{equation}
   Let $\de >0$.  Since there exists $B>0$, such that for all $n \in \Nz$, $\|\mu_n\|_{C^{\al}} \leq B$, in particular, if we choose $\ep$ small enough:
   \begin{itemize}
      \item $\|\ph_{\ep}^{x,y}\|_{L^1} = 2$
      \item for all $n \in \Nz$, $\biggl| \int_{\Rz^n} \ph_{\ep}^{x,y} \der \mu_n - \bigl(\mu_n(x) - \mu_n(y)\bigr)\biggr| \leq \de$.
   \end{itemize} 
   In particular, it follows that for all $x,y \in \Rz^N$ with $x \neq y$,
   \begin{equation}
      \int_{\Rz^n} \ph_{\ep}^{x,y} \der \mu = \lim_{n \to \infty} \int_{\Rz^n} \ph_{\ep}^{x,y} \der \mu_n \leq B|x-y|^{\al} + \de.
   \end{equation}
   Hence, $\mu \in \Cc^{\al}(\Rz^N)$ and $\|\mu\|_{\Cc^{\al}} \leq B$.
\end{proof}

On some functional spaces of stronger regularity that are naturally embedded into $\Wz_1(\Rz^N)$, the passage from convergence in $\Wz_1(\Rz^N)$ to convergence in the stronger norms rests on the following interpolation inequalities.

\begin{lem}[Interpolation]\label{le:interpol}
    Let $N \geq 1$, $\al \in (0,1)$ and $c_N \coloneqq \int_{B(0,1)} 1 - \|z\| \der z$. Then:
    \begin{itemize}
        \item for all $\nu \in W^{1,\infty}(\Rz^N) \cap \Wz_1(\Rz^N)$,
        \begin{equation}\label{eq:interpol}
            \|\nu\|_{\infty} \leq \max\biggl( \frac{2}{c_N}\|\nu\|_{\Wz_1}, \Bigl(\frac{2(4\|\nu\|_{\Lip})^{N+1}}{c_N}\|\nu\|_{\Wz_1}\Bigr)^{\frac{1}{N+2}}\biggr);
        \end{equation}
        \item there exists $C(N,\al) > 0$ such that for all $u \in \Cc^{2,\al}_b(\Rz^N)$,
        \begin{equation}\label{eq:interpolholder}
            \|\nabla u\|_{\infty} + \|\nabla^2 u\|_{\infty} \leq C(N,\al)\Bigl(\|u\|_{\infty} + \|u\|_{\infty}^{\frac{\al}{2+\al}}\{\nabla^2 u\}_{\Cc^{\al}}^{\frac{2}{2 + \al}}\Bigr).
        \end{equation}
    \end{itemize}
    In particular, if $(\nu_n)_n$ is bounded in $\Cc^{2,\al}_b(\Rz^N)$ and converges in $\Wz_1(\Rz^N)$, then it converges in $\Cc^2_b(\Rz^N)$.
\end{lem}
\begin{proof}
   For the first point, let $a < \|\nu\|_{\infty}$ and w.l.o.g $x_0 \in \Rz^N$ with $\nu(x_0) \geq a$.
   Let $\rh \coloneqq \min\bigl(1, \frac{a}{2\|\nu\|_{\Lip}}\bigr)$. Then, $\nu \geq \frac{a}{2}$ on $B(x_0, \rh)$, and testing the $\Wz_1$--norm with 
   \begin{equation}
      \begin{aligned}
         f: \Rz^N &\to \Rz\\
         x &\mapsto \max(0, \rh - \|x - x_0\|),
      \end{aligned}
   \end{equation}
    which satisfies $\max(\|f\|_{\infty}, \|f\|_{\Lip}) \leq 1$, gives $\|\nu\|_{\Wz_1} \geq \frac{a}{2}c_N \rh^{N+1}$. If $\rh = 1$, this gives $a \leq \frac{2}{c_N}\|\nu\|_{\Wz_1}$. Otherwise, $a^{N+2} \leq 2(2\|\nu\|_{\Lip})^{N+1}c_N^{-1}\|\nu\|_{\Wz_1}$. Hence, \eqref{eq:interpol} follows. 
    
   The second point is the interpolation inequality for Hölder norms \cite[Lemma 6.32]{GilbargTrudinger2001}, applied on balls of radius one. The last assertion follows by applying both points to $\nu_n - \nu_m$.
\end{proof}

While the topologies induced by $\Wz_1$ and $\Wc_1$ have fundamentally different properties, both distances are related by a relatively simple formula.
\begin{theorem}
   Let $\mu \in \Wz_1(M)$. Then,
   \begin{equation}
      \|\mu\|_{\Wz_1} = \inf \biggl\{ \biggl(|\mu^+ - \nu| + |\mu^- - m| + \Wc_1(\nu, m)\biggr) \mid m,\nu \in \Mc_+(M): |m| = |\nu|\biggr\}.
   \end{equation}
\end{theorem}
\begin{proof}
   In \cite{PiccoliRossiTournus2023}, the right--hand side is the definition of the norm $\|\mu\|^{1,1}$, see \cite[Definitions 11 and 19]{PiccoliRossiTournus2023}, and the equality with the supremum of Definition~\ref{de:Wassnorm} is the Kantorovich--type duality \cite[Proposition 23]{PiccoliRossiTournus2023}.
\end{proof}

\subsection{Compact subsets of Wasserstein spaces}

First, we introduce the concept of tightness in the $p$--th moment, which, as we shall see, is equivalent to relative compactness.

\begin{defi}
\small
Let $p > 0$ and $(M,d)$ be a metric space. Fix $x_0 \in M$. A family $(\mu_i) \in \Mc_p(M)^I$ is called tight in the $p$--th moment if for all $\epsilon > 0$, there exists a compact set $K$ such that for all $i\in I$, $\int_{M\setminus K} d(x,x_0)^p \der \bigl|\mu_i\bigr|(x) \leq \epsilon$.
\end{defi}

\begin{lem}\label{le:q-tight}
Let $(M,d)$ be a metric space. Take $x_0 \in M$. Let $(\mu_i)_{i \in I}$ be a family of measures such that there exists $B > 0$ such that for all $i\in I$, $\Mc_p(\mu_i) \leq B^p$. Then for all $R > 0$ and all $i \in I$, 
\begin{equation}
\bigl|\mu_i\bigr|(M \setminus B(x_0, R)) \leq \min\biggl(\bigl|\mu_i\bigr|,\frac{B^p}{R^p}\biggr).
\end{equation}
Furthermore, if $(M,d)$ is proper, $(\mu_i)_{i \in I}$ is tight in the $q$--th moment for all $q < p$. More precisely, for $\ep > B^q\frac p{p-q}$ all $R\geq R_{\ep}$, with
\begin{equation}
R_{\ep} \coloneqq 
   \biggl(B^p\frac{p}{p-q}\biggr)^{\frac{1}{p-q}}\ep^{\frac{1}{q-p}} 
\end{equation}
we obtain
\begin{equation}
\int_{M \setminus B(x_0, R_{\ep})} d(x_0,x)^q \der \mu_i(x) \leq \ep.
\end{equation}
\end{lem}
\begin{proof}
The first point is a straightforward consequence of the Markov inequality.
Now, let $q<p$, $i\in I$, $\ep>0$ and $R\geq B$.
\begin{equation}
\begin{aligned}
&\int_{M \setminus B(x_0, R)} d(x_0,x)^q \der \bigl|\mu_i\bigr|(x) \\
&\qquad = \int_{R^q}^{\infty} \bigl|\mu_i\bigr|(d(x_0,x)^q>\la) \der \la + \int_0^{R^q} \bigl|\mu_i\bigr|(M\setminus \ol{B(x_0,R^q)}) \der \la\\
&\qquad\leq \int_{R^q}^{\infty} \bigl|\mu_i\bigr|(M\setminus B(x_0,\la^{\frac{1}{q}})) \der \la + R^q\frac{B^p}{R^p}\\
&\qquad\leq \int_{R^q}^{\infty} \frac{B^p}{\la^{\frac{p}{q}}} \der \la + \frac{B^p}{R^{p-q}} = \frac{B^p}{R^{p-q}}\biggl(\frac1{\frac{p}{q}-1}+1\biggr).
\end{aligned}
\end{equation}
Hence, if 
\begin{equation}
R \geq \biggl(\frac{B^pp}{p-q}\biggr)^{\frac{1}{p-q}}\ep^{\frac{1}{q-p}},
\end{equation}
we obtain 
\begin{equation}
\int_{M \setminus B(x_0, R_{\ep})} d(x_0,x)^q \der \bigl|\mu_i\bigr|(x) \leq \ep.
\end{equation}
Therefore, the $(\mu_i)_{i\in I}$ are tight in the $q$--th moment.
\end{proof}

The following characterization of compact subsets of Wasserstein spacesonly hold when the base space is proper.

\begin{lem}\label{le:tightconv}
Let $(M,d)$ be a proper metric space. Let $x_0 \in M$ and $p >0$. Let $(\mu_i) \in \Wc_p(M)^I$ be tight in the $p$--th moment. Then $(\mu_i)$ is tight. Furthermore, let $(i_n)_{n\in\Nz}$ be a sequence in $I$ such that $(\mu_{i_n})$ converges in the weak--*--topology. Then, $(\mu_{i_n})$ converges with respect to the $\Wz_1$--norm.
\end{lem}
\begin{proof}
Let $\eps > 0$ and take $K$ a compact set such that for all $i\in I$, 
\begin{equation}
\int_{M \setminus K} d(x,x_0)^p \der \bigl|\mu_i\bigr|(x) < \eps.
\end{equation}
Then, in particular $\int_{M \setminus (K \cup B(x_0,1))} d(x,x_0)^p \der \bigl|\mu_i\bigr|(x) < \eps$. Hence, we have 
\begin{equation}
\bigl|\mu_i\bigr|\bigl(M \setminus (K \cup B(x_0,1))\bigr) \leq \eps
\end{equation}
 and we conclude by compactness of $K \cup B(x_0,1)$ that $(\mu_i)$ is tight.

Now, suppose $(\mu_{k_n})_{n \in \Nz}$ is a sequence in $(\mu_i)$ that converges weak-* to a signed Borel--measure $\mu$. Then, since $\bigl(\mu_i\bigr)_{i \in I}$ is tight, $\bigl(\mu_{k_n}\bigr)_{n \in \Nz}$ converges narrowly to $\mu$, and, in particular, $\lim_{n \to \infty} \mu_{k_n}(M) = \mu(M)$.
Suppose $(\mu_{k_n})_{n \in \Nz}$ is not converging to $\mu$ in $\Wz_1(M)$, i.e., there exists a subsequence $(\mu_{i_n})_{n \in \Nz}$ of $(\mu_{k_n})_{n \in \Nz}$, such that for all $n \in \Nz$, $\|\mu - \mu_{i_n}\|_{\Wz_1} > 4\ep$.

Since  $(\mu_{i_n})_{n \in \Nz}$ is a tight family, there exists $R > 0$ such that $K \subset B(x_0,R)$ and for all $n \in \Nz$ a map $f_n \in L^{\infty}\bigl(B(x_0,R)\bigr) \cap \Lip\bigl(B(x_0,R)\bigr)$, such that $\max(\|f_n\|_{\infty}, \|f_n\|_{\Lip}) \leq 1$ and 
\begin{equation}
   \int_M f_n \der\bigl(\mu - \mu_{i_n}\bigr) \geq \|\mu - \mu_{i_n}\|_{\Wz_1} - \ep.
\end{equation}
Since $B(x_0, R)$ is compact and the $f_n$ are $1$--Lipschitz and uniformly bounded, the Theorem of Arzelà--Ascoli yields 
the existence of a subsequence $\bigl(f_{j_n}\bigr)_{n \in \Nz}$ which uniformly converges to a function $f$
supported on $B(x_0, R)$ satisfying $\max(\|f\|_{\infty}, \|f\|_{\Lip}) \leq 1$. In particular, there exists $N_0 \in \Nz$,
such that for all $n \geq N_0$, 
\begin{equation}
   \biggl| \int_M (f - f_{j_n}) \der\bigl(\mu - \mu_{i_{j_n}} \bigr) \biggr| \leq \ep.
\end{equation}
On the other hand, by narrow convergence 
\begin{equation}
   \lim_{n \to \infty} \int_M f \der \mu_{i_{j_n}} = \int_M f \der \mu.
\end{equation}
Thus, there exists $N_1 \geq N_0$, such that for all $n \geq N_1$,
\begin{equation}
   \biggl| \int_M f_{j_n} \der \bigl(\mu - \mu_{i_{j_n}}\bigr) \biggr| \leq 2\ep.
\end{equation}
In particular, $\|\mu - \mu_{i_n}\|_{\Wz_1} \leq 3\ep$, which yields a contradiction. Therefore, $\bigl(\mu_{k_n}\bigr)_{n \in \Nz}$ converges to $\mu$ in 
$\Wz_1(M)$.
\end{proof}

\begin{co}\label{co:compacttight}
	Let $p>0$ and $(\mu_i) \in \Wc_p(M)^I$ be tight in the $p$--th moment. Then $(\mu_i)$ is relatively compact in $\Wz_1(M)$.
\end{co}

\begin{re}
   \citet{Schichl2025} proves stronger analogues of the previous results for compact subsets of the classical Wasserstein--spaces $\Wc_p(M)$ for $p \geq 1$.
\end{re}

\subsection{Weighted Hölder spaces}\label{sbsec:weightedholder}

This section is dedicated to weighted Hölder spaces, the natural setting when working with the Navier-Stokes equations. We recall their definition and some basic properties we shall use in the further course of this work. Let $N \in \Nz$. For all $x \in \Rz^N$, we write
\begin{equation}
    \langle x \rangle \coloneqq \bigl(1 + \|x\|^2\bigr)^{\frac{1}{2}}.
\end{equation}
Note that $\langle \cdot \rangle$ is smooth, $1$--Lipschitz, and $\max(1, \|x\|) \leq \langle x \rangle \leq 1 + \|x\|$.

\begin{defi}[Weighted Hölder norms]\label{de:weightednorms}
    Let $N \geq 2$, $\al \in (0,1)$, $\be, \ga > 0$. 
    \begin{itemize}
    \item Let $E$ be a finite--dimensional normed vector space and $G: \Rz^N \to E$ be bounded. We define
    \begin{equation}
      \begin{aligned}
         |G|_{\be} &\coloneqq \sup_{x \in \Rz^N} \langle x \rangle^{\be}|G(x)|,\\  [G]_{\al, \ga} &\coloneqq \sup_{x \neq y} \frac{|G(x) - G(y)|}{\|x - y\|^{\al}\,\omega_{\ga}(x,y)},  
      \end{aligned}
    \end{equation}
    where $\omega_{\ga}(x,y) \coloneqq \frac{1}{2}\bigl(\langle x \rangle^{-\ga} + \langle y \rangle^{-\ga}\bigr)$ for all $(x,y) \in \bigl(\Rz^N \bigr)^2$.
    \item Fix $\ep \in (0,1)$. For all $\nu \in \Cc^2(\Rz^N, \Rz^N)$ we set
    \begin{equation}
        \|\nu\|_{\mathcal E_{\al,\ep}} \coloneqq \sum_{k = 0}^{2} |\nabla^k \nu|_{N + 1 + k - \ep} + [\nabla^2 \nu]_{\al, N + 3 + \al - \ep},
    \end{equation}
    and denote by $\mathfrak E_{\al,\ep}$ the space of all $\nu \in \Cc^2(\Rz^N, \Rz^N)$ with $\|\nu\|_{\mathcal E_{\al,\ep}} < \infty$, and by $\mathcal E_{\al,\ep} \subset \mathfrak E_{\al,\ep}$ the subset of divergence--free fields. 
   \end{itemize}
    To keep notation short, if $\al$ and $\ep$ are clear in the context, we shall simply write $\Ef$, $\Ec$ and $\|\cdot\|_{\Ec}$ instead. For $F \in \Cc^1(\Rz^N)$ we set
    \begin{equation}
        \|F\|_{\mathcal F} \coloneqq |F|_{N + 2 - \ep} + |\nabla F|_{N + 3 - \ep} + [\nabla F]_{\al, N + 3 + \al - \ep}.
    \end{equation}
\end{defi}

\begin{re}
Since $\omega_{\ga}(x,y) \leq \min(\langle x \rangle, \langle y \rangle)^{-\ga} \leq 2\omega_{\ga}(x,y)$, the seminorm $[\cdot]_{\al,\ga}$ is equivalent to a weighted Hölder seminorm with weight $\min(\langle x \rangle, \langle y \rangle)^{\ga}$. 
\end{re}

\begin{lem}\label{le:weights}
    Let $\ga > 0$. There exists $C = C(N,\ga) > 0$, non--decreasing in $\ga$, such that the following holds.
    \begin{itemize}
        \item For all $x, a \in \Rz^N$, $\langle x + a \rangle \geq \langle x \rangle / (1 + \|a\|)$.
        \item For all $x, a \in \Rz^N$ with $\|a\| \leq 1$,
        \begin{equation}\label{eq:symweight}
            \frac{1}{2}\bigl(\langle x + a \rangle^{-\ga} + \langle x - a \rangle^{-\ga}\bigr) \leq \bigl(1 + C\|a\|^2\bigr)\langle x \rangle^{-\ga}.
        \end{equation}
        \item Let $E$ be a finite--dimensional normed space and $G: \Rz^N \to E$ be bounded. If $|G|_{\ga - \al} < \infty$ for some $\al \in (0,1)$, then for all $x \neq y$ with $\|x - y\| > \frac{1}{2}\min(\langle x \rangle, \langle y \rangle)$,
        \begin{equation}\label{eq:farpairs}
            |G(x) - G(y)| \leq 2^{2 + \al}|G|_{\ga - \al}\,\|x - y\|^{\al}\,\omega_{\ga}(x,y).
        \end{equation}
        In particular, 
        \begin{equation}
            [G]_{\al,\ga} \leq 2^{2+\al}|G|_{\ga - \al} + \sup_{\|x - y\| \leq \frac{1}{2}\min(\langle x \rangle, \langle y \rangle)}\frac{|G(x) - G(y)|}{\|x-y\|^{\al}\omega_{\ga}(x,y)}.
        \end{equation}
    \end{itemize}
\end{lem}
\begin{proof}
    The first point is Peetre's inequality \cite[Section 7.1]{Hormander1990}. For the second point, let $\phi \coloneqq \langle \cdot \rangle^{-\ga}$, which is smooth with $|\nabla^2 \phi(x)| \leq C(\ga)\langle x \rangle^{-\ga - 2}$ for all $x \in \Rz^N$. By the first point, $\langle x \pm r a \rangle \geq \langle x \rangle/2$ for $r \in [0,1]$, hence Taylor's formula gives
    \begin{equation}
        \begin{aligned}
            &\phi(x + a) + \phi(x - a) - 2\phi(x)\\
            &\qquad = \int_0^1 (1 - r)\Bigl(\nabla^2 \phi(x + ra)(a,a) + \nabla^2\phi(x - ra)(a,a)\Bigr) \der r\\
            &\qquad \leq C(\ga) 2^{\ga + 2}\|a\|^2\langle x \rangle^{-\ga - 2}.
        \end{aligned}
    \end{equation}
    For the third point, w.l.o.g let $\langle x \rangle \leq \langle y \rangle$. Then 
    \begin{equation}
      \begin{aligned}
         &|G(x) - G(y)| \leq 2|G|_{\ga - \al}\langle x \rangle^{-\ga + \al} \leq 2|G|_{\ga - \al}\langle x \rangle^{-\ga}(2\|x - y\|)^{\al}\\ &\qquad \leq 2^{2 + \al}|G|_{\ga - \al}\|x - y\|^{\al}\omega_{\ga}(x,y).
      \end{aligned}
    \end{equation}
\end{proof}

\begin{lem}
   Let $\nu \in \mathfrak E$ with $\div \nu = 0$. Then, $\nu$ has mean zero. Moreover, for $\nu \in \mathfrak E$ and $q \in [0, 1 - \ep)$,
   \begin{equation}\label{eq:Etomoments}
      |\nu|(\Rz^N) + \Mc_q(\nu) \leq C(N, q, \ep)\|\nu\|_{\mathcal E}, \qquad \|\nabla \nu\|_{L^1} \leq C(N,\ep)\|\nu\|_{\mathcal E}.
   \end{equation}
\end{lem}
\begin{proof}
   To prove the first point, we note that $\hat\nu$ is continuous with $\xi \cdot \hat \nu(\xi) = 0$ for $\xi \neq 0$, hence, $\hat \nu(0) = 0$. For the second point, we use
   \begin{equation}
      |\nu|(\Rz^N) + \Mc_q(\nu) \leq \int_{\Rz^N} \bigl(1 + \|x\|^q\bigr)\langle x \rangle^{-N-1+\ep} \der x \; |\nu|_{N+1-\ep}.
   \end{equation}
\end{proof}

\begin{lem}[Weighted flow estimates]\label{le:flowweights}
    Let $\al \in (0,1)$ and $\ep > 0$. Let $h \in \mathcal E$ and $t > 0$ with $t\|h\|_{\mathcal E} \leq 1$. There exists $C = C(N, \al, \ep) > 0$ such that for all $x, y \in \Rz^N$:
    \begin{itemize}
        \item $\rh_t^h$ preserves the Lebesgue measure, $\|\rh_t^h(x) - x\| \leq t\|h\|_{\mathcal E}$, and $\langle \rh_t^h(x) \rangle \geq \langle x \rangle/(1 + t\|h\|_{\mathcal E})$,
        \item $\|D\rh_t^h(x)\| \leq e$, such that
        \begin{equation}
            \|D\rh_t^h(x) - I\| \leq C t\|h\|_{\mathcal E}\langle x \rangle^{-N-2+\ep}, \qquad \|D^2\rh_t^h(x)\| \leq C t\|h\|_{\mathcal E}\langle x \rangle^{-N-3+\ep},
        \end{equation}
        \item for $\|x - y\| \leq \frac{1}{2}\min(\langle x \rangle, \langle y \rangle)$, then $\|D\rh_t^h(x) - D\rh_t^h(y)\| \leq Ct\|h\|_{\mathcal E}\|x - y\|\langle x \rangle^{-N-3+\ep}$ and
        \begin{equation}\label{eq:D2rhoholder}
            \|D^2 \rh_t^h(x) - D^2 \rh_t^h(y)\| \leq C t\|h\|_{\mathcal E}\,\|x - y\|^{\al}\,\omega_{N + 3 + \al - \ep}(x,y).
        \end{equation}
    \end{itemize}
\end{lem}
\begin{proof}
    The fact that $\rh^h_r$ preserves the Lebesgue measure is classical, see \cite[Chapter 4]{MajdaBertozzi2002}. A simple Taylor--expansion together with point one of Lemma~\ref{le:weights} proves the first point of the result. 
    
    To simplify notation, write $\rh_r \coloneqq \rh^h_r$ for $r \in [0,t]$. First, Grönwall's inequality gives $$\|D\rh_r\| \leq e^{r\|\nabla h\|_{\infty}} \leq e.$$
    Moreover, the first point yields $\langle \rh_r(x) \rangle \geq \langle x \rangle/2$, such that $$|\nabla^k h(\rh_r(x))| \leq \|h\|_{\mathcal E}\langle \rh_r(x) \rangle^{-N-1-k+\ep} \leq 2^{N+3}\|h\|_{\mathcal E}\langle x \rangle^{-N-1-k+\ep}$$ for $k \in \{1,2\}$. 
    Hence, $\|D\rh_t(x) - I\| \leq \int_0^t |\nabla h(\rh_r(x))|e \der r \leq Ct\|h\|_{\mathcal E}\langle x \rangle^{-N-2+\ep}$ and $\|D^2\rh_t(x)\| \leq e^{t\|\nabla h\|_{\infty}}\int_0^t |\nabla^2h(\rh_r(x))|e^2 \der r \leq Ct\|h\|_{\mathcal E}\langle x \rangle^{-N-3+\ep}$, which is the second point.

    For the third point, let $\|x - y\| \leq \frac{1}{2}\min(\langle x \rangle, \langle y \rangle)$, so that $\langle z \rangle \geq \langle x \rangle/2$ for all $z \in [x,y]$ and $\langle x \rangle \leq 3\langle y \rangle$. 
    The estimate on $D\rh(x) - D\rh(y)$ follows from the bound on $D^2\rh$. 
    
    For \eqref{eq:D2rhoholder}, set $\ga \coloneqq N + 3 + \al - \ep$ and $\Delta(r): r \mapsto \|D^2\rh_r(x) - D^2\rh_r(y)\|\,\|x-y\|^{-\al}\omega_{\ga}(x,y)^{-1}$. Differentiating in $r$ and using $\|\rh_r(x) - \rh_r(y)\| \leq e\|x - y\|$, $\omega_{\ga}(\rh_r(x), \rh_r(y)) \leq 2^{\ga}\omega_{\ga}(x,y)$ and $\omega_{\ga}(x,y) \geq \frac{1}{2}\langle x \rangle^{-\ga}$, we estimate the four resulting differences:
    \begin{equation*}
      \begin{aligned}
         &|\nabla^2h(\rh_r(x)) - \nabla^2 h(\rh_r(y))| \leq [\nabla^2 h]_{\al,\ga}e^{\al}2^{\ga}\|x - y\|^{\al}\omega_{\ga}(x,y),\\
         &|\nabla h(\rh_r(y))|\,\|D^2\rh_r(x) - D^2\rh_r(y)\| \leq \|h\|_{\mathcal E}\Delta(r)\|x-y\|^{\al}\omega_{\ga}(x,y),\\
         &|\nabla^2h(\rh_r(y))|\,\|D\rh_r(x) - D\rh_r(y)\| \leq C\|h\|_{\mathcal E}^2 r\|x-y\|\langle x \rangle^{-2N-6+2\ep}\\
         &\qquad \leq C\|h\|_{\mathcal E}\|x - y\|^{\al}\langle x \rangle^{-\ga},\\
         &|\nabla h(\rh_r(x)) - \nabla h(\rh_r(y))|\,\|D^2\rh_r(x)\|\\
         &\qquad \leq C\|h\|_{\mathcal E}\langle x \rangle^{-N-3+\ep}\|x - y\|\cdot C r\|h\|_{\mathcal E}\langle x \rangle^{-N-3+\ep},
      \end{aligned}
    \end{equation*}
    where we used $r\|h\|_{\mathcal E} \leq 1$ and $\|x - y\|^{1-\al} \leq \langle x \rangle^{1 - \al}$. Hence, $\Delta(t) \leq \int_0^t C\|h\|_{\mathcal E} + \|h\|_{\mathcal E}\Delta(r) \der r$, and Grönwall's inequality yields $\Delta(t) \leq Ct\|h\|_{\mathcal E}e^{t\|h\|_{\mathcal E}} \leq Cet\|h\|_{\mathcal E}$, and the result follows.
\end{proof}

\section{Evolution operators on $\Wz_1$ and the incompressible Navier--Stokes equations}\label{ch:parstructwas}
In this section, we transfer the differential structure of \citet{Schichl2025} to the normed Wasserstein space $\Wz_1$ and apply it to the incompressible Navier--Stokes equations. In Section~\ref{sbs:ops}, we define the analogues on $\Wz_1(M)$ of the averaged flow dispersion operators and of the expectational flow operators of \citet{Schichl2025}, together with the source operators, which encode a general source term, and the dievergence operators, which encode the pressure term of the Navier--Stokes equations. Section~\ref{sbs:divfree} introduces divergence--free tangent vectors and the incompressible flow vectors, for which the expectational flow operators act on densities by composition with a measure--preserving flow. Section~\ref{sbsec:flowestims} collects the estimates on all these operators in $L^{\infty}$, Hölder, $L^1$, moment and weighted Hölder norms, which are the technical core of this work. Finally, in Section~\ref{sbs:NS}, we build the linear scheme (LIFAS) and the non--linear scheme (IFAS) for the incompressible Navier--Stokes equations and prove our main results: the convergence of the (LIFAS) to any given classical solution (Theorem~\ref{th:LIFAS}), the resulting a priori estimates for classical solutions, and the convergence of the (IFAS) to the unique classical solution on an explicit time interval (Theorem~\ref{th:IFAS}), with its consequences for the lifespan of solutions and the inviscid limit (Corollary~\ref{co:lifespan}).

\subsection{Evolution operators on $\Wz_1$}\label{sbs:ops}

Let $(M,g)$ be a complete Riemannian manifold. We define analogues of the flow operators in \citet{Schichl2025} on $\Wz_1(M)$. Recall that a vector--field probability is an element $V \in \Wc_p(\Xc^{2,\infty}(M))$ for some $p>2$. The definition of the Average vector--field $\ol{V}$ can be found in \citet{Schichl2025}.
\begin{defi}[Averaged flow operators]\label{de:avflowop}
Let $p>2$, $(M,g)$ be a Riemannian manifold and $V \in \Wc_p(\Xc^{2,\infty}(M))$. First, we define the \emph{averaged flow dispersion operators associated to $V$} by
\begin{equation}
\begin{aligned}
f^{V}_t: \Wz_1(M) &\to \Wz_1(M) \\
\nu &\mapsto \int_{\Xc(M)} (\rh^{X-\ol{V}}_{\sqrt t})_*\nu \der V(X)
\end{aligned}
\end{equation}
for all $t\geq 0$. Furthermore, define the \emph{expectational flow operators}
\begin{equation}
\begin{aligned}
e^{V}_t: \Wz_1(M) &\to \Wz_1(M) \\
\nu &\mapsto (\rh^{\ol{V}}_t)_* \nu
\end{aligned}
\end{equation}
for all  $t \geq 0$.
\end{defi}
Let us show that these operators are well--defined and establish some generalbasic properties.
\begin{lem}\label{le:flowopwd}
    Let $p>2$, $(M,g)$ be a complete Riemannian manifold. Take a vector--field probability $V \in \Wc_p(\Xc^{2,\infty}(M))$ and $t,s \geq 0$. Then, $f^V_t$ and $e^V_t$ are well--defined linear operators, which map positive measures onto positive measures, and for all $\nu \in \Wz_1(M)$,
    \begin{itemize}
        \item $\bigl|f^V_t(\nu)\bigr|(M) \leq |\nu|(M)$ and $\bigl|e^V_t(\nu)\bigr|(M) = |\nu|(M)$,
        \item $\|f^V_t(\nu) - f^V_s(\nu)\|_{\Wz_1} \leq 2\bigl|\sqrt{t} - \sqrt{s}\bigr|\Mc_1(V)|\nu|(M)$,
        \item $\|e^V_t(\nu) - e^V_s(\nu)\|_{\Wz_1} \leq |t-s|\Mc_1(V)|\nu|(M)$.
    \end{itemize}
\end{lem}
\begin{proof}
    By \cite[Lemma 3.9]{Schichl2025}, $\ol{V} \in \Xc^{2,\infty}(M)$ with $\|\ol{V}\|_{\infty} \leq \Mc_1(V)$. In particular, for all $X \in \Xc^{2,\infty}(M)$ the flows $\rh^{\ol{V}}_t$ and $\rh^{X - \ol V}_{\sqrt t}$ are globally defined bi--Lipschitz homeomorphisms of $M$. Hence, $e^V_t(\nu)$ is a well--defined finite signed Borel--measure, satisfying $\bigl(e^V_t(\nu)\bigr)^{\pm} = (\rh^{\ol V}_t)_*\nu^{\pm}$ and, thus, $$\bigl|e^V_t(\nu)\bigr|(M) = |\nu|(M).$$
    Let $\nu \in \Wz_1(M)$ and $f \in \Cc_b(M)$. Then,
    \begin{equation*}
        \begin{aligned}
            M \times \Xc^{2,\infty}(M) &\to \Rz\\
            (x, X) &\mapsto f\bigl(\rh^{X - \ol V}_{\sqrt t}(x)\bigr)
        \end{aligned}
    \end{equation*}
    is continuous and bounded by $\|f\|_{\infty}$. Hence, as in \cite[Remark 3.19]{Schichl2025}, Fubini--Tonelli's theorem yields that
    \begin{equation}\label{eq:fdef}
        \begin{aligned}
            \Cc_b(M) &\to \Rz\\
            f &\mapsto \int_{\Xc^{2,\infty}(M)} \int_M f \circ \rh^{X - \ol V}_{\sqrt t} \der \nu \der V(X)
        \end{aligned}
    \end{equation}
    is a well--defined linear functional on $\Cc_b(M)$, with operator norm bounded by $|\nu|(M)$, which is positive if $\nu$ is. Now, we readily conclude using the Riesz representation theorem.

    Now, let $f \in L^{\infty}(M) \cap \Lip(M)$ with $\max(\|f\|_{\infty}, \|f\|_{\Lip}) \leq 1$ and $X \in \Xc^{2,\infty}(M)$. Then, for all $x \in M$,
    \begin{equation}
        \begin{aligned}
            &\Bigl| f\bigl(\rh^{X - \ol V}_{\sqrt t}(x)\bigr) - f\bigl(\rh^{X - \ol V}_{\sqrt s}(x)\bigr)\Bigr| \leq d\Bigl(\rh^{X - \ol V}_{\sqrt t}(x), \rh^{X - \ol V}_{\sqrt s}(x)\Bigr)\\
            &\qquad \leq \bigl|\sqrt t - \sqrt s\bigr| \bigl(\|X\|_{\infty} + \|\ol V\|_{\infty}\bigr).
        \end{aligned}
    \end{equation}
    Integrating with respect to $|\nu|$ and $V$, and using $\|X\|_{\infty} \leq \|X\|_{W^{2,\infty}}$ and $\|\ol V\|_{\infty} \leq \Mc_1(V)$, we obtain
    \begin{equation}
        \begin{aligned}
                &\biggl| \int_M f \der f^V_t(\nu) - \int_M f \der f^V_s(\nu) \biggr|\\
                &\qquad \leq \bigl|\sqrt t - \sqrt s\bigr| \bigl(\Mc_1(V) + \|\ol V\|_{\infty}\bigr)|\nu|(M) \leq 2\bigl|\sqrt t - \sqrt s\bigr|\Mc_1(V)|\nu|(M).
        \end{aligned}
    \end{equation}
    Taking the supremum over $f$ yields the second point. The third point follows in the same way from $d\bigl(\rh^{\ol V}_t(x), \rh^{\ol V}_s(x)\bigr) \leq |t-s|\|\ol V\|_{\infty}$, compare \cite[Lemma 3.9 and Lemma 3.31]{Schichl2025}.
\end{proof}
We further introduce a class of operators that map $\Wz_1(M)$ into $\Wz_1(M)$, but in general do not map probability measures onto probability measures. In particular, they have no analogue in \cite{Schichl2025}. Unlike the averaged flow dispersion operators and the expectational flow operators, these operators do not encode a transport of mass, but rather a source term. 
\begin{defi}[Source operators]
    Let $T>0$, $N \in \Nz$ and $d \in \Cc^1\bigl([0,T) \times M, \Rz^N\bigr)$. We define the \emph{source operators} associated to $d$ by
    \begin{equation}
        \begin{aligned}
            s^{d,l}_{t,i}: \Wz_1(M) &\to \Wz_1(M) \\
            \nu &\mapsto \nu + \int_l^{t+l} \nabla d(s)_i \der s,
        \end{aligned}
    \end{equation}
    for all $t,l >0$ such that $t + l < T$ and $i \in \llbracket 0, N-1 \rrbracket$.
\end{defi}
Introducing the source operators enables to mimic the procedure introduced in \citet{Schichl2025} for non--linear parabolic PDEs with source term. We shall call the three families of operators introduced above \emph{general evolution operators}. 

Finally, we define one further class of operators, the dievergence operator, which we shall use to build approximations of solutions of the Navier--Stokes equations in particular. 
\begin{defi}[Dievergence operator]
    Let $N \in \Nz$ and $\theta \in [0,1]$. We define the \emph{dievergence operator} by
    \begin{equation}
        \begin{aligned}
            s_{\theta}: \Sc\bigl(\Rz^N, \Rz^N\bigr) &\to \Sc\bigl(\Rz^N, \Rz^N\bigr) \\
            \nu &\mapsto \nu - \theta\nabla \Delta^{-1} \div \nu.
        \end{aligned}
    \end{equation}
    We shall simply write $s \coloneqq s_1$.
\end{defi}
The dievergence operator is a priori defined on the space of Schwarz--functions. However, we shall see that it can be extended to a bounded linear operator on $L^p$ for $p \in (1,\infty)$.

\subsection{Divergence--free tangent vectors}\label{sbs:divfree}

To obtain stronger estimates on the flow operators than those of Lemma~\ref{le:flowopwd}, we shall restrict the class of vector--field probabilities to those that are divergence--free. This is in particular motivated by the fact that the incompressible Navier--Stokes equations are a system of parabolic PDEs with source term, whose solutions are divergence--free vector--fields.

\begin{defi}[Divergence--free vectors]
    Let $v \in \Wc_p(\Xc^{2,\infty}(M))$ such that for all $X \in \supp(v)$, $X$ is divergence-free. Then $v$ is called a divergence--free vector. A probability vector--field
    $V: \Wz_1(M) \to \Wc_p(\Xc^{2,\infty}(M))$ is \emph{divergence--free} if for all $\mu \in \Wz_1(M)$, $V[\mu]$ is divergence--free.
\end{defi}
First, we show that the average vector--field of a divergence--free vector is itself divergence--free.
\begin{lem}\label{le:avdivfree}
    Let $V \in \Wc_p(\Xc^{2,\infty}(M))$ by divergence--free. Then, $\ol{V}$ is a divergence--free vector--field.
\end{lem}
\begin{proof}
    By \cite[Lemma 3.9]{Schichl2025}, $\ol{V} \in \Xc^{2,\infty}(M)$. Let $x \in M$ and fix local coordinates $(x_0, \dots, x_{n-1})$ around $x$ in which the volume form reads $\sqrt{\det g}\der x$. 
    As in the proof of \cite[Lemma 3.9, (40)]{Schichl2025}, we obtain for all $i \in \llbracket 0, n-1 \rrbracket$,
    \begin{equation}
        \begin{aligned}
            \frac{\der \ol V(x)}{\der x_i} = \int_{\Xc^{2,\infty}(M)} \pi_x\biggl(\frac{\der X}{\der x_i}\biggr) \der V(X).
        \end{aligned}
    \end{equation}
    Therefore,
    \begin{equation}
        \begin{aligned}
            \div \ol V(x) &= \frac{1}{\sqrt{\det g}} \sum_{i=0}^{n-1} \frac{\der}{\der x_i}\Bigl(\sqrt{\det g}\, \ol V_i\Bigr)(x)\\
            &= \int_{\Xc^{2,\infty}(M)} \frac{1}{\sqrt{\det g}} \sum_{i=0}^{n-1} \frac{\der}{\der x_i}\Bigl(\sqrt{\det g}\, X_i\Bigr)(x) \der V(X)\\
            &= \int_{\Xc^{2,\infty}(M)} \div X(x) \der V(X) = 0,
        \end{aligned}
    \end{equation}
    since $V$ is divergence--free.
\end{proof}

The particular vector--field probabilities we shall use the flow operators introduced in Section~\ref{sbs:ops} with are the so--called Incompressible Flow Vectors.
\begin{defi}[Incompressible Flow vectors]
    Let $h \in \Xc^{2,\infty}(\Rz^N)$ such that $h$ is divergence--free. Let $(e_i)_{i \in \llbracket 0, N-1\rrbracket}$ be the canonical base of $\Rz^N$ and denote for all $i \in \llbracket 0, N-1\rrbracket$ by $X_i$ the constant vector--field equal to $e_i$. Define,
    \begin{equation}
        \Phi(h) \coloneqq \frac{1}{2N} \sum_{i=0}^{N-1} \bigl(\delta_{\eta X_i - h} + \delta_{- \eta X_i -  h}\bigr).
    \end{equation}
    In particular, if $T >0$ and $(g,d)$ is a solution of \eqref{eq:NS} on $[0,T)$, we shall write for all $t \in [0,T)$,
    \begin{equation}
        \Phi(g,t) \coloneqq \Phi\bigl(g(t)\bigr).
    \end{equation}
\end{defi}
Note that the Incompressible Flow vectors do indeed satisfy the assumptions of Definition~\ref{de:avflowop}.
\begin{co}
    For all $h \in \Xc^{2,\infty}(\Rz^N)$, $\Phi(h)$ has bounded support.
\end{co}
\begin{proof}
This follows directly from Proposition~\ref{pr:regclasssol}.
\end{proof}
The following result is an immediate consequence of the definition of the average vector--field of a vector--field probability.
\begin{co}\label{co:avPhi}
    For all $h \in \Xc^{2,\infty}(\Rz^N)$, $\ol{\Phi(h)} = -h$.
\end{co}

\subsection{Evolution operators: properties and estimates}\label{sbsec:flowestims}

This section is devoted to a list of technical estimates for the general evolution operators introduced in Section~\ref{sbs:ops} with divergence--free tangent vectors or---even more specifically---with Incompressible Flow vectors. These estimates are crucial to prove the main results of this work.

\subsubsection{Properties of the general evolution operators with divergence--free tangent vectors}

Unlike in \citet{Schichl2025}, in this setting, only the expectational flow operators carry
a non--linearity, while the averaged flow operators just describe a simple isotropic diffusion. Hence, the crucial estimates are those for the expectational flow operators, which we give in the following lemma.
The statements in Lemma~\ref{le:eflowop} heavily rely on the fact that the vector--fields $g(t)$  are divergence--free for all $t \in [0,T)$.

\begin{lem}[Properties of the expectational flow operators]\label{le:eflowop}
Let $N \in \Nz$, $p > 2$ and let $v \in \Wc_p(\Xc^{2,\infty}(\Rz^N))$ be a divergence--free vector, so that $\ol v \in \Xc^{2,\infty}(\Rz^N)$ is divergence--free by Lemma~\ref{le:avdivfree}. Then, for all $t,s >0$ and all $\mu \in \Wz_1(\Rz^N)$, the following holds:
\begin{itemize}
    \item $e^{v}_t(\mu) \in \Wz_1(\Rz^N)$, and $\bigl|e^{v}_t(\mu)\bigr|(\Rz^N) = |\mu|(\Rz^N)$.
    \item $\|e^{v}_s(\mu) - e^{v}_t(\mu)\|_{\Wz_1} \leq |t-s| \|\ol v\|_{\infty}|\mu|(\Rz^N)$.
    \item if $\mu \in L^{\infty}(\Rz^N)$, then $e^{v}_t(\mu) \in L^{\infty}(\Rz^N)$, $\|e^{v}_t(\mu)\|_{\infty} = \|\mu\|_{\infty}$, and for almost every $x \in \Rz^N$,
    \begin{equation}
        e^{v}_t(\mu)(x) = \mu\bigl(\rho^{-\ol v}_t(x)\bigr).
    \end{equation}
    \item if there exists $\al \in (0,1]$ such that $\mu \in \Cc^{\al}(\Rz^N)$, then $e^{v}_t(\mu) \in \Cc^{\al}(\Rz^N)$ and $\bigl\{e^{v}_t(\mu)\bigr\}_{\Cc^{\al}} \leq \{\mu\}_{\Cc^{\al}}e^{t\|\ol v\|_{\Lip} \al}$.
    \item let $q \geq 1$, and suppose $\mu$ has finite $q$--moment, then $e^{v}_t(\mu)$ has finite $q$--moment.
    \item let $q \in (0,1)$, and suppose $\mu$ has finite $q$--moment and is bounded, then $e^{v}_t(\mu)$ has finite $q$--moment and
        \begin{equation}
            \begin{aligned}
                &\int_{\Rz^N} \|x\|^q \der \bigl|e^{v}_t\bigr|(\mu)(x) \leq \Mc_q(\mu) + \|\ol v\|_{\infty} \bigl(|\mu|(\Rz^N) + C(N,q)\|\mu\|_{\infty}\bigr)t,
            \end{aligned}
        \end{equation}
        where $C(N,q) = \int_{\|x\| \leq 1} \frac{1}{\|x\|^{1-q}} \der x < \infty$.
    \item suppose $\mu \in W^{1,\infty}(\Rz^N)$ and $ \nabla \mu \in L^1(\Rz^N)$.  Then, $e^{v}_t(\mu) \in W^{1,\infty}(\Rz^N)$, $ \nabla e^{v}_t(\mu) \in L^1(\Rz^N)$ and
    \begin{equation}
        \|\nabla e^{v}_t(\mu)\|_{L^1} \leq \|\nabla \mu \|_{L^1}e^{t\|\ol v\|_{\Lip}}.
    \end{equation}
\end{itemize}
In particular, since $\ol{\Phi(h)} = -h$ by Corollary~\ref{co:avPhi}, the above holds for the incompressible flow vectors $v = \Phi(h)$, $h \in \Xc^{2,\infty}(\Rz^N)$ divergence--free, with $\rh^{-\ol v}_t = \rh^{h}_t$ and $\|\ol v\|_{\infty} = \|h\|_{\infty}$, $\|\ol v\|_{\Lip} = \|h\|_{\Lip}$.

Finally, assume $N \geq 2$. Let $\al, \ep \in (0,1)$. There exists $K_0 = K_0(N, \al, \ep) \geq 1$ such that for all $h \in \mathcal E$, $\mu \in \mathfrak E$ and all $t > 0$ such that $t\|h\|_{\mathcal E} \leq 1$:
    \begin{itemize}
        \item $e^{\Ph(h)}_t(\mu) \in \mathfrak E$, $\|e^{\Ph(h)}_t(\mu)\|_{\mathcal E} \leq (1 + K_0 t\|h\|_{\mathcal E})\|\mu\|_{\mathcal E}$ and $\int_{\Rz^N} e^{\Ph(h)}_t(\mu) \der x = \int_{\Rz^N} \mu \der x$;
        \item if we further assume $\div \mu = 0$, then $F \coloneqq \div\bigl(e^{\Ph(h)}_t(\mu)\bigr)$ satisfies $\|F\|_{\mathcal F} \leq K_0 t\|h\|_{\mathcal E}\|\mu\|_{\mathcal E}$, and $\int_{\Rz^N} F \der x = 0$, $\int_{\Rz^N} x_j F \der x = -\int_{\Rz^N} \mu_j \der x$ for all $j \in \llbracket 0, N-1 \rrbracket$.
    \end{itemize}
\end{lem}
\begin{proof}

The first point follows immediately since $e^{v}_t(\mu)$ is a push--forward of $\mu$ under the homeomorphism $\rh^{\ol v}_t$.

For the second point, note that for all $f \in L^{\infty}(\Rz^N) \cap \Lip(\Rz^N)$, such that $\max(\|f\|_{\infty}, \|f\|_{\Lip}) \leq 1$,
\begin{equation}
    \begin{aligned}
        &\biggl|\int_{\Rz^N} f \der e^{v}_t(\mu) - \int_{\Rz^N} f \der e^{v}_s(\mu)\biggr| = \biggl|\int_{\Rz^N} f \circ \rho^{\ol v}_t - f \circ \rho^{\ol v}_s \der \mu\biggr| \\
        &\qquad \leq \int_{\Rz^N}\bigl|\rho^{\ol v}_t - \rho^{\ol v}_s\bigr| \der |\mu| \leq \int_{\Rz^N} \int_s^t \bigl|\ol v\bigl(\rho^{\ol v}_r(x)\bigr)\bigr| \der r \der |\mu|(x)\\
        &\qquad \leq |t-s| \|\ol v\|_{\infty}|\mu|(\Rz^N).
    \end{aligned}
\end{equation}

For the third point, suppose $\mu \in L^{\infty}(\Rz^N)$. Since $\ol v$ is a divergence--free Lipschitz vector--field, its flow preserves the Lebesgue measure. Let $\ep >0$ and let $\ph_{\ep}$ be a standard mollifier supported in $B(0,\ep)$. Then, for all $x \in \Rz^N$,
\begin{equation}
   \begin{aligned}
        &\int_{\Rz^N} \ph_{\ep}(x-y) \der e^{v}_t(\mu)(y) = \int_{\Rz^N} \ph_{\ep}\bigl(x- \rho^{\ol v}_t(y)\bigr) \mu(y) \der y\\
        &\qquad = \int_{\Rz^N} \ph_{\ep}(x- z) \mu\bigl(\rho^{-\ol v}_t(z)\bigr) \der z,
   \end{aligned}
\end{equation}
by the change of variables $z = \rho^{\ol v}_t(y)$. Since $\ep>0$ and $x \in \Rz^N$ were chosen arbitrarily, it follows that $e^{v}_t(\mu)$ has the density $\mu \circ \rho^{-\ol v}_t \in L^{\infty}(\Rz^N)$, whose $L^{\infty}$--norm equals $\|\mu\|_{\infty}$ by measure preservation.

For the fourth point, note that for all $x,y \in \Rz^N$, by Grönwall's inequality,
\begin{equation}
    \biggl\|\rho^{-\ol v}_t(x) - \rho^{-\ol v}_t(y)\biggr\| \leq e^{t\|\ol v\|_{\Lip}}\|x-y\|.
\end{equation}
Therefore,
\begin{equation}
    \begin{aligned}
        &\biggl|e^{v}_t(\mu)(x) - e^{v}_t(\mu)(y)\biggr| = \biggl|\mu\bigl(\rho^{-\ol v}_t(x)\bigr) - \mu\bigl(\rho^{-\ol v}_t(y)\bigr)\biggr|\\
        &\qquad \leq \{\mu\}_{\Cc^{\al}} \bigl\|\rho^{-\ol v}_t(x) - \rho^{-\ol v}_t(y)\bigr\|^{\al} \leq \{\mu\}_{\Cc^{\al}} e^{t\al\|\ol v\|_{\Lip}}\|x-y\|^{\al}.
    \end{aligned}
\end{equation}

Let us prove the moment estimates. Let $q \geq1$ and suppose $\mu$ has finite $q$--moment. Then,
\begin{equation}
    \begin{aligned}
        &\biggl(\int_{\Rz^N} \|x\|^q \der \bigl|e^{v}_t(\mu)\bigr|(x) \biggr)^{\frac{1}{q}}= \biggl(\int_{\Rz^N} \bigl\|\rho^{\ol v}_t(x)\bigr\|^q \der |\mu|(x)\biggr)^{\frac{1}{q}}\\
        &\qquad \leq \Mc_q(\mu)^{\frac{1}{q}} + \biggl(\int_{\Rz^N} \biggl\|\int_0^t \ol v\bigl(\rho^{\ol v}_s(x)\bigr)\der s\biggr\|^q \der |\mu|(x)\biggr)^{\frac{1}{q}}\\
        &\qquad \leq \Mc_q(\mu)^{\frac{1}{q}} + t\|\ol v\|_{\infty}\bigl|\mu\bigr|(\Rz^N)^{\frac{1}{q}}.
    \end{aligned}
\end{equation}
If $q \in (0,1)$, and $\mu \in L^{\infty}(\Rz^N)$, then, since $y \mapsto \|y\|^q$ is $q$--Hölder and $\bigl|\|x+y\|^q - \|x\|^q\bigr| \leq q\|y\|\,\|x\|^{q-1}$ for $\|y\| \leq \|x\|/2$, we obtain
\begin{equation}
    \begin{aligned}
        &\int_{\Rz^N} \|x\|^q \der \bigl|e^{v}_t(\mu)\bigr|(x) = \int_{\Rz^N} \biggl\|x + \int_0^t \ol v( \rho^{\ol v}_s(x))\der s\biggr\|^q \der |\mu|(x)\\
        &\qquad \leq \Mc_q(\mu) + t \int_{\Rz^N} \frac{\|\ol v\|_{\infty}}{\|x\|^{1-q}} \der |\mu|(x)\\
        &\qquad \leq \Mc_q(\mu) + \|\ol v\|_{\infty} \bigl(|\mu|\bigl(\Rz^N\bigr) + C(N,q)\|\mu\|_{\infty}\bigr)t,
    \end{aligned}
\end{equation}
where $C(N,q) = \int_{\|x\|\leq 1} \frac{1}{\|x\|^{1-q}} \der x < \infty$; here we used $|\mu| \leq \|\mu\|_{\infty}\der x$ on $B(0,1)$ and $\|x\|^{q-1} \leq 1$ outside. This concludes the proof of the moment estimates.

For the seventh point, let $\mu \in W^{1,\infty}(\Rz^N)$ and $\nabla \mu \in L^1(\Rz^N)$. Then, by the third point, for almost every $x \in \Rz^N$,
\begin{equation}
    \nabla e^{v}_t(\mu)(x) = \nabla \mu\bigl(\rho^{-\ol v}_t(x)\bigr)\circ d \rho^{-\ol v}_t(x).
\end{equation}
Hence, again by Grönwall,
\begin{equation}
    \|\nabla e^{v}_t(\mu)(x)\| \leq \|\nabla \mu\bigl(\rho^{-\ol v}_t(x)\bigr)\|e^{t\|\ol v\|_{\Lip}}.
\end{equation}
Since $\ol v$ is divergence--free, $\rho^{-\ol v}_t$ preserves the Lebesgue measure, and it follows $$\|\nabla e^{v}_t(\mu)\|_{L^1} \leq \|\nabla \mu \|_{L^1}e^{t\|\ol v\|_{\Lip}}.$$

Let us prove the last two points. To simplify notation, let $\rh \coloneqq \rh^h_t$ and
 $\la \coloneqq t\|h\|_{\mathcal E} \leq 1$. Let $\be \leq N + 5$. By Lemma~\ref{le:flowweights}, $\langle x \rangle^{\be} \leq (1 + \la)^{\be}\langle \rh(x) \rangle^{\be} \leq (1 + C\la)\langle \rh(x)\rangle^{\be}$ for some $C$ depending only on $N$, such that $$|\mu\circ \rh|_{N+1-\ep} \leq (1 + C\la)|\mu|_{N+1-\ep}.$$ 
Therefore, the second point of Lemma~\ref{le:flowweights} yields,
    \begin{equation}\label{eq:auxflowe1}
        \begin{aligned}
            &|\nabla(\mu\circ\rh)|_{N+2-\ep} \leq (1 + C\la)^2|\nabla\mu|_{N+2-\ep},\\
            &|\nabla^2(\mu\circ\rh)|_{N+3-\ep} \leq (1 + C\la)^3|\nabla^2\mu|_{N+3-\ep}\\
            &\qquad + \sup_x \langle x \rangle^{N+3-\ep}\,C\langle x \rangle^{-N-2+\ep}|\nabla \mu|_{N+2-\ep}\,C\la\langle x \rangle^{-N-3+\ep} \leq (1 + C\la)\|\mu\|_{\mathcal E}.
        \end{aligned}
    \end{equation}
    Now, set $\ga \coloneqq N + 3 + \al - \ep$. It remains to bound $[\nabla^2(\mu\circ\rh)]_{\al,\ga}$. For all $x \neq y$,
    \begin{equation}
        \begin{aligned}
            &\bigl|\nabla^2\mu(\rh(x))\bigl(D\rh(x), D\rh(x)\bigr) - \nabla^2\mu(\rh(y))\bigl(D\rh(y), D\rh(y)\bigr)\bigr|\\
            &\qquad \leq |\nabla^2\mu(\rh(y))|\,\bigl\|D\rh(x)\bigr\|\bigl\| D\rh(x) - D\rh(y) \bigr\|\\
            &\qquad\quad  + |\nabla^2\mu(\rh(y))|\,\bigl\|D\rh(y)\bigr\|\bigl\| D\rh(x) - D\rh(y) \bigr\|\\
            &\qquad\quad  + |\nabla^2\mu(\rh(x)) - \nabla^2\mu(\rh(y))|\,\|D\rh(x)\|^2.
        \end{aligned}
    \end{equation}
    The last term is bounded by 
    \begin{equation}
        \begin{aligned}
            &[\nabla^2\mu]_{\al,\ga}\|\rh(x) - \rh(y)\|^{\al}\omega_{\ga}(\rh(x),\rh(y))(1 + C\la)^2\\
            &\qquad \leq (1 + C\la)[\nabla^2\mu]_{\al,\ga}\|x - y\|^{\al}\omega_{\ga}(x,y), 
        \end{aligned}
    \end{equation} 
    since $\|\rh(x) - \rh(y)\| \leq e^{t\|\nabla h\|_{\infty}}\|x - y\| \leq (1 + C\la)\|x - y\|$ and $\omega_{\ga}(\rh(x),\rh(y)) \leq (1 + \la)^{\ga}\omega_{\ga}(x,y)$ by Lemma~\ref{le:flowweights}. 
        
    For the first two terms, recall that, by the third point of Lemma~\ref{le:weights} and the bound on $|\nabla^2(\mu\circ\rh)|_{N+3-\ep} = |\nabla^2(\mu \circ \rh)|_{\ga - \al}$ obtained in \eqref{eq:auxflowe1}, it suffices to consider pairs with $\|x - y\| \leq \frac{1}{2}\min(\langle x \rangle, \langle y \rangle)$.  
    Then, $\langle x \rangle \leq 3 \langle y \rangle$ and $\omega_{\ga}(x,y) \geq \frac{1}{2}\langle x \rangle^{-\ga}$. 
    Lemma~\ref{le:flowweights} gives $\|D\rh(x) - D\rh(y)\| \leq C\la\|x - y\|\langle x \rangle^{-N-3+\ep}$ and $|\nabla^2\mu(\rh(y))| \leq C|\nabla^2\mu|_{N+3-\ep}\langle x \rangle^{-N-3+\ep}$, such that both terms are bounded by $$C\la\|\mu\|_{\mathcal E}\|x - y\|^{\al}\langle x \rangle^{1-\al}\langle x \rangle^{-2N-6+2\ep} \leq C\la\|\mu\|_{\mathcal E}\|x-y\|^{\al}\omega_{\ga}(x,y),$$ where we use $\|x - y\|^{1-\al} \leq \langle x \rangle^{1-\al}$ and $2N + 5 + \al - 2\ep \geq \ga$. 
    
    The term $(\nabla\mu)\circ\rh\,D^2\rh$ is treated in the same way, using $$|\nabla\mu(\rh(x)) - \nabla\mu(\rh(y))| \leq C|\nabla^2\mu|_{N+3-\ep}\langle x \rangle^{-N-3+\ep}\|x - y\|,$$ the bounds on $D^2\rh$ and \eqref{eq:D2rhoholder}.
    
    Hence, in total $[\nabla^2(\mu\circ\rh)]_{\al,\ga} \leq (1 + C\la)\|\mu\|_{\mathcal E}$, which proves the norm estimate. The integral is preserved since $\rh$ preserves the Lebesgue measure.

    For the last point, since $\div \mu = 0$, $F = \div(\mu \circ \rh) = \operatorname{tr}\bigl((\nabla\mu)\circ\rh\,(D\rh - I)\bigr)$, and
    \begin{equation}
        \nabla F = \operatorname{tr}\bigl((\nabla^2\mu)\circ\rh[D\rh, \cdot]\,(D\rh - I)\bigr) + \operatorname{tr}\bigl((\nabla \mu)\circ\rh\,D^2\rh\bigr),
    \end{equation}
    where every term contains a factor $D\rh - I$ or $D^2\rh$, which are bounded by $C\la\langle x \rangle^{-N-2+\ep}$, resp. $C\la\langle x \rangle^{-N-3+\ep}$, and whose differences are controlled by the third point of Lemma~\ref{le:flowweights}. Exactly as above, this yields $|F|_{N+2-\ep} + |\nabla F|_{N+3-\ep} + [\nabla F]_{\al,\ga} \leq C\la\|\mu\|_{\mathcal E}$. 
    
    Finally, $\mu\circ\rh \in W^{1,1}$ gives $\int_{\Rz^N} F \der x = 0$, and $\int_{B(0,R)} x_j F \der x = \int_{\partial B(0,R)} x_j\, (\mu \circ \rh)\cdot n \der S - \int_{B(0,R)} (\mu\circ\rh)_j \der x$, where $CR\cdot R^{-N-1+\ep}R^{N-1} \to 0$ bounds the boundary term, such that $\int_{\Rz^N} x_j F \der x = -\int_{\Rz^N} \mu_j \der x$.
\end{proof}
While the Averaged flow dispersion operators just describe an isotropic diffusion in our setting, useful estimates are still somewhat difficult to obtain, as unlike in \citet{Schichl2025}, we have only very mild decay assumptions on the initial data and the source term. However, we work only with divergence--free vector--fields, 
which shall allow us to obtain even stronger estimates than for general AFAS with stronger assumptions on the initial data and zero source term.  
\begin{lem}[Properties of the averaged flow dispersion operators]\label{le:fflowop}
Let $N \in \Nz$, $p>2$, $t,s, \eta >0$ and $\mu \in \Wz_1(\Rz^N)$. Let $v \in \Wc_p(\Xc^{2,\infty}(\Rz^N))$ be a divergence--free vector. Then, the following holds:
\begin{itemize}
    \item $f^{v}_t(\mu) \in \Wz_1(\Rz^N)$, and $\bigl|f^{v}_t(\mu)\bigr|(\Rz^N) \leq |\mu|(\Rz^N)$.
    \item if $\mu \in L^{\infty}(\Rz^N)$, then $f^{v}_t(\mu) \in L^{\infty}(\Rz^N)$ and $\|f^{v}_t(\mu)\|_{\infty} \leq \|\mu\|_{\infty}$.
\end{itemize}
Furthermore, if $v = \Ph(h)$ for some divergence--free $h \in \Xc^{2,\infty}(\Rz^N)$, the following holds:
\begin{itemize}
    \item if $\mu \in W^{1,\infty}(\Rz^N)$ and $\nabla \mu \in L^1(\Rz^N)$, 
    \begin{equation}
        \|f^{\Ph(h)}_s(\mu) - f^{\Ph(h)}_t(\mu)\|_{\Wz_1} \leq \frac{3}{2}\eta^2\bigl\|\frac{\der} {\der x_i} \mu \bigr\|_{L^1(\Rz^N)} |t-s|.
    \end{equation}
    \item if there exists $\al \in (0,1]$ such that $\mu \in \Cc^{\al}(\Rz^N)$, then $f^{\Ph(h)}_t(\mu) \in \Cc^{\al}(\Rz^N)$ and $\bigl\{f^{\Ph(h)}_t(\mu)\bigr\}_{\Cc^{\al}} \leq \{\mu\}_{\Cc^{\al}}$.
    \item suppose $N>2$, let $q \in (0,1)$, $t<1$, and suppose $\mu$ has finite $q$--moment and is bounded, then $f^{\Ph(h)}_t(\mu)$ has finite $q$--moment and 
        \begin{equation}
            \begin{aligned}
                &\Mc_q\bigl(f^{\Ph(h)}_t(\mu)\bigr)\\
                &\qquad \leq \Mc_q(\mu) + 2^{1-q}\eta^2 q \bigl(D(N,q)(1 + \eta^q)\|\mu\|_{\infty} + |\mu|(\Rz^N)\bigr)t,
            \end{aligned}
        \end{equation} 
        where $D(N,q) >0$ is a constant depending only on $q$ and $N$.
    \item suppose $\mu \in W^{1,\infty}(\Rz^N)$ and $\nabla \mu \in L^1(\Rz^N)$. Then, $f^{\Ph(h)}_t(\mu) \in W^{1,\infty}(\Rz^N)$, $\nabla f^{\Ph(h)}_t(\mu) \in L^1(\Rz^N)$ and $\bigl\| \nabla f^{\Ph(h)}_t(\mu)\bigr\|_{L^1} \leq \bigl\| \nabla \mu\bigr\|_{L^1}$. 
\end{itemize}
Finally, suppose $N \geq 2$. Let $\al, \ep \in (0,1)$. There exists $K_0 = K_0(N, \al, \ep) \geq 1$ such that for all $h \in \mathcal E$, $\mu \in \mathfrak E$ and $t > 0$, such that $\eta\sqrt{t} \leq 1$, $f^{\Ph(h)}_t(\mu) \in \mathfrak E$, \begin{equation}\|f^{\Ph(h)}_t(\mu)\|_{\mathcal E} \leq (1 + K_0\eta^2 t)\|\mu\|_{\mathcal E},\end{equation} $\int_{\Rz^N} f^{\Ph(h)}_t(\mu) \der x = \int_{\Rz^N} \mu \der x$, and $\div f^{\Ph(h)}_t(\mu) = 0$ if $\div \mu = 0$.

\end{lem}
\begin{proof}

The first point follows immediately from Lemma~\ref{le:flowopwd}, since $f^{v}_t(\mu)$ is an average of push--forwards of $\mu$ with respect to the probability measure $v$.

For the second point, let $\mu \in L^{\infty}(\Rz^N)$. For $v$--almost every $X \in \Xc^{2,\infty}(\Rz^N)$, the vector--field $X - \ol v$ is divergence--free and Lipschitz by Lemma~\ref{le:avdivfree}. Hence, as in the proof of the third point of Lemma~\ref{le:eflowop}, $(\rh^{X - \ol v}_{\sqrt t})_*\mu$ has the density $\mu \circ \rh^{-(X - \ol v)}_{\sqrt t}$, whose $L^{\infty}$--norm equals $\|\mu\|_{\infty}$. Hence, we readily show that $f^{v}_t(\mu) \in L^{\infty}(\Rz^N)$ and $\|f^{v}_t(\mu)\|_{\infty} \leq \|\mu\|_{\infty}$.

Now, let $v = \Ph(h)$. Then, by Corollary~\ref{co:avPhi}, for all $i \in \llbracket 0, N-1 \rrbracket$, $\pm\eta X_i - h - \ol{\Ph(h)} = \pm \eta X_i$, hence
\begin{equation}
    f^{\Ph(h)}_t(\mu) = \frac{1}{2N}\sum_{i = 0}^{N-1} \bigl(\rh^{\eta X_i}_{\sqrt t}\bigr)_*\mu + \bigl(\rh^{-\eta X_i}_{\sqrt t}\bigr)_*\mu.
\end{equation}

For the third point, let $t > s >0$ and denote by $\bigl(e_i\bigr)_{i \in \llbracket 0, N-1 \rrbracket}$ the canonical basis of $\Rz^N$ and let $i \in \llbracket 0, N-1 \rrbracket$. Let $f \in L^{\infty}\bigl(\Rz^N\bigr) \cap \Lip\bigl(\Rz^N\bigr)$, such that $\max(\|f\|_{\infty}, \|f\|_{\Lip}) \leq 1$.
Then, 
\begin{equation}
    \begin{aligned}
        &\int_{\Rz^N} f \der f^{\Ph(h)}_t(\mu) - \int_{\Rz^N} f \der f^{\Ph(h)}_s(\mu)\\
        &\qquad \leq \frac{1}{2N} \sum_{i = 0}^{N-1} \int_{\Rz^N} f \circ \rho^{\eta X_i}_{\sqrt{t}} - f \circ \rho^{\eta X_i}_{\sqrt{s}} \der \mu + \int_{\Rz^N} f \circ \rho^{-\eta X_i}_{\sqrt{t}} - f \circ \rho^{-\eta X_i}_{\sqrt{s}} \der \mu.\\
    \end{aligned}
\end{equation}
Now, for all $i \in \llbracket 0, N-1 \rrbracket$, $X_i$ is a constant vector--field, and therefore, in particular, divergence--free. Hence, 
as in the proof of the third point of Lemma~\ref{le:eflowop}, we obtain
\begin{equation}
    \int_{\Rz^N} f \circ \rho^{\eta X_i}_{\sqrt{t}} \der \mu = \int_{\Rz^N} f(x) \mu\bigl(\rho^{-\eta X_i}_{\sqrt{t}}(x)\bigr) \der x. 
\end{equation}
It follows that, 
\begin{equation}
    \begin{aligned}
        &0 = \int_{\Rz^N} f \circ \rho^{\eta X_i}_{\sqrt{t}}(x) \mu \circ \rho^{\eta X_i}_{\sqrt{t}}(x) \der x  - \int_{\Rz^N} f \circ \rho^{\eta X_i}_{\sqrt{s}}(x) \mu \circ \rho^{\eta X_i}_{\sqrt{s}}(x) \der x\\
        & =\int_{\Rz^N} \biggl(f \circ \rho^{\eta X_i}_{\sqrt{t}}(x) - f \circ \rho^{\eta X_i}_{\sqrt{s}}(x)\biggr)\biggl(\mu \circ \rho^{\eta X_i}_{\sqrt{t}}(x) - \mu \circ \rho^{\eta X_i}_{\sqrt{s}}(x)\biggr) \der x\\
        &\quad + \int_{\Rz^N} \biggl(f \circ \rho^{\eta X_i}_{\sqrt{t}}(x) - f \circ \rho^{\eta X_i}_{\sqrt{s}}(x)\biggr)\mu \circ \rho^{\eta X_i}_{\sqrt{s}}(x) \der x\\
        &\quad + \int_{\Rz^N} f \circ \rho^{\eta X_i}_{\sqrt{s}}(x)\biggl(\mu \circ \rho^{\eta X_i}_{\sqrt{t}}(x) - \mu \circ \rho^{\eta X_i}_{\sqrt{s}}(x)\biggr) \der x.\\
    \end{aligned}
\end{equation}
Now, using Fubini--Tonelli's theorem, we obtain
\begin{equation}
    \begin{aligned}
        & \biggl|\int_{\Rz^N} \biggl(f \circ \rho^{\eta X_i}_{\sqrt{t}}(x) - f \circ \rho^{\eta X_i}_{\sqrt{s}}(x)\biggr)\biggl(\mu \circ \rho^{\eta X_i}_{\sqrt{t}}(x) - \mu \circ \rho^{\eta X_i}_{\sqrt{s}}(x)\biggr) \der x \biggr|\\
        &\qquad \biggl| \leq \int_{\Rz^N} \eta^2 \bigl(\sqrt{t} - \sqrt{s}\bigr)\int_{\sqrt{s}}^{\sqrt{t}} \frac{\der}{\der x_i} \mu(x + r \eta e_i) \der r \der x \biggr|\\
        &\qquad \leq \eta^2 \bigl|\sqrt{t} - \sqrt{s}\bigr| \int_{\sqrt{s}}^{\sqrt{t}} \int_{\Rz^N} \biggl| \frac{\der}{\der x_i} \mu(x + r \eta e_i) \biggr|\der r \der x\\
        &\qquad \leq \eta^2 (t- s)\bigl\|\frac{\der} {\der x_i} \mu \bigr\|_{L^1(\Rz^N)}.
    \end{aligned}
\end{equation}
Furthermore, 
\begin{equation}
    \begin{aligned}
        &\biggl|\int_{\Rz^N} \biggl(f \circ \rho^{\eta X_i}_{\sqrt{t}}(x) - f \circ \rho^{\eta X_i}_{\sqrt{s}}(x)\biggr)\mu \circ \rho^{\eta X_i}_{\sqrt{s}}(x) \der x\\
        &\quad - \int_{\Rz^N} \biggl(f \circ \rho^{\eta X_i}_{\sqrt{t}}(x) - f \circ \rho^{\eta X_i}_{\sqrt{s}}(x)\biggr)\mu (x) \der x\biggr|\\
        &= \biggl|\int_{\Rz^N} \biggl(f \circ \rho^{\eta X_i}_{\sqrt{t}}(x) - f \circ \rho^{\eta X_i}_{\sqrt{s}}(x)\biggr)\biggl(\mu \circ \rho^{\eta X_i}_{\sqrt{s}}(x) - \mu(x) \biggr) \der x\biggr|\\
        &\leq \eta^2 \bigl|\sqrt{t}- \sqrt{s}\bigr|\sqrt{s}\bigl\|\frac{\der} {\der x_i} \mu \bigr\|_{L^1(\Rz^N)} \leq \eta^2 \bigl|t -s\bigr|\bigl\|\frac{\der} {\der x_i} \mu \bigr\|_{L^1(\Rz^N)}.
    \end{aligned}
\end{equation}
In total,  we obtain
\begin{equation}
    \begin{aligned}
        &\biggl|\int_{\Rz^N} f \der f^{\Ph(h)}_t(\mu) - \int_{\Rz^N} f \der f^{\Ph(h)}_s(\mu)\biggr|\\
        &\qquad \leq \frac{1}{2N} \sum_{i = 0}^{N-1} \biggl| \int_{\Rz^N} f \circ \rho^{\eta X_i}_{\sqrt{t}} - f \circ \rho^{\eta X_i}_{\sqrt{s}} \der \mu + \int_{\Rz^N} f \circ \rho^{-\eta X_i}_{\sqrt{t}} - f \circ \rho^{-\eta X_i}_{\sqrt{s}} \der \mu\biggr|\\
        &\qquad \leq \frac{1}{2N} \sum_{i = 0}^{N-1} 3 \eta^2 (t-s)\bigl\|\frac{\der} {\der x_i} \mu \bigr\|_{L^1(\Rz^N)} = \frac{3}{2}\eta^2\bigl\|\frac{\der} {\der x_i} \mu \bigr\|_{L^1(\Rz^N)}(t-s). \\
    \end{aligned}
\end{equation}

The fourth point follows immediately since $f^{\Ph(h)}_t(\mu)$ is a convex combination of translates of $\mu$.

For the fifth point, take $q \in (0,1)$, $t \in [0,1)$, and suppose $\mu \in L^{\infty}(\Rz^N)$. 
Then, for all $i \in \llbracket 0, N-1 \rrbracket$ and all $x \in \Rz^N \setminus \ol{B(0, 2\eta\sqrt{t})}$, $\bigl\|\rho^{\eta X_i}_{\sqrt{t}}(x)\bigr\| > \frac{1}{2}\|x\| > 0$.
In particular,  $t \mapsto \bigl\|\rho^{\eta X_i}_{t}(x)\bigr\|^q$ is $\Cc^2$ in $0$. Now, a Taylor's expansion of second order with integral remainder yields 
\begin{equation}
    \begin{aligned}
        &\int_{\Rz^N \setminus \ol{B(0, 3\eta\sqrt{t})}} \|x\|^q \der \bigl|f^{\Ph(h)}_t(\mu)\bigr|(x)\\
        &\qquad  \leq \frac{1}{2N} \sum_{i=0}^{N-1} \int_{\Rz^N \setminus \ol{B(0, 2\eta\sqrt{t})}} \bigl\|\rho^{\eta X_i}_{\sqrt{t}}(x)\bigr\|^q + \bigl\|\rho^{- \eta X_i}_{\sqrt{t}}(x)\bigr\|^q \der |\mu|(x)\\
        &\qquad \leq \Mc_q(\mu) + \frac{\eta^2 q }{2N} \sum_{i=0}^{N-1} \int_{\Rz^N \setminus \ol{B(0, 2\eta\sqrt{t})}} \int_0^{\sqrt{t}}(\sqrt{t} - s)\biggl(\bigl\|\rho^{\eta X_i}_{s}(x) \bigr\|^{q-2}\\
        &\qquad \qquad \qquad + \bigl\|\rho^{-\eta X_i}_{s}(x) \bigr\|^{q-2}\biggr)\der s \der|\mu|(x)\\
        &\qquad \leq \Mc_q(\mu) + \frac{\eta^2 q }{N}\\
        &\qquad\qquad \cdot \sum_{i=0}^{N-1} \int_{\Rz^N \setminus \ol{B(0, 2\eta\sqrt{t})}} \int_0^{\sqrt{t}}(\sqrt{t} - s)\biggl(\bigl(\bigl\|x\bigr\| - \eta \sqrt{t} \bigr)^{q-2}\biggr)\der s \der|\mu|(x)\\
        &\qquad \leq \Mc_q(\mu) + 2^{1-q}\eta^2 q\bigl(C(N,q)\|\mu\|_{\infty} + |\mu|(\Rz^N)\bigr)t,
    \end{aligned}
\end{equation}
where $C(N,q)>0$ depends only on $N$ and $q$. Note that the first order term cancels out when adding the expansions of $t \mapsto \bigl\|\rho^{\eta X_i}_{t}(x)\bigr\|^q$ and $t \mapsto \bigl\|\rho^{-\eta X_i}_{t}(x)\bigr\|^q$.
Finally,
\begin{equation}
    \begin{aligned}
        &\int_{\ol{B(0, 3\eta\sqrt{t})}} \|x\|^q \der \bigl|f^{\Ph(h)}_t\bigr|(\mu)(x)\\
        &\qquad \leq 2\bigl(4\eta \sqrt{t}\bigr)^q \|u\|_{\infty} \bigl|\ol{B(0, 4\eta\sqrt{t})} \bigr|\\
        &\qquad \leq C(N,q)\eta^q \|u\|_{\infty} t^{\frac{n - 1 + q}{2}}.
    \end{aligned}
\end{equation}
Therefore, the estimate for the $q$--moment of $f^{\Ph(h)}_t(\mu)$ follows.

The sixth point follows immediately since $f^{\Ph(h)}_t(\mu)$ is a convex combination of push--forwards of $\mu$ under the flows of constant vector--fields.

Finally let us prove the last point. We set $a_i \coloneqq \eta\sqrt{t}e_i$ for all $i \in \llbracket 0, N-1 \rrbracket$. We readily show that for $k \leq 2$, 
\begin{equation}
    \begin{aligned}
        \langle x \rangle^{\be}|\nabla^k f^{\Ph(h)}_t(\mu)(x)| &\leq |\nabla^k\mu|_{\be}\,\frac{1}{2N}\sum_i \langle x \rangle^{\be}\bigl(\langle x - a_i\rangle^{-\be} + \langle x + a_i \rangle^{-\be}\bigr)\\ &\leq (1 + C\eta^2 t)|\nabla^k \mu|_{\be}
    \end{aligned}
\end{equation} 
by \eqref{eq:symweight}. Likewise, for $x \neq y$ and $\ga \coloneqq N + 3 + \al - \ep$,
    \begin{equation}
        \begin{aligned}
            &\frac{1}{2N}\sum_i \bigl|\nabla^2\mu(x - a_i) - \nabla^2 \mu(y - a_i)\bigr| + \bigl|\nabla^2\mu(x + a_i) - \nabla^2 \mu(y + a_i)\bigr|\\
            &\qquad \leq [\nabla^2 \mu]_{\al,\ga}\|x-y\|^{\al}\frac{1}{2N}\sum_i\bigl(\omega_{\ga}(x - a_i, y - a_i) + \omega_{\ga}(x + a_i, y + a_i)\bigr).
        \end{aligned}
    \end{equation}
Now, \eqref{eq:symweight} applied at $x$ and at $y$ gives $\frac{1}{2}\bigl(\omega_{\ga}(x - a, y - a) + \omega_{\ga}(x + a, y + a)\bigr) \leq (1 + C\eta^2t)\omega_{\ga}(x,y)$. The rest follows immediately, since translations commute with $\div$ and preserve the integral.
\end{proof}

Now, similar estimates for the source operators follow readily.

\begin{lem}[Properties of the source operators]\label{le:sop}
Let $T>0$, $N \in \Nz$ and $d \in \Cc^1\bigl([0,T) \times \Rz^N\bigr)$. For all $T_0 \in [0,T)$, $t,s >0$, $i \in \llbracket 0, N-1 \rrbracket$ and all $\mu \in \Wz_1(\Rz^N)$, the following holds:
\begin{itemize}
    \item $s^{d,T_0}_{t,i}(\mu) \in \Wz_1(\Rz^N)$, and 
    \begin{equation}
        \bigl|s^{d,T_0}_{t,i}(\mu)\bigr|(\Rz^N) = |\mu|(\Rz^N) + \int_{T_0}^{T_0 + t} \|\nabla d(s)_i\|_{L^1} \der s.
    \end{equation}
    \item $\|s^{d,T_0}_{s,i}(\mu) - s^{d,T_0}_{t,i}(\mu)\|_{\Wz_1} \leq \int_{T_0 + s}^{T_0 + t} \bigl\|\nabla d(r)_i\bigr\|_{L^1} \der r$.
    \item if $\mu \in L^{\infty}(\Rz^N)$, then $s^{d,T_0}_{t,i}(\mu) \in L^{\infty}(\Rz^N)$ and 
    \begin{equation}
        \|s^{d,T_0}_{t,i}(\mu)\|_{\infty} \leq \|\mu\|_{\infty} + \int_{T_0}^{T_0 + t} \|\nabla d(s)_i\|_{\infty} \der s.
    \end{equation}
    \item if there exists $\al \in (0,1]$ such that $\mu \in \Cc^{\al}(\Rz^N)$, then $s^{d,T_0}_{t,i}(\mu) \in \Cc^{\al}(\Rz^N)$ and 
    \begin{equation}
        \bigl\{s^{d,T_0}_{t,i}(\mu)\bigr\}_{\Cc^{\al}} \leq \{\mu\}_{\Cc^{\al}} + \int_{T_0}^{T_0 + t} \{\nabla d(s)_i\}_{\Cc^{\al}} \der s.
    \end{equation}
    \item let $q \in (0,1)$, and suppose $\mu$ has finite $q$--moment and is bounded, then $s^{d,T_0}_{t,i}(\mu)$ has finite $q$--moment and 
        \begin{equation}
            \begin{aligned}
                &\Mc_q\bigl(s^{d,T_0}_{t,i}(\mu)\bigr) \leq \Mc_q(\mu) + \int_{T_0}^{T_0 + t} \Mc_q\bigl(\nabla d(s)_i\bigr) \der s.
            \end{aligned}
        \end{equation} 
    \item suppose $\mu \in W^{1,\infty}$ and $ \nabla \mu \in L^1(\Rz^N)$.  Then, $s^{d,T_0}_{t,i}(\mu) \in W^{1,\infty}$, $ \nabla s^{d,T_0}_{t,i}(\mu) \in L^1(\Rz^N)$ and
    \begin{equation}
        \|\nabla s^{d,T_0}_{t,i}(\mu)\|_{L^1} \leq \|\nabla \mu \|_{L^1} + \int_{T_0}^{T_0 + t} \bigl\|\nabla\bigl(\nabla d(s)_i\bigr)\bigr\|_{L^1} \der s.
    \end{equation}
\end{itemize}
\end{lem}
\begin{proof}
The only point that does not follow immediately from the definition of $s^{d,T_0}_{t,i}(\mu)$ is the second one. Let $t > s >0$ and let $f \in L^{\infty}\bigl(\Rz^N\bigr) \cap \Lip\bigl(\Rz^N\bigr)$, such that $\max(\|f\|_{\infty}, \|f\|_{\Lip}) \leq 1$. Then,
\begin{equation}
    \begin{aligned}
        &\biggl|\int_{\Rz^N} f \der s^{d,T_0}_{t,i}(\mu) - \int_{\Rz^N} f \der s^{d,T_0}_{s,i}(\mu)\biggr| = \biggl|\int_{T_0 + s}^{T_0 + t} \int_{\Rz^N} f \der \bigl(\nabla d(r)_i\bigr) \der r\biggr|\\
        &\qquad \leq \int_{T_0 + s}^{T_0 + t} \int_{\Rz^N} |\nabla d(r)_i(x)| \der x \der r \leq \int_{T_0 + s}^{T_0 + t} \bigl\|\nabla d(r)_i\bigr\|_{L^1} \der r.
    \end{aligned}  
\end{equation}
Thus, the statement follows.
\end{proof}

\subsubsection{Properties of the dievergence operator}

First we establish we establish some basic identities for the dievergence operators and study on which subspaces of $L^1(\Rz^n, \Rz^n)$ they are bounded. 
The points in Lemma~\ref{le:dievwd} are mostly well-known results from Calderon--Zygmund theory.  

\begin{lem}\label{le:dievwd}
    Let $N \in \Nz$, $N \geq 2$, $\theta \in [0,1]$, and let $(R_j)_{j \in \llbracket 0, N-1 \rrbracket}$ be the Riesz transforms on $\Rz^N$, see \cite[Chapter III, \S 1]{Stein1970} and \cite[Section 5.1]{Grafakos2014}. Then:
    \begin{itemize}
        \item for all $\nu \in \Sc(\Rz^N, \Rz^N)$ and $j \in \llbracket 0, N-1 \rrbracket$,
        \begin{equation}\label{eq:dievriesz}
            s_{\theta}(\nu)_j = \nu_j + \theta \sum_{k=0}^{N-1} R_jR_k \nu_k, \qquad \widehat{s_{\theta}(\nu)}(\xi) = \hat\nu(\xi) - \theta\,\frac{\xi\,\bigl(\xi \cdot \hat\nu(\xi)\bigr)}{|\xi|^2}.
        \end{equation}
        In particular, $s_{\theta}(\nu) \in \Cc^{\infty}(\Rz^N, \Rz^N) \cap L^p(\Rz^N, \Rz^N)$ for all $p \in (1, \infty]$.
        \item For all $p \in (1,\infty)$, $s_{\theta}$ extends uniquely to a bounded linear operator on $L^p(\Rz^N,\Rz^N)$, and there exists $C(N,p) > 0$ such that $$\|s_{\theta}(\nu)\|_{L^p} \leq (1 + \theta C(N,p))\|\nu\|_{L^p}$$ for all $\nu \in L^p(\Rz^N, \Rz^N)$. Moreover, 
        \begin{itemize}
            \item $s_{\theta}$ commutes with translations and with differentiation,
            \item $\div s_{\theta}(\nu) = (1-\theta)\div \nu$ in the sense of distributions,
            \item $s_{\theta}(\nu) = \nu$ whenever $\div \nu = 0$,
            \item $s = s_1$ is the Leray--Helmholtz projection, i.e., $s \circ s = s$, $s(\nabla \psi) = 0$ for all $\psi$ with $\nabla \psi \in L^p(\Rz^N, \Rz^N)$, and $s(\nu) = \nu$ if and only if $\div \nu = 0$.
        \end{itemize}
    \end{itemize}
\end{lem}
\begin{proof}
    Let $\nu \in \Sc(\Rz^N, \Rz^N)$. Then, $\div \nu \in \Sc(\Rz^N)$ and $\widehat{\div \nu}(\xi) = i\,\xi \cdot \hat\nu(\xi)$. Hence, $\widehat{\Delta^{-1} \div \nu}(\xi) = -i\,\xi\cdot\hat\nu(\xi)/|\xi|^2$, which is locally integrable since $N \geq 2$. Taking the gradient yields \eqref{eq:dievriesz}, since the symbol of $R_jR_k$ is $-\xi_j\xi_k/|\xi|^2$. 
    Now, $\xi \mapsto \xi_j\xi_k/|\xi|^2$ is smooth and homogeneous of degree $0$ on $\Rz^N \setminus \{0\}$, hence $s_{\theta}(\nu)$ is smooth and belongs to $L^p$ for all $p \in (1,\infty]$ by the $L^p$--boundedness of the Riesz transforms \cite[Chapter III, \S1]{Stein1970}, \cite[Section 5.1 and Corollary 5.2.8]{Grafakos2014}, together with the Sobolev embedding for $p = \infty$.

    The second point is the content of the Calder\'on--Zygmund theory for the Riesz transforms, \cite[Chapter II, Theorem 3 and Chapter III, Theorem 4]{Stein1970}, \cite[Corollary 5.2.8]{Grafakos2014}. The properties of $s = s_1$ are those of the Helmholtz--Weyl decomposition $L^p(\Rz^N,\Rz^N) = \{\nu \in L^p : \div \nu = 0\} \oplus \{\nabla \psi \in L^p\}$, see \cite[Theorem III.1.2]{Galdi2011}. The identity $\div s_{\theta}(\nu) = \div \nu - \theta \Delta \Delta^{-1} \div \nu = (1 - \theta)\div \nu$ holds on $\Sc$ and extends by density and continuity.
\end{proof}

\begin{re}
    In particular, $s_i$ is well--defined on $L^1(\Rz^N, \Rz^N) \cap L^p(\Rz^N, \Rz^N)$ with values in $L^p(\Rz^N, \Rz^N)$. However, for $i > 0$, $s_i$ does not extend to a bounded operator on $L^1(\Rz^N, \Rz^N)$, and $s_i(\nu) \notin L^1(\Rz^N,\Rz^N)$ in general. 
    The Riesz transforms are not bounded on $L^1(\Rz^N)$, \cite[Chapter II, \S 4]{Stein1970}, \cite[Section 5.1]{Grafakos2014}. Indeed, if $0 \leq \nu_k \in L^1(\Rz^N) \cap L^{\infty}(\Rz^N)$ is compactly supported with $\int_{\Rz^N} \nu_k \der x \neq 0$, then $R_jR_k\nu_k(x) \sim c(N)\frac{\delta_{jk}|x|^2 - Nx_jx_k}{|x|^{N+2}}\int_{\Rz^N} \nu_k \der x$ as $|x| \to \infty$, which is not integrable. However, the operator $s_i$ is bounded on the Hardy space $\Hc^1(\Rz^N, \Rz^N)$, see for instance \cite[Chapter III, \S 3]{Stein1993}.
\end{re}

We complement Lemma~\ref{le:dievwd} by the two elementary estimates on the dievergence operators, similar to those for the source operator in Lemma~\ref{le:sop}. 

\begin{lem}\label{le:dievop}
    Let $N \geq 2$, $i,j \in [0,1]$ and $\nu \in W^{1,1}(\Rz^N, \Rz^N) \cap W^{1,\infty}(\Rz^N, \Rz^N)$. Then, the following holds:
    \begin{itemize}
        \item If $\nabla \Delta^{-1}\div \nu \in L^1(\Rz^N, \Rz^N)$, then $s_i(\nu) \in L^1(\Rz^N, \Rz^N) \subset \Wz_1(\Rz^N)^N$, 
        \begin{equation}
            \begin{aligned}
                &|s_i(\nu)|(\Rz^N) \leq \|\nu\|_{L^1} + i\|\nabla \Delta^{-1}\div \nu\|_{L^1},\\  &\|s_i(\nu) - s_j(\nu)\|_{\Wz_1} \leq |i - j|\,\|\nabla \Delta^{-1}\div \nu\|_{L^1},
            \end{aligned}
        \end{equation}
        where the $\Wz_1$--norm is taken componentwise.
        \item $s_i(\nu) \in L^{\infty}(\Rz^N, \Rz^N)$ and
        \begin{equation}
            \|s_i(\nu)\|_{\infty} \leq \|\nu\|_{\infty} + i\,C(N)\|\div \nu\|_{\infty}^{1 - \frac{1}{N}}\|\div \nu\|_{L^1}^{\frac{1}{N}}.
        \end{equation}
    \end{itemize}
\end{lem}
\begin{proof}
    Recall that for $\nu \in W^{1,1}(\Rz^N, \Rz^N) \cap W^{1,\infty}(\Rz^N, \Rz^N)$, $\div \nu \in L^1(\Rz^N) \cap L^{\infty}(\Rz^N)$, and $\nabla \Delta^{-1} \div \nu = \nabla \Gamma * \div \nu$ with $\Gamma$ the fundamental solution of the Laplacian.

    The first point follows from $s_i(\nu) - s_j(\nu) = (j - i)\nabla\Delta^{-1}\div \nu$ and from $\|f\|_{\Wz_1} \leq \|f\|_{L^1}$ for $f \in L^1(\Rz^N)$.

    For the second point, write $F \coloneqq \div \nu \in L^1 \cap L^{\infty}$. Since $|\nabla \Gamma(y)| \leq C(N)\|y\|^{1-N}$, for all $r > 0$ and $x \in \Rz^N$,
    \begin{equation}
        \begin{aligned}
            &\bigl|\nabla \Gamma * F(x)\bigr| \leq C(N)\int_{\|y\| < r}\|y\|^{1-N}|F(x-y)| \der y\\
            &\qquad + C(N)\int_{\|y\| \geq r}\|y\|^{1-N}|F(x-y)| \der y \leq C(N)\bigl(r\|F\|_{\infty} + r^{1-N}\|F\|_{L^1}\bigr),
        \end{aligned}
    \end{equation}
    and the choice $r = \bigl(\|F\|_{L^1}/\|F\|_{\infty}\bigr)^{\frac{1}{N}}$ yields the result.
\end{proof}

The following lemma is the Euclidean, Hölder--space version of the Fredholm theory of the Laplacian on weighted spaces, see \citet{McOwen1979}, \citet[Sections 1--2]{Bartnik1986} and \citet{LockhartMcOwen1985}: for decay rates between $|x|^{-N}$ and $|x|^{1-N}$, the range of $\Delta$ consists exactly of the functions which are orthogonal to the harmonic polynomials of degree at most one, which is the condition \eqref{eq:Fmoments} below. Since the weighted Hölder norms of Definition~\ref{de:weightednorms} differ slightly from the ones used there, we give a short proof, which rests on the pointwise decay of the Newtonian potential and on the scaled interior Schauder estimates of \cite[Theorems 6.2 and 6.17]{GilbargTrudinger2001}.

\begin{lem}[Dievergence operator on $\mathfrak E$]\label{le:sweighted}
    Let $N \geq 2$, $\al, \ep \in (0,1)$. There exists $K_1 = K_1(N,\al,\ep) > 0$ such that for all $F \in \Cc^1(\Rz^N)$ with $\|F\|_{\mathcal F} < \infty$, satisfying for all $ j \in \llbracket 0, N-1 \rrbracket$,
    \begin{equation}\label{eq:Fmoments}
        \begin{aligned}
            &\int_{\Rz^N} F \der x = 0,\\ &\int_{\Rz^N} x_j F \der x = 0, 
        \end{aligned}
    \end{equation}
    the field $\nabla \Gamma * F = \nabla\Delta^{-1}F$ belongs to $\mathfrak E$, has mean zero, and
    \begin{equation}
        \|\nabla \Gamma * F\|_{\mathcal E} \leq K_1 \|F\|_{\mathcal F}.
    \end{equation}
    Therefore, if $\nu \in \mathfrak E$ has mean zero and $F \coloneqq \div \nu$ satisfies $\|F\|_{\mathcal F} < \infty$, for all $\theta \in [0,1]$, $s_{\theta}(\nu) \in \mathfrak E$ has mean zero, 
    \begin{equation}
        \|s_{\theta}(\nu)\|_{\mathcal E} \leq \|\nu\|_{\mathcal E} + \theta K_1\|F\|_{\mathcal F},
    \end{equation}
    and $s(\nu) \in \mathcal E$.
\end{lem}
\begin{proof}
    Throughout the proof, $C$ denotes constants depending only on $N$, $\al$ and $\ep$, and we write $\be \coloneqq N - \ep$, so that $|F(x)| \leq \|F\|_{\mathcal F}\langle x \rangle^{-\be - 2}$, $|\nabla F(x)| \leq \|F\|_{\mathcal F}\langle x \rangle^{-\be-3}$, and $\ga \coloneqq \be + 3 + \al$. Let $u \coloneqq \Gamma * F$. Since $F \in \Cc^{1,\al}(\Rz^N) \cap W^{1,1}(\Rz^N)$, \cite[Lemmas 4.2 and 4.4]{GilbargTrudinger2001}, applied to $F$ and to $\nabla F$, show that $u \in \Cc^{3,\al}_{\mathrm{loc}}(\Rz^N)$, $\Delta u = F$ and $\nabla u = \nabla \Gamma * F \in \Cc^{2,\al}(\Rz^N, \Rz^N)$.

    \bigskip

    \textbf{Step 1: } 

    \bigskip
    
    We show that $|u(x)| \leq C\|F\|_{\mathcal F}\langle x \rangle^{-\be}$ for all $x \in \Rz^N$.
    For $\|x\| \leq 1$ this follows from $|\Gamma(z)| \leq C(\|z\|^{2-N} + |\log\|z\||)$ for all $z \in \Rz^N$, which ensures integrability. For $\|x\| > 1$, \eqref{eq:Fmoments} allows us to write
    \begin{equation}
        u(x) = \int_{\Rz^N}\Bigl(\Gamma(x - y) - \Gamma(x) + \nabla\Gamma(x)\cdot y\Bigr)F(y) \der y.
    \end{equation}
    On $\{\|y\| \leq \|x\|/2\}$, the integrand is bounded by $C\|y\|^2\|x\|^{-N}|F(y)|$ by Taylor's formula. Now, since $\int_{\|y\| \leq \|x\|/2}\|y\|^2\langle y \rangle^{-\be-2}\der y \leq C\|x\|^{\ep}$, we obtain the bound $C\|F\|_{\mathcal F}\|x\|^{-\be}$. 
    
    On $\{\|y\| > \|x\|/2\}$, we have $|\Gamma(x - y) - \Gamma(x)| \leq C\bigl(\|x - y\|^{2-N} + \|x\|^{2-N}\bigr)$ for $N \geq 3$, $|\Gamma(x-y) - \Gamma(x)| \leq C\bigl(1 + |\log(\|x - y\|/\|x\|)|\bigr)$ for $N = 2$, and $|\nabla\Gamma(x)\cdot y| \leq C\|x\|^{1-N}\|y\|$. Now, substituting $y = \|x\|w$ and using $\langle y \rangle^{-\be-2} \leq \|x\|^{-\be-2}\|w\|^{-\be-2}$, we obtain
    \begin{equation}
        \begin{aligned}
            &C\|F\|_{\mathcal F}\|x\|^{-\be}\int_{\|w\| > 1/2}\Bigl(\bigl\|\tfrac{x}{\|x\|} - w\bigr\|^{2-N}\\
            &\qquad \qquad \qquad  + 1 + \bigl|\log\bigl\|\tfrac{x}{\|x\|} - w\bigr\|\bigr| + \|w\|\Bigr)\|w\|^{-\be-2} \der w \leq C\|F\|_{\mathcal F}\|x\|^{-\be}.
        \end{aligned}
    \end{equation}

    \bigskip

    \textbf{Step 2:} 
    
    \bigskip

    Let $x \in \Rz^N$ and $R \coloneqq \frac{3}{8}\langle x \rangle$, so that $\frac{1}{4}\langle x \rangle \leq \langle z \rangle \leq \frac{7}{4}\langle x \rangle$ for all $z \in B(x, 2R)$. The function $v: z \mapsto u(x + Rz)$ satisfies $\Delta v(z) = R^2 F(x + Rz)$ for $z \in B(0,2)$. Thus, $\Delta \nabla v(z) = R^3(\nabla F)(x + Rz)$, and the interior Schauder estimates \cite[Theorems 6.2 and 6.17]{GilbargTrudinger2001}, applied to $v$ and to $\nabla v$, yield after rescaling
    \begin{equation}\label{eq:scaledschauder}
        \begin{aligned}
            &\sum_{j = 1}^{3} R^j|\nabla^j u(x)| + R^{3 + \al}\{\nabla^3 u\}_{\Cc^{\al}(B(x,R))}\\
            &\qquad \leq C\Bigl(\sup_{B(x,2R)}|u| + R^2\sup_{B(x,2R)}|F| + R^3\sup_{B(x,2R)}|\nabla F| + R^{3 + \al}\{\nabla F\}_{\Cc^{\al}(B(x,2R))}\Bigr).
        \end{aligned}
    \end{equation}
    By Step 1 and using $\omega_{\ga}(y,y') \leq C\langle x \rangle^{-\ga}$ for $y,y' \in B(x,2R)$, the right--hand side is bounded by $C\|F\|_{\mathcal F}\langle x \rangle^{-\be}$. Hence \eqref{eq:scaledschauder} yields $|\nabla^{k}(\nabla u)|_{\be + 1 + k} \leq C\|F\|_{\mathcal F}$ for $k \in \llbracket 0, 2 \rrbracket$. Moreover, for all $y \in B(x,R)$, $|\nabla^3 u(x) - \nabla^3 u(y)| \leq C\|F\|_{\mathcal F}\langle x \rangle^{-\ga}\|x - y\|^{\al} \leq C\|F\|_{\mathcal F}\|x - y\|^{\al}\omega_{\ga}(x,y)$. 
    
    For the remaining pairs, with $\|x - y\| > \frac{3}{8}\min(\langle x \rangle, \langle y \rangle)$, the argument of the third point of Lemma~\ref{le:weights} gives $|\nabla^3u(x) - \nabla^3u(y)| \leq C|\nabla^3 u|_{\ga - \al}\|x - y\|^{\al}\omega_{\ga}(x,y)$. Thus $[\nabla^2(\nabla u)]_{\al,\ga} \leq C\|F\|_{\mathcal F}$, and altogether $\|\nabla u\|_{\mathcal E} \leq K_1\|F\|_{\mathcal F}$.

    \bigskip

    \textbf{Step 3: } 

    \bigskip
    
    By Step 2, $\nabla u \in L^1(\Rz^N, \Rz^N)$, and $\int_{B(0,R)}\nabla u \der x = \int_{\partial B(0,R)}u\,n \der S$ is bounded by $CR^{N-1}\sup_{\partial B(0,R)}|u| \leq C\|F\|_{\mathcal F}R^{-1+\ep}$ by Step 1, hence tends to zero as $R \to \infty$.

    Finally, let $\nu \in \mathfrak E$ have mean zero and set $F \coloneqq \div\nu$. Then $\int_{\Rz^N} F \der x = 0$, as $\nu \in W^{1,1}$, and $\int_{\Rz^N} x_j F \der x = -\int_{\Rz^N} \nu_j \der x = 0$ by the boundary argument of Lemma~\ref{le:eflowop}. Hence, the above applies to $F$, and the statements on $s_{\theta}(\nu) = \nu - \theta\nabla\Gamma * F$ follow, $\div s(\nu) = 0$ being Lemma~\ref{le:dievwd}.
\end{proof}

\subsection{The incompressible Navier-Stokes equations and Incompressible flow approximation series (IFAS)}\label{sbs:NS}

In this section, we use the evolution operators defined previously to build approximation series for solutions of the incompressible Navier--Stokes equations. We shall apply the estimates we have established for the evolution operators in the previous section to prove that these approximations uniformly converge to a solution on a given time interval. First, we state the incompressible Navier--Stokes equations and recall some well--known properties of their solutions.

Let $N \in \Nz$ with $N \geq 2$, $\eta > 0$ and $T > 0$. For a time--dependent vector--field $g \in \Cc^1\bigl([0,T) \times \Rz^N, \Rz^N\bigr)$, we apply the Lie derivative along $g(t)$ componentwise and write
\begin{equation}
    \Lc_g g \coloneqq \bigl(\Lc_{g(t)} g(t)_i\bigr)_{i \in \llbracket 0, N-1 \rrbracket}.
\end{equation}
The incompressible Navier--Stokes equations with viscosity $\frac{\eta^2}{2N}$, velocity $g : [0,T) \times \Rz^N \to \Rz^N$, pressure $d : [0,T) \times \Rz^N \to \Rz$ and initial condition $g_0$ read
\begin{equation}\label{eq:NS}
    \begin{cases}
        \partial_t g = \Lc_g g + \dfrac{\eta^2}{2N}\Delta g + \nabla d & \text{on } [0,T) \times \Rz^N,\\[1ex]
        \div g = 0 & \text{on } [0,T) \times \Rz^N,\\[0.5ex]
        g(0) = g_0.
    \end{cases}
\end{equation}
A pair $(g, d)$ with $g \in \Cc^1\bigl([0,T) \times \Rz^N, \Rz^N\bigr) \cap \Cc\bigl([0,T), \Cc^2(\Rz^N,\Rz^N)\bigr)$ and $d \in \Cc\bigl([0,T), \Cc^1(\Rz^N)\bigr)$ satisfying \eqref{eq:NS} pointwise is called a classical solution of \eqref{eq:NS} on $[0,T)$. 

\begin{re}
Note that $(g,d)$ is a classical solution of \eqref{eq:NS} with initial condition $g_0$ if and only if $(-g, d)$ is a classical solution of the Navier--Stokes equations in their usual form $\partial_t \tilde g + \Lc_{\tilde g} \tilde g = \frac{\eta^2}{2N}\Delta \tilde g - \nabla d$, $\div \tilde g = 0$, with initial condition $-g_0$. The sign convention in \eqref{eq:NS} is adapted to the operators of Section~\ref{sbs:ops}.
\end{re}
Taking the divergence of \eqref{eq:NS}, the pressure satisfies
\begin{equation}\label{eq:pressure}
    \Delta d = - \div\bigl(\Lc_g g\bigr) = - \sum_{j,k=0}^{N-1} \frac{\der g_k}{\der x_j}\frac{\der g_j}{\der x_k} \qquad \text{on } [0,T) \times \Rz^N.
\end{equation}

\begin{pr}[Classical solutions with Schwarz initial condition]\label{pr:regclasssol}
    Let $T, \eta >0$, $g_0 \in \Sc(\Rz^N, \Rz^N)$ with $\div g_0 = 0$ and let $(g, d)$ be a classical solution of the incompressible Navier--Stokes equations \eqref{eq:NS} on $[0,T)$ with initial condition $g_0$. Suppose furthermore that $g$ is bounded on $[0,T_0] \times \Rz^N$ for all $T_0 < T$ and $d(t,x) \to 0$ for $|x| \to \infty$. Then,
    $(g,d)$ satisfies the following regularity properties on $[0,T)$:
    \begin{itemize}
        \item $g \in \Cc^{\infty}([0,T) \times \Rz^N)$ and $g \in \Cc\bigl([0,T_0], W^{k,1}(\Rz^N) \cap W^{k,\infty}(\Rz^N)\bigr)$ for all $T_0 < T$ and all $k \in \Nz$. In particular, $g \in \Cc\bigl([0,T_0], \Cc^2(\Rz^N)\bigr)$ and $g \in \Cc\bigl([0,T_0], W^{1,1}(\Rz^N)\bigr)$ for all $T_0 < T$. Moreover, for all $0 < q < 1$, all $i \in \llbracket 0, N-1 \rrbracket$ and all $t \in [0,T)$, $g(t)_i$ has finite $q$--moment.
        \item $d \in \Cc^{\infty}([0,T) \times \Rz^N)$, $\nabla d \in \Cc\bigl([0,T_0], W^{1,1}(\Rz^N) \cap W^{1,\infty}(\Rz^N)\bigr)$ and $d \in \Cc\bigl([0,T_0], \Cc^2(\Rz^N)\bigr)$ for all $T_0 < T$. Moreover, for $0<q<1$, all $i \in \llbracket 0, N-1 \rrbracket$ and all $t \in [0,T)$, $\nabla d(t)_i$ has finite $q$--moment, and $t \mapsto \Mc_q(\nabla d(t)_i)$ is continuous on $[0,T)$.
    \end{itemize}
\end{pr}
\begin{proof}
    Local existence, uniqueness in the class of bounded classical solutions, and $\Cc^{\infty}$--smoothness of solutions with initial data in $\bigcap_m H^m(\Rz^N)$, together with $g \in \Cc([0,T_0], H^m(\Rz^N))$ for all $m \in \Nz$, are classical, see \cite[Chapter 3, Theorems 3.4 and 3.6]{MajdaBertozzi2002} and \cite{Kato1972}; the Sobolev embedding yields $g \in \Cc([0,T_0], W^{k,\infty}(\Rz^N))$ for all $k$. 
    
    The spatial decay $|g(t,x)| + |\nabla d(t,x)| \leq C(T_0)(1 + |x|)^{-N-1}$, together with the corresponding decay of the spatial derivatives and the fact that the decay rate $|x|^{-N-1}$ is optimal, follow from the asymptotic profiles of \cite{DobrokhotovShafarevich1994, Brandolese2004, BrandoleseVigneron2007}, in particular, the $L^1$--integrability of $\nabla^k g(t)$, $\nabla d(t)$ and $\nabla^2 d(t)$. 
    
    The finiteness of the $q$--moments of $g(t)_i$ and $\nabla d(t)_i$ for $q < 1$ follow directly. The continuity in time of the $q$--moments follows from the continuity of $t \mapsto \nabla d(t)$ in the weighted spaces used in \cite{BrandoleseVigneron2007}. Finally, the regularity of $d$ follows from that of $g$ via \eqref{eq:pressure} and Lemma~\ref{le:dievwd}, since $\nabla d = -\nabla \Delta^{-1}\div(\Lc_g g)$ is the only solution of \eqref{eq:pressure} vanishing at infinity.
\end{proof}

\begin{re}
In general, $g(t) \notin \Sc(\Rz^N, \Rz^N)$ for $t > 0$. Even if $g_0 \in \Sc(\Rz^N,\Rz^N)$, $g(t,x)$ and $\nabla d(t,x)$ decay like $|x|^{-N-1}$ as $|x| \to \infty$, and not faster. Let $\al \in (0,1)$, $\ep > 0$. Note that $\mathfrak E \subset \Cc^{2,\al}_b(\Rz^N, \Rz^N) \cap W^{2,1}(\Rz^N, \Rz^N)$, that every divergence--free $g_0 \in \Sc(\Rz^N, \Rz^N)$ belongs to $\mathcal E$, and that $\mathfrak E$ is exactly the class of fields having the decay of Proposition~\ref{pr:regclasssol} up to the loss $\ep$. 
\end{re}

Now, we shall build approximations of solutions of the incompressible Navier-Stokes equations using the operators introduced in Section~\ref{sbs:ops}. The construction is inspired by the Average Flow Approximation Series (AFAS) introduced in \citet{Schichl2025}.
To approximate solutions of the equation, we introduce the following analogon of AFAS, which we call Incompressible Flow Approximation Series (IFAS).

\begin{defi}[Incompressible Flow Approximation Series]\label{de:IFAS}
Let $\de>0$ and $\Pf(\de)$ be a partition of $[0,T]$ with maximal step--length $\pf(\Pf(\de)) \leq \de$. Let $M \coloneqq |\Pf(\de)|$ and index the elements of $\Pf(\de)$ in increasing order. Define for all $i \in \llbracket 0, N-1\rrbracket$, 
$$
\mu^{\Pf(\de)}(0)_i = \bigl(g_0\bigr)_i,
$$
and for all $l \in \llbracket 0,M-1 \rrbracket$,
$$
\mu^{\Pf(\de)}(x_{l+1}) = s \circ \biggl[\biggl( e^{\Phi\bigl(\mu^{\Pf(\de)}(x_l)\bigr)}_{x_{l+1} - x_l}\circ f^{\Phi\bigl(\mu^{\Pf(\de)}(x_l)\bigr)}_{x_{l+1} - x_l}\biggr)\bigl(\mu^{\Pf(\de)}(x_l)_i\bigr)\biggr]_{i \in \llbracket 0, N-1 \rrbracket}.
$$
Furthermore, define $\mu^{\Pf(\de)}(t)$ by 
$$
\mu^{\Pf(\de)}(t)_i = f^{\Phi\bigl(\mu^{\Pf(\de)}(x_l)\bigr)}_{3(t - x_l)}\bigl(\mu^{\Pf(\de)}(x_l)_i\bigr),
$$
for all $t \in \bigl[x_l, \frac{x_{l+1} + 2 x_l}{3}\bigr]$, 
$$
\mu^{\Pf(\de)}(t)_i = \biggl(e^{\Phi\bigl(\mu^{\Pf(\de)}(x_l)\bigr)}_{3\bigl(t - \frac{x_{l+1} + 2x_l}{3}\bigr)}\circ f^{\Phi\bigl(\mu^{\Pf(\de)}(x_l)\bigr)}_{x_{l+1} - x_l}\biggr)\bigl(\mu^{\Pf(\de)}(x_l)_i\bigr),
$$
for all $t \in \bigl[ \frac{x_{l+1} + 2 x_l}{3}, \frac{2 x_{l+1} + x_l}{3} \bigr]$, and 
$$
\mu^{\Pf(\de)}(t) = s_{3\frac{t - \frac{2x_{l +1} + x_l}{3}}{x_{l+1} - x_l}} \circ \biggl[\biggl( e^{\Phi\bigl(\mu^{\Pf(\de)}(x_l)\bigr)}_{x_{l+1} - x_l}\circ f^{\Phi\bigl(\mu^{\Pf(\de)}(x_l)\bigr)}_{x_{l+1} - x_l}\biggr)\bigl(\mu^{\Pf(\de)}(x_l)_i\bigr)\biggr]_{i \in \llbracket 0, N-1 \rrbracket}
$$
for all $t \in \bigl[ \frac{2 x_{l+1} + x_l}{3}, x_{l+1} \bigr]$.
\end{defi}

We shall see that IFAS converge uniformly to a classical solution of the incompressible Navier--Stokes equations on a time interval depending only on the $\Ec$--norm of the initial condition and the viscosity. But first, we introduce a linearized version of IFAS, which lead to interesting estimates for diverse norms of the solution of the incompressible Navier--Stokes equations.

\subsubsection{Linear incompressible flow approximation series (LIFAS)}

Analogously to the Linear average fow approximation series in \citet{Schichl2025}, we introduce the Linear incompressible flow approximation series (LIFAS), an approximation series of the linearized Incompressible Navier--Stokes equations. 

Let $(g,d)$ be a classical solution of \eqref{eq:NS} on $[0,T)$. Since $\div g = 0$, we have 
\begin{equation}\label{eq:divlieprop}
     \Lc_g g_i = \bigl<g,\nabla g_i\bigr> = \div(g\, g_i),
\end{equation}
such that the components $u_i \coloneqq g_i$ of $g$ solve the \emph{linear} equation
\begin{equation}\label{eq:FP}
    \partial_t u_i = \frac{\eta^2}{2N}\Delta u_i + \div(g\,u_i) + \frac{\der d}{\der x_i} \qquad \text{on } [0,T) \times \Rz^N,
\end{equation}
in which the velocity $g$ and the pressure $d$ are regarded as given. Equation \eqref{eq:FP} is a Fokker--Planck equation with diffusion coefficient $\frac{\eta^2}{2N}$, bounded Lipschitz drift $-g$ and source term $\frac{\der d}{\der x_i}$.

\begin{re}
    Without the source term, \eqref{eq:FP} is the forward Kolmogorov equation for the law of the diffusion process $\der X_t = -g(t, X_t)\der t + \frac{\eta}{\sqrt N}\der B_t$, and its natural solutions are curves of finite (signed) measures $t \mapsto \mu(t)$, which solve \eqref{eq:FP} in the weak sense if for all $\ph \in \Cc^{\infty}_c(\Rz^N)$ and $t \in [0,T)$,
    \begin{equation}\label{eq:FPweak}
        \begin{aligned}
            &\int_{\Rz^N}\ph \der \mu(t) - \int_{\Rz^N}\ph \der \mu(0)\\
            &\qquad = \int_0^t \biggl(\int_{\Rz^N}\Bigl(\frac{\eta^2}{2N}\Delta \ph - \Lc_{g(r)}\ph\Bigr)\der\mu(r) + \int_{\Rz^N}\ph\,\frac{\der d(r)}{\der x_i} \der x\biggr)\der r,
        \end{aligned}
    \end{equation}
    see \cite[Chapter 6 and Chapter 9]{BKRS2015} for the general theory and for the uniqueness of measure--valued solutions.
\end{re}

\begin{defi}[Linear Incompressible Flow Approximation Series]\label{de:LIFAS}
Let $\de>0$ and $\Pf(\de)$ be a partition of $[0,T]$ with maximal step--length $\pf(\Pf(\de)) \leq \de$. Let $M \coloneqq |\Pf(\de)|$ and index the elements of $\Pf(\de)$ in increasing order. Define for all $i \in \llbracket 0, N-1\rrbracket$, 
$$
\mu^{\Pf(\de)}(0)_i = \bigl(g_0\bigr)_i,
$$
and for all $l \in \llbracket 0,M-1 \rrbracket$,
$$
\mu^{\Pf(\de)}(x_{l+1})_i = \biggl(s^{d, x_l}_{x_{l+1} - x_l,i} \circ e^{\Phi(x_l,g)}_{x_{l+1} - x_l}\circ f^{\Phi(x_l,g)}_{x_{l+1} - x_l}\biggr)\bigl(\mu^{\Pf(\de)}(x_l)_i\bigr).
$$
Furthermore, define $\mu^{\Pf(\de)}(t)$ by 
$$
\mu^{\Pf(\de)}(t)_i = f^{\Phi(x_l,g)}_{3(t - x_l)}\bigl(\mu^{\Pf(\de)}(x_l)_i\bigr),
$$
for all $t \in \bigl[x_l, \frac{x_{l+1} + 2 x_l}{3}\bigr]$, 
$$
\mu^{\Pf(\de)}(t)_i = \biggl(e^{\Phi(x_l,g)}_{3\bigl(t - \frac{x_{l+1} + 2x_l}{3}\bigr)}\circ f^{\Phi(x_l,g)}_{x_{l+1} - x_l}\biggr)\bigl(\mu^{\Pf(\de)}(x_l)_i\bigr),
$$
for all $t \in \bigl[ \frac{x_{l+1} + 2 x_l}{3}, \frac{2 x_{l+1} + x_l}{3} \bigr]$, and 
$$
\mu^{\Pf(\de)}(t)_i = \biggl(s^{d, x_l}_{3\bigl(t - \frac{2x_{l+1} + x_l}{3}\bigr),i} \circ e^{\Phi(x_l,g)}_{x_{l+1} - x_l}\circ f^{\Phi(x_l,g)}_{x_{l+1} - x_l}\biggr)\bigl(\mu^{\Pf(\de)}(x_l)_i\bigr).
$$
for all $t \in \bigl[ \frac{2 x_{l+1} + x_l}{3}, x_{l+1} \bigr]$.
\end{defi}

The three operators of Section~\ref{sbs:ops} correspond exactly to the three terms of \eqref{eq:FP}:
\begin{itemize}
    \item  the averaged flow dispersion operators generate the diffusion $\frac{\eta^2}{2N}\Delta$,
    \item the expectational flow operators with $\ol{\Ph(g(t))} = -g(t)$ generate the transport $\div(g\,\cdot)$,
    \item the source operators add $\int \frac{\der d}{\der x_i}$.
\end{itemize}
The (LIFAS) defined above is thus a splitting scheme for \eqref{eq:FP}, and Theorem~\ref{th:LIFAS} states that it converges to the solution $g$ of \eqref{eq:FP}, in the $\Wz_1$--norm. The identification of the limit rests on the uniqueness of measure--valued solutions of \eqref{eq:FPweak}.
While we expect LIFAS to have limited use cases in numerics, since we are a priori required to know a solution of the Navier--Stokes equations to build them, they are nevertheless 
useful to provide intuitive proofs for estimates of several different norms of that solution. They shall be based on the estimates for LIFAS proved below.


\begin{lem}\label{le:massmom}
    Let $T_0 \in [0,T)$ and let $\Pf$ be a partition of $[0,T_0]$ with $\pf(\Pf) < 1$, and let $\mu^{\Pf}$ be the LIFAS of Definition~\ref{de:LIFAS} with initial value $g_0$. Then, for all $i \in \llbracket 0, N-1 \rrbracket$ and all $t \in [0,T_0]$, 
    \begin{itemize}
        \item we have
        \begin{equation}
            \bigl|\mu^{\Pf}(t)_i\bigr|(\Rz^N) \leq \|g(0)\|_{L^1} + \int^{T_0}_{0} \|\nabla d(s)_i\|_{L^1} \der s.
        \end{equation}
        \item $\mu^{\Pf}(t)_i \in L^{\infty}(\Rz^N)$ and
        \begin{equation}\label{eq:LIFASsupnorm}
            \|\mu^{\Pf}(t)_i\|_{\infty} \leq \|g(0)_i\|_{\infty} + \int_{0}^{t} \|\nabla d(s)_i\|_{\infty} \der s.
        \end{equation}
        \item for all $0<q<1$, $\mu^{\Pf}(t)_i$ has finite $q$--moment and 
        \begin{equation}\label{eq:qmombound}
            \Mc_q\bigl(\mu^{\Pf}(t)_i\bigr) \leq \Mc_q\bigl(g(0)_i\bigr) + \int_0^{T_0} \Mc_q\bigl(\nabla d(s)_i\bigr) \der s + K T_0,
        \end{equation}
        where
        \begin{equation}
            \begin{aligned}
                K &\coloneqq \biggl(\|g\|_{L^{\infty}([0,T_0] \times \Rz^N)}\bigl(1 + C(N,q)\bigr)\\
                &\qquad + \eta^2 q\,2^{1-q}\bigl(D(N,q)(1 + \eta^{N+q-2}) + 1\bigr)\biggr)\\
                &\qquad \cdot \biggl(\|g(0)\|_{L^1} + \|g(0)_i\|_{\infty} + \int_0^{T_0} \|\nabla d(s)_i\|_{L^1} + \|\nabla d(s)_i\|_{\infty} \der s\biggr),
            \end{aligned}
        \end{equation}
        with $C(N,q)$ and $D(N,q)$ the constants of Lemmas~\ref{le:eflowop} and \ref{le:fflowop}.
    \end{itemize}
    In particular, for any sequence of partitions $\Pf_n$ of $[0,T_0]$ such that $\lim_{n \to \infty} \pf(\Pf_n) = 0$, for all $i \in \llbracket 0, N-1 \rrbracket$ and all $t \in [0,T_0]$, the sequence $\bigl(\mu^{\Pf_n}(t)_i\bigr)_{n \in \Nz}$ is tight and uniformly bounded.
\end{lem}
\begin{proof}
    Let $(x_l)_{l \in \llbracket 0, M \rrbracket}$ be the points of $\Pf$ in increasing order. We proceed by induction over $l$ and prove the three estimates simultaneously. All three estimates hold for $t = 0$ since $g_0 \in \Sc(\Rz^N, \Rz^N)$. Assume they hold at $t = x_l$. Now, let $t \in [x_l, x_{l+1}]$. 

    \begin{itemize}

    \item Mass. By the first points of Lemmas~\ref{le:fflowop}, \ref{le:eflowop} and \ref{le:sop}, for all $r < x_{l+1} - x_l$ the mass does not increase under $f^{\Ph(x_l,g)}_{r}$, is preserved by $e^{\Ph(x_l,g)}_{r}$, and increases by at most $\int_{x_l}^{x_{l+1}}\|\nabla d(s)_i\|_{L^1} \der s$ under $s^{d,x_l}_{r,i}$. Summing over the steps yields the first point.

    \item Sup--norm. By the second point of Lemma~\ref{le:fflowop}, the third point of Lemma~\ref{le:eflowop} and the third point of Lemma~\ref{le:sop}, we obtain \eqref{eq:LIFASsupnorm} at time $t$ in the same way.

    \item $q$--moments. By the fifth point of Lemma~\ref{le:fflowop} (which applies since $N \geq 2$ and $x_{l+1} - x_l < 1$), the sixth point of Lemma~\ref{le:eflowop} with $\|\ol{\Ph(x_l,g)}\|_{\infty} = \|g(x_l)\|_{\infty}$, and the fifth point of Lemma~\ref{le:sop}, together with the mass bound and \eqref{eq:LIFASsupnorm} already established at the intermediate stages, the $q$--moment increases by at most $K(x_{l+1} - x_l) + \int_{x_l}^{x_{l+1}} \Mc_q(\nabla d(s)_i) \der s$ between $x_l$ and $t$. Summing over the steps yields \eqref{eq:qmombound}, since $\sum_{k < l}(x_{k+1} - x_k) \leq T_0$. Note that all the quantities appearing in $K$ are finite by Proposition~\ref{pr:regclasssol}.

    \end{itemize}

    Finally, tightness follows from the uniform $q$--moment bound and the first point of Lemma~\ref{le:q-tight}.
\end{proof}
Obtaining bounds of the $\Cc^{\al}$--norms of the LIFAS for $\al \in (0,1]$ is slightly more technical.
\begin{lem}\label{th:Linfholder}
    Let $T_0 \in [0,T)$, $\al \in (0,1]$ and let $\Pf$ be a partition of $[0,T_0]$. Then, for all $i \in \llbracket 0, N-1 \rrbracket$ and all $t \in [0,T_0]$, $\mu^{\Pf}(t)_i \in \Cc^{\al}(\Rz^N)$ and there exists $\de(T_0) > 0$ such that for all $\Pf$ satisfying $\pf(\Pf) < \de(T_0)$,
        \begin{equation}\label{eq:holderbound}
            \{\mu^{\Pf}(t)_i\}_{\Cc^{\al}} \leq \biggl(\{g(0)_i\}_{\Cc^{\al}} + \int_{0}^{t} \{\nabla d(s)_i\}_{\Cc^{\al}} \der s\biggr)e^{2\al\int_{0}^{t} \|g(s)\|_{\Lip} \der s}.
        \end{equation}
\end{lem}
\begin{proof}
    Let $(x_l)_{l \in \llbracket 0, M \rrbracket}$ be the points of $\Pf$ in increasing order. Since $s \mapsto \|g(s)\|_{\Lip}$ is continuous on $[0,T_0]$ by Proposition~\ref{pr:regclasssol}, we may choose $\de(T_0) > 0$ such that for every partition $\Pf$ of $[0,T_0]$ with $\pf(\Pf) < \de(T_0)$ and every $l \in \llbracket 0, M \rrbracket$,
    \begin{equation}\label{eq:riemann}
        \sum_{k < l} (x_{k+1} - x_k)\|g(x_k)\|_{\Lip} \leq 2 \int_0^{x_l} \|g(s)\|_{\Lip} \der s.
    \end{equation}
    Indeed, if $\int_0^{T_0}\|g(s)\|_{\Lip} \der s = 0$, then $g \equiv 0$ on $[0,T_0]$, as $g(s)$ is constant in space and integrable, and there is nothing to prove. Otherwise, the left--hand side is a Riemann sum which converges to the integral uniformly in $l$ as $\pf(\Pf) \to 0$. Now, let $t \in [x_l, x_{l+1}]$. 
    
    As in the proof of Lemma~\ref{le:massmom}, $\mu^{\Pf}(t)_i$ is obtained from $\mu^{\Pf}(x_l)_i$ by applying $f^{\Ph(x_l,g)}_r$, then possibly $e^{\Ph(x_l,g)}_{r'}$, then possibly $s^{d,x_l}_{r'',i}$, with $r, r', r'' \in [0, x_{l+1} - x_l]$. By the fourth points of Lemmas~\ref{le:fflowop}, \ref{le:eflowop} and \ref{le:sop}, using $\|\ol{\Ph(x_l,g)}\|_{\Lip} = \|g(x_l)\|_{\Lip}$,
    \begin{equation}\label{eq:holderstep}
        \{\mu^{\Pf}(t)_i\}_{\Cc^{\al}} \leq \{\mu^{\Pf}(x_l)_i\}_{\Cc^{\al}} e^{\al (x_{l+1} - x_l)\|g(x_l)\|_{\Lip}} + \int_{x_l}^{t} \{\nabla d(s)_i\}_{\Cc^{\al}} \der s.
    \end{equation}
    Setting $a_l \coloneqq \{\mu^{\Pf}(x_l)_i\}_{\Cc^{\al}}$ and $b_l \coloneqq \int_{x_l}^{x_{l+1}} \{\nabla d(s)_i\}_{\Cc^{\al}} \der s$, \eqref{eq:holderstep} at $t = x_{l+1}$ reads $a_{l+1} \leq a_l e^{\al(x_{l+1} - x_l)\|g(x_l)\|_{\Lip}} + b_l$, and the discrete Grönwall inequality yields by induction
    \begin{equation}
        \begin{aligned}
            a_l &\leq \biggl(a_0 + \sum_{k < l} b_k\biggr) \prod_{k < l} e^{\al(x_{k+1} - x_k)\|g(x_k)\|_{\Lip}}\\ &\leq \biggl(\{g(0)_i\}_{\Cc^{\al}} + \int_0^{x_l}\{\nabla d(s)_i\}_{\Cc^{\al}} \der s\biggr)e^{2\al \int_0^{x_l} \|g(s)\|_{\Lip} \der s},
        \end{aligned}
    \end{equation}
    where we used \eqref{eq:riemann}. Inserting this into \eqref{eq:holderstep} and using once more \eqref{eq:riemann} for the last partial step, we obtain \eqref{eq:holderbound}.
\end{proof}
Now, we are ready to prove that the LIFAS converge to $g$.
\begin{theorem}\label{th:LIFAS}
    Let $T>0$, $N \in \Nz$, $N \geq 2$ and $\eta > 0$. Let
    \begin{itemize} 
        \item $g_0 \in \Sc(\Rz^N, \Rz^N)$ with $\div g_0 = 0$
        \item $(g,d)$ be a classical solution of \eqref{eq:NS} on $[0,T) \times \Rz^N$ with initial condition $g_0$ satisfying the assumptions of Proposition~\ref{pr:regclasssol}
        \item $T' < T$ and $(\Pf_n)_{n \in \Nz}$ be a sequence of partitions of $[0,T']$ which satisfies $\lim_{n \to \infty} \pf(\Pf_n) = 0$. 
    \end{itemize}
    Then the (LIFAS) $\mu^{\Pf_n}$ with initial value $g_0$ take values in $\Cc^1(\Rz^N, \Rz^N)$ and converge uniformly to $g$ on $[0,T']$ with respect to the $\Wz_1$--norm, i.e., for all $i \in \llbracket 0, N-1 \rrbracket$,
    \begin{equation}
        \lim_{n \to \infty} \sup_{t \in [0,T']} \bigl\| \mu^{\Pf_n}(t)_i - g(t)_i \bigr\|_{\Wz_1} = 0.
    \end{equation}
\end{theorem}
\begin{proof}
    Fix $i \in \llbracket 0, N-1 \rrbracket$. Throughout the proof, $\Pf$ denotes a partition of $[0,T']$ with points $(x_l)_{l \in \llbracket 0, M \rrbracket}$ in increasing order and $\pf(\Pf) \leq 1$, and $\mu^{\Pf}$ the associated LIFAS. To alleviate notation, we shall write $\Ph_l \coloneqq \Ph(x_l, g)$ and $\tau_l \coloneqq x_{l+1} - x_l$ for all $l \in \llbracket 0, M-1 \rrbracket$. 
    
    By Proposition~\ref{pr:regclasssol}, $g \in \Cc\bigl([0,T'], \Cc^2(\Rz^N)\bigr) \cap \Cc\bigl([0,T'], W^{1,1}(\Rz^N)\bigr)$ and $\nabla d \in \Cc\bigl([0,T'], W^{1,1}(\Rz^N) \cap W^{1,\infty}(\Rz^N)\bigr)$.
    Furthermore, by definition $g(x_l) \in \Xc^{2,\infty}(\Rz^N)$ is divergence--free, such that $\Ph_l$ is a well--defined incompressible flow vector with $\ol{\Ph_l} = - g(x_l)$ by Corollary~\ref{co:avPhi}.
    
    We denote by $\omega_g$ a common modulus of continuity of the maps $t \mapsto g(t) \in \Cc^1_b(\Rz^N, \Rz^N)$ and $t \mapsto \|g(t)\|_{\Lip}$ on $[0,T']$, and define for all $\de >0$
    \begin{equation}\label{eq:omegad}
        \omega_d(\de) \coloneqq \sup\biggl\{ \int_A \|\nabla d(r)_i\|_{L^1} \der r \;\Big|\; A \subset [0,T'] \text{ measurable}, |A| \leq \de \biggr\},
    \end{equation}
    which satisfies $\lim_{\de \to 0} \omega_d(\de) = 0$ by absolute continuity of the Lebesgue integral. The proof consists of five steps.
    
    \bigskip

    \textbf{Step 1: Regularity and uniform bounds.}

    \bigskip

    Let $l \in \llbracket 0,M \rrbracket$. Then, Lemmas~\ref{le:eflowop}, \ref{le:fflowop} and \ref{le:sop} yield $\mu^{\Pf}(t)_i \in \Cc^1(\Rz^N) \cap W^{1,\infty}(\Rz^N)$ with $\nabla \mu^{\Pf}(t)_i \in L^1(\Rz^N)$, and 
    \begin{equation}
        \|\nabla \mu^{\Pf}(t)_i\|_{L^1} \leq \|\nabla \mu^{\Pf}(x_l)_i\|_{L^1}e^{\tau_l \|g(x_l)\|_{\Lip}} + \int_{x_l}^{x_{l+1}} \bigl\|\nabla \bigl(\nabla d(s)_i\bigr)\bigr\|_{L^1} \der s.
    \end{equation}
    By induction over $l$, we obtain 
    \begin{equation}\label{eq:gradL1bound}
        \begin{aligned}
            &\|\nabla \mu^{\Pf}(t)_i\|_{L^1}\\
            &\qquad \leq K_1 \coloneqq \biggl(\|\nabla g(0)_i\|_{L^1} + \int_0^{T'} \bigl\|\nabla \bigl(\nabla d(s)_i\bigr)\bigr\|_{L^1} \der s\biggr) e^{T' \sup_{s \in [0,T']}\|g(s)\|_{\Lip}},
        \end{aligned}
    \end{equation}
    since $\sum_{k < l} \tau_k \|g(x_k)\|_{\Lip} \leq T' \sup_{s \in [0,T']} \|g(s)\|_{\Lip}$. Furthermore, by Lemma~\ref{le:massmom}, for all $t \in [0,T']$,
    \begin{equation}\label{eq:massbound}
        \bigl|\mu^{\Pf}(t)_i\bigr|(\Rz^N) \leq K_2 \coloneqq \|g(0)\|_{L^1} + \int_0^{T'} \|\nabla d(s)_i\|_{L^1} \der s,
    \end{equation}
    and, for any fixed $q \in (0,1)$, $\sup_{t \in [0,T']} \Mc_q\bigl(\mu^{\Pf}(t)_i\bigr) \leq K_3$ for a constant $K_3$ independent of $\Pf$.

    \bigskip

    \textbf{Step 2: Equicontinuity in $\Wz_1(\Rz^N)$.}

    \bigskip

    Next, we show that the LIFAS are equicontinuous on $[0,T']$. Let $l \in \llbracket 0, M-1 \rrbracket$ and $t,s \in [x_l, x_{l+1}]$, such that $t$ and $s$ lie in the same of the three subintervals $I_1 \coloneqq \bigl[x_l, \frac{x_{l+1} + 2x_l}{3}\bigr]$, $I_2 \coloneqq \bigl[\frac{x_{l+1} + 2x_l}{3}, \frac{2x_{l+1} + x_l}{3}\bigr]$, $I_3 \coloneqq \bigl[\frac{2x_{l+1} + x_l}{3}, x_{l+1}\bigr]$. 
    Then, Definition~\ref{de:LIFAS} together with
    \begin{itemize} 
        \item  \eqref{eq:gradL1bound} and the third point of Lemma~\ref{le:fflowop} applied to $\mu^{\Pf}(x_l)_i$,
        \item  the second point of Lemma~\ref{le:eflowop} applied to $f^{\Ph_l}_{\tau_l}(\mu^{\Pf}(x_l)_i)$, whose mass is bounded by $K_2$,
        \item  the second point of Lemma~\ref{le:sop},
    \end{itemize}
    yield respectively
    \begin{equation}
        \bigl\|\mu^{\Pf}(t)_i - \mu^{\Pf}(s)_i\bigr\|_{\Wz_1} \leq 
        \begin{cases}
            \frac{9}{2}\eta^2 K_1 |t-s| & \text{for } s,t \in I_1,\\
            3 \|g\|_{L^{\infty}([0,T'] \times \Rz^N)} K_2 |t-s| & \text{for } s,t \in I_2, \\
            \omega_d(3|t-s|) & \text{for } s,t \in I_3.
        \end{cases}
    \end{equation}
    For arbitrary $s < t$ in $[0,T']$, we insert all points of $\Pf$ and all endpoints of the above subintervals lying in $[s,t]$ and apply the triangle inequality. The pieces lying in first subintervals have total length at most $|t-s|$, those in second subintervals as well, and the pieces lying in third subintervals correspond to disjoint subsets of $[0,T']$ of total measure at most $3|t-s|$. By \eqref{eq:omegad}, we conclude that for all $s,t \in [0,T']$,
    \begin{equation}\label{eq:equicont}
        \bigl\|\mu^{\Pf}(t)_i - \mu^{\Pf}(s)_i\bigr\|_{\Wz_1} \leq K_4|t-s| + \omega_d(3|t-s|),
    \end{equation}
    where $K_4 \coloneqq \frac{9}{2}\eta^2 K_1 + 3\|g\|_{L^{\infty}([0,T'] \times \Rz^N)} K_2$ does not depend on $\Pf$.

    \bigskip

    \textbf{Step 3: The LIFAS satisfy the Fokker--Planck \eqref{eq:FP} equation up to a vanishing error.}

    \bigskip

    In this step, we show that the LIFAS satisfy the linearised Navier--Stokes equations in the limit. Let $\ph \in \Cc^{\infty}_c(\Rz^N)$ and $t \in [0,T']$. 
    Let $l \in \llbracket 0, M-1 \rrbracket$. We decompose
    \begin{equation}\label{eq:threesteps}
        \begin{aligned}
            &\int_{\Rz^N} \ph \der \bigl(\mu^{\Pf}(x_{l+1})_i - \mu^{\Pf}(x_l)_i\bigr)\\
            &\qquad = \int_{\Rz^N} \ph \der \Bigl(f^{\Ph_l}_{\tau_l}\bigl(\mu^{\Pf}(x_l)_i\bigr) - \mu^{\Pf}(x_l)_i\Bigr)\\
            &\qquad\quad + \int_{\Rz^N} \ph \der \Bigl( e^{\Ph_l}_{\tau_l}\bigl(f^{\Ph_l}_{\tau_l}(\mu^{\Pf}(x_l)_i)\bigr) - f^{\Ph_l}_{\tau_l}\bigl(\mu^{\Pf}(x_l)_i\bigr)\Bigr)\\
            &\qquad\quad + \int_{x_l}^{x_{l+1}} \int_{\Rz^N} \ph\, \nabla d(r)_i \der x \der r,
        \end{aligned}
    \end{equation}
    where the last line is the exact contribution of the source operator $s^{d,x_l}_{\tau_l,i}$. For the first term, by the proof of Lemma~\ref{le:fflowop} and Taylor's formula, for all $x \in \Rz^N$ and $j \in \llbracket 0, N-1 \rrbracket$,
    \begin{equation}
        \Bigl| \ph(x + \eta\sqrt{\tau_l} e_j) + \ph(x - \eta\sqrt{\tau_l}e_j) - 2\ph(x) - \eta^2 \tau_l \frac{\der^2 \ph}{\der x_j^2}(x) \Bigr| \leq \frac{\eta^3 \tau_l^{\frac{3}{2}}}{3}\|\nabla^3 \ph\|_{\infty},
    \end{equation}
    hence, by \eqref{eq:massbound},
    \begin{equation}\label{eq:fstep}
        \begin{aligned}
            &\biggl| \int_{\Rz^N} \ph \der \Bigl(f^{\Ph_l}_{\tau_l}\bigl(\mu^{\Pf}(x_l)_i\bigr) - \mu^{\Pf}(x_l)_i\Bigr) - \frac{\eta^2 \tau_l}{2N} \int_{\Rz^N} \Delta \ph \der \mu^{\Pf}(x_l)_i \biggr|\\
            &\qquad \leq \frac{\eta^3 \tau_l^{\frac{3}{2}}}{6} \|\nabla^3 \ph\|_{\infty} K_2.
        \end{aligned}
    \end{equation}
    For the second term, note that $e^{\Ph_l}_{\tau_l} = (\rh^{-g(x_l)}_{\tau_l})_*$ by Corollary~\ref{co:avPhi}, and that for all $x \in \Rz^N$, since $\frac{\der}{\der r}\rh^{-g(x_l)}_r(x) = - g(x_l)\bigl(\rh^{-g(x_l)}_r(x)\bigr)$,
    \begin{equation}
        \begin{aligned}
            &\Bigl| \ph\bigl(\rh^{-g(x_l)}_{\tau_l}(x)\bigr) - \ph(x) + \tau_l \Lc_{g(x_l)} \ph(x) \Bigr| \\
            &\qquad \leq \frac{\tau_l^2}{2} \|\ph\|_{\Cc^2}\Bigl(\|g(x_l)\|_{\infty}^2 + \|g(x_l)\|_{\infty}\|g(x_l)\|_{\Lip}\Bigr).
        \end{aligned}
    \end{equation}
    Moreover, the function $\psi \coloneqq \Lc_{g(x_l)} \ph $ is bounded and Lipschitz satisfying $\|\psi\|_{\Lip} \leq \|\ph\|_{\Cc^2}\bigl(\|g(x_l)\|_{\infty} + \|g(x_l)\|_{\Lip}\bigr)$, so that, by definition of $f^{\Ph_l}_{\tau_l}$,
    \begin{equation}
        \biggl| \int_{\Rz^N} \psi \der \Bigl(f^{\Ph_l}_{\tau_l}\bigl(\mu^{\Pf}(x_l)_i\bigr) - \mu^{\Pf}(x_l)_i\Bigr) \biggr| \leq \eta\sqrt{\tau_l}\|\psi\|_{\Lip} K_2.
    \end{equation}
    Combining the last two estimates with \eqref{eq:massbound} and the first point of Lemma~\ref{le:fflowop}, we obtain
    \begin{equation}\label{eq:estep}
        \begin{aligned}
            &\biggl| \int_{\Rz^N} \ph \der \Bigl( e^{\Ph_l}_{\tau_l}\bigl(f^{\Ph_l}_{\tau_l}(\mu^{\Pf}(x_l)_i)\bigr) - f^{\Ph_l}_{\tau_l}\bigl(\mu^{\Pf}(x_l)_i\bigr)\Bigr)\\
            &\qquad + \tau_l \int_{\Rz^N} \Lc_{g(x_l)} \ph  \der \mu^{\Pf}(x_l)_i \biggr| \leq K_5\|\ph\|_{\Cc^2} \tau_l^{\frac{3}{2}},
        \end{aligned}
    \end{equation}
    where $K_5 \coloneqq K_2\bigl(1 + \eta\bigr)\bigl(1 + \sup_{s \in [0,T']}\|g(s)\|_{\Cc^1}\bigr)^2$. Denote by $l^{\Pf}(t)$ the maximal $l \in \llbracket 0, M \rrbracket$ with $x_l \leq t$. Summing \eqref{eq:threesteps}, \eqref{eq:fstep} and \eqref{eq:estep} over $l < l^{\Pf}(t)$ yields
    \begin{equation}\label{eq:sum}
        \begin{aligned}
            &\biggl| \int_{\Rz^N} \ph \der \bigl(\mu^{\Pf}(x_{l^{\Pf}(t)})_i - g(0)_i\bigr)\\
            &\qquad - \sum_{l < l^{\Pf}(t)} \tau_l \biggl( \frac{\eta^2}{2N} \int_{\Rz^N} \Delta \ph \der \mu^{\Pf}(x_l)_i - \int_{\Rz^N} \Lc_{g(x_l)} \ph  \der \mu^{\Pf}(x_l)_i \biggr)\\
            &\qquad - \int_0^{x_{l^{\Pf}(t)}} \int_{\Rz^N} \ph\, \nabla d(r)_i \der x \der r \biggr| \leq T' \Bigl(\frac{\eta^3}{6}\|\nabla^3 \ph\|_{\infty}K_2 + \|\ph\|_{\Cc^2}K_5\Bigr)\pf(\Pf)^{\frac{1}{2}}.
        \end{aligned}
    \end{equation}
    Now, for $r \in [x_l, x_{l+1}]$, since $\Delta \ph$ and $\Lc_{g(x_l)} \ph$ are bounded and Lipschitz, the definition of the $\Wz_1$--norm and \eqref{eq:equicont}, \eqref{eq:massbound} give
    \begin{equation}\label{eq:replace}
        \begin{aligned}
            &\biggl| \int_{\Rz^N} \Delta \ph \der \bigl(\mu^{\Pf}(r)_i - \mu^{\Pf}(x_l)_i\bigr) \biggr| \leq \|\ph\|_{\Cc^3}\bigl(K_4 \pf(\Pf) + \omega_d(3\pf(\Pf))\bigr),\\
            &\biggl| \int_{\Rz^N} \Lc_{g(r)} \ph  \der \mu^{\Pf}(r)_i - \int_{\Rz^N} \Lc_{g(x_l)} \ph \der \mu^{\Pf}(x_l)_i \biggr|\\
            &\qquad \leq \|\ph\|_{\Cc^1}\omega_g\bigl(\pf(\Pf)\bigr)K_2\\
            &\qquad + \|\ph\|_{\Cc^2}\Bigl(1 + \sup_{s \in [0,T']}\|g(s)\|_{\Cc^1}\Bigr)\bigl(K_4\pf(\Pf) + \omega_d(3\pf(\Pf))\bigr).
        \end{aligned}
    \end{equation}
    Hence, the sum in \eqref{eq:sum} differs from 
    $$
    \int_0^{x_{l^{\Pf}(t)}}\biggl(\frac{\eta^2}{2N}\int \Delta \ph \der \mu^{\Pf}(r)_i - \int \Lc_{g(r)} \ph  \der \mu^{\Pf}(r)_i\biggr)\der r
    $$ 
    by at most $T'$ times the sum of the right--hand sides above.
    By \eqref{eq:equicont}, $$\bigl|\int \ph \der (\mu^{\Pf}(t)_i - \mu^{\Pf}(x_{l^{\Pf}(t)})_i)\bigr| \leq \|\ph\|_{\Cc^1}\bigl(K_4\pf(\Pf) + \omega_d(3\pf(\Pf))\bigr),$$ and the integrand 
    $$
    \frac{\eta^2}{2N}\int_{\Rz^N} \Delta \ph \der \mu^{\Pf}(r)_i - \int_{\Rz^N} \Lc_{g(r)} \ph  \der \mu^{\Pf}(r)_i + \int_{\Rz^N} \ph\, \nabla d(r)_i \der x 
    $$ 
    is bounded by $\|\ph\|_{\Cc^2}\bigl(\frac{\eta^2}{2N} + \sup_s \|g(s)\|_{\infty}\bigr)K_2 + \|\ph\|_{\infty}\|\nabla d(r)_i\|_{L^1}$, so that its integral over $[x_{l^{\Pf}(t)}, t]$ is bounded by a constant times $\pf(\Pf) + \omega_d(\pf(\Pf))$. 
    
    Finally, summarizing the above, there exists a nondecreasing function $\Omega : (0,1] \to \Rz_+$ with $\lim_{\de \to 0}\Omega(\de) = 0$, depending only on $\ph$, $g$, $d$, $T'$, $\eta$ and $N$, such that
    \begin{equation}\label{eq:consistency}
        \begin{aligned}
            &\biggl| \int_{\Rz^N} \ph \der \mu^{\Pf}(t)_i - \int_{\Rz^N} \ph\, g(0)_i \der x - \int_0^t \biggl( \frac{\eta^2}{2N}\int_{\Rz^N} \Delta \ph \der \mu^{\Pf}(r)_i\\
            &\qquad - \int_{\Rz^N} \Lc_{g(r)} \ph  \der \mu^{\Pf}(r)_i + \int_{\Rz^N} \ph\, \nabla d(r)_i \der x \biggr) \der r \biggr| \leq \Omega\bigl(\pf(\Pf)\bigr).
        \end{aligned}
    \end{equation}

    \bigskip

    \textbf{Step 4: Compactness.}

    \bigskip

    Let $(\Pf_n)_{n \in \Nz}$ be as in the statement; without loss of generality, $\pf(\Pf_n) \leq 1$ for all $n$. 
    
    By Step 1, for all $t \in [0,T']$, the family $\bigl(\mu^{\Pf_n}(t)_i\bigr)_{n \in \Nz}$ has uniformly bounded mass and uniformly bounded $q$--moments. Hence, by Lemma~\ref{le:q-tight}, it is tight in the $q'$--th moment for all $q' < q$, and, by Corollary~\ref{co:compacttight}, relatively compact in $\Wz_1(\Rz^N)$. By Step 2, the curves $t \mapsto \mu^{\Pf_n}(t)_i$ are equicontinuous in $\Wz_1(\Rz^N)$. Therefore, the Theorem of Arzel\`a--Ascoli \cite[Chapter 7, Theorem 17]{Kelley1955} yields a subsequence $\bigl(\mu^{\Pf_{n_k}}\bigr)_{k \in \Nz}$ and a continuous curve $\mu_i : [0,T'] \to \Wz_1(\Rz^N)$ with $\mu_i(0) = g(0)_i$, such that
    \begin{equation}\label{eq:unifconv}
        \lim_{k \to \infty} \sup_{t \in [0,T']} \bigl\|\mu^{\Pf_{n_k}}(t)_i - \mu_i(t)\bigr\|_{\Wz_1} = 0.
    \end{equation}
    Moreover, the mass and moment bounds of Step 1 pass to the limit. 
    
    Indeed, let $f \in \Cc_c(\Rz^N)$ with $\|f\|_{\infty} \leq 1$, and let $f_{\de} \coloneqq f * \ph_{\de}$ with $\ph_{\de}$ a standard mollifier with support in $\bigl[-\frac{\de}{2}, \frac{\de}{2}\bigr]$. Then $f_{\de} \in \Cc_b(\Rz^N) \cap \Lip(\Rz^N)$, $\|f_{\de}\|_{\infty} \leq 1$ and $f_{\de} \to f$ uniformly. By Lemma~\ref{le:weakerconv} and \eqref{eq:massbound}, $\int f_{\de} \der \mu_i(t) = \lim_k \int f_{\de} \der \mu^{\Pf_{n_k}}(t)_i \leq K_2$, and $\int f \der \mu_i(t) = \lim_{\de \to 0}\int f_{\de} \der \mu_i(t)$ since $\mu_i(t)$ is a finite measure. Taking the supremum over $f$ gives $|\mu_i(t)|(\Rz^N) \leq K_2$. 
    
    In the same way, for $f \in \Cc_c(\Rz^N)$ with $|f| \leq \min(\|x\|^q, R)$, $R > 0$, we have $|f_{\de}| \leq \min(\|x\|^q, R) + \de^q$, hence $\int f_{\de} \der \mu^{\Pf_{n_k}}(t)_i \leq K_3 + \de^q K_2$, and letting $k \to \infty$, $\de \to 0$, taking the supremum over $f$ and finally $R \to \infty$ (monotone convergence), we obtain, for all $t \in [0,T']$, $|\mu_i(t)|(\Rz^N) \leq K_2$ and $\Mc_q\bigl(\mu_i(t)\bigr) \leq K_3$.

    \textbf{Step 5: Identification of the limit.}
    For all $\ph \in \Cc^{\infty}_c(\Rz^N)$ and $r \in [0,T']$, the functions $\Delta \ph$ and $\Lc_{g(r)} \ph$ are bounded and Lipschitz, uniformly in $r$. Hence, by the definition of the $\Wz_1$--norm, the integrands in \eqref{eq:consistency} converge uniformly in $r \in [0,t]$ along $(\Pf_{n_k})_k$ by \eqref{eq:unifconv}, such that we may pass to the limit in \eqref{eq:consistency} and obtain for all $\ph \in \Cc^{\infty}_c(\Rz^N)$ and all $t \in [0,T']$,
    \begin{equation}\label{eq:limiteq}
        \begin{aligned}
            &\int_{\Rz^N} \ph \der \mu_i(t) - \int_{\Rz^N} \ph\, g(0)_i \der x\\
            &\qquad = \int_0^t \biggl( \int_{\Rz^N} \Bigl(\frac{\eta^2}{2N}\Delta \ph - \Lc_{g(r)} \ph\Bigr) \der \mu_i(r) + \int_{\Rz^N} \ph\, \nabla d(r)_i \der x \biggr) \der r.
        \end{aligned}
    \end{equation}
    On the other hand, \eqref{eq:divlieprop} yields that $g_i$ itself satisfies \eqref{eq:limiteq}. Hence, $w(t) \coloneqq \mu_i(t) - g(t)_i \der x \in \Wz_1(\Rz^N)$ satisfies $w(0) = 0$, $\sup_{t \in [0,T']}|w(t)|(\Rz^N) < \infty$ by Step 4. Furthermore, $t \mapsto w(t)$ is continuous in $\Wz_1(\Rz^N)$ by Step 4 and since $g \in \Cc([0,T'], L^1(\Rz^N))$ by Proposition~\ref{pr:regclasssol}, and for all $\ph \in \Cc^{\infty}_c(\Rz^N)$ and $t \in [0,T']$,
    \begin{equation}\label{eq:homog}
        \int_{\Rz^N} \ph \der w(t) = \int_0^t \int_{\Rz^N} \Bigl(\frac{\eta^2}{2N}\Delta \ph - \Lc_{g(r)}\ph\Bigr) \der w(r) \der r.
    \end{equation}
    In other words, $w$ is a measure--valued solution of the Fokker--Planck equation \eqref{eq:FP} without source term and with vanishing initial condition. Recall that its drift $-g$ is bounded and Lipschitz, uniformly in time. Uniqueness of measure--valued solutions of such equations is classical, see \cite[Chapter 9]{BKRS2015} and \cite[Section 2]{Figalli2008}. Our class of solutions is slightly different, i.e. signed measures, continuity with respect to $\Wz_1$, but the proof is analogous. 

    \textbf{Conclusion.}
    We have shown that every subsequence of $\bigl(\mu^{\Pf_n}_i\bigr)_{n \in \Nz}$ admits a further subsequence converging uniformly on $[0,T']$ in $\Wz_1(\Rz^N)$ to the same limit $g_i$. Therefore, the whole sequence converges uniformly to $g_i$, which concludes the proof.
\end{proof}

\subsubsection{Estimates for solutions of the Navier--Stokes equations}

In this section, $(g,d)$ denotes a classical solution of \eqref{eq:NS} on $[0,T)$ with initial condition $g_0 \in \Sc(\Rz^N, \Rz^N)$, $\div g_0 = 0$, as in Proposition~\ref{pr:regclasssol}, and $\|\cdot\|_{\Cc^{\al}} = \|\cdot\|_{\infty} + \{\cdot\}_{\Cc^{\al}}$ denotes the full Hölder norm.

\begin{co}\label{co:solestim}
    For all $i \in \llbracket 0, N-1 \rrbracket$, all $\al \in (0,1]$ and all $t \in [0,T)$:
    \begin{itemize}
        \item the following estimate holds:
        \begin{equation}
            \|g(t)_i\|_{\infty} \leq \|g(0)_i\|_{\infty} + \int_{0}^{t} \|\nabla d(s)_i\|_{\infty} \der s.
        \end{equation}
        \item the following estimate holds:
        \begin{equation}
            \{g(t)_i\}_{\Cc^{\al}} \leq \biggl(\{g(0)_i\}_{\Cc^{\al}} + 2\int_{0}^{t} \{\nabla d(s)_i\}_{\Cc^{\al}} \der s\biggr)e^{2\al\int_{0}^{t} \|g(s)\|_{\Lip} \der s}.
        \end{equation}
        \item in particular,
        \begin{equation}
            \|g(t)_i\|_{\Lip} \leq \biggl(\|g(0)_i\|_{\Lip} + 2\int_{0}^{t} \|\nabla d(s)_i\|_{\Lip} \der s\biggr)e^{2\int_{0}^{t} \|g(s)\|_{\Lip} \der s}.
        \end{equation}
    \end{itemize}
\end{co}
\begin{proof}
    Let $T_0 \in [0, T)$ and, for all $m \in \Nz$, let $\Pf_m$ be a partition of $[0, T_0]$ with $\pf(\Pf_m) < \min\bigl(\frac{1}{m}, \de(T_0)\bigr)$, where $\de(T_0)$ is the constant of Lemma~\ref{th:Linfholder}. By Theorem~\ref{th:LIFAS}, $\mu^{\Pf_m}(t)_i \to g(t)_i$ in $\Wz_1(\Rz^N)$ as $m \to \infty$, for all $t \in [0,T_0]$, and by Lemma~\ref{th:Linfholder} the $\mu^{\Pf_m}(t)_i$ satisfy the two claimed bounds uniformly in $m$. Proposition~2.5 yields the first two points for all $t \in [0,T_0]$. The third point is the second with $\al = 1$. Since $T_0 < T$ was arbitrary, the claim follows.
\end{proof}

The next corollary has no analogue in the probability setting of \cite{Schichl2025}, where the mass of the solution is conserved: the source term $\nabla d$ creates and destroys mass, and the $L^1$--norm and the moments of the solution are controlled by those of the source.

\begin{co}\label{co:L1estim}
    For all $i \in \llbracket 0, N-1 \rrbracket$, all $q \in (0,1)$ and all $t \in [0,T)$,
    \begin{itemize}
        \item the following estimate holds: 
        \begin{equation}
            \|g(t)_i\|_{L^1} \leq \|g(0)\|_{L^1} + \int_0^t \|\nabla d(s)_i\|_{L^1} \der s,
        \end{equation}
        \item the following estimate holds: 
        \begin{equation}
            \Mc_q\bigl(g(t)_i\bigr) \leq \Mc_q\bigl(g(0)_i\bigr) + \int_0^t \Mc_q\bigl(\nabla d(s)_i\bigr) \der s + K t,
        \end{equation}
    \end{itemize}
    where $K$ is the constant of Lemma~\ref{le:massmom}.
\end{co}
\begin{proof}
    Let $T_0 \in [0, T)$ and, for all $m \in \Nz$, let $\Pf_m$ be a partition of $[0, T_0]$ with $\pf(\Pf_m) < \min\bigl(\frac{1}{m}, \de(T_0)\bigr)$, where $\de(T_0)$ is the constant of Lemma~\ref{th:Linfholder}. 
    By Lemma~\ref{le:massmom}, the $\mu^{\Pf_m}(t)_i$ satisfy the two bounds uniformly in $m$, with $T_0$ in place of $t$ in the integrals. 
    Since $\mu^{\Pf_m}(t)_i \to g(t)_i$ in $\Wz_1(\Rz^N)$ with uniformly bounded masses, the mass and the $q$--moment of the limit are bounded by the limits inferior of those of the sequence, analogously to Step 4 of the proof of Theorem~\ref{th:LIFAS}. Thus, the result follows.
\end{proof}

\begin{re}\label{re:riccati}
    Corollary~\ref{co:solestim} gives a bound on the Lipschitz norm of $g(t)$ in terms of $\int_0^t\|\nabla^2 d(s)\|_{\infty} \der s$, and with \eqref{eq:pressure}, one may thus hope to close the estimate in the form $\|g(t)\|_{\Lip} \leq \|g(0)\|_{\Lip} + C(N)\int_0^t\|g(s)\|_{\Lip}^2 \der s$, which, by the generalized Grönwall inequality \cite{Bihari1956}, would yield $\|g(t)\|_{\Lip} \leq \|g(0)\|_{\Lip}\bigl(1 - C(N)\|g(0)\|_{\Lip}t\bigr)^{-1}$ on $\bigl[0, \min(T, \frac{1}{C(N)\|g(0)\|_{\Lip}})\bigr)$. However, his would require an estimate $\|\nabla^2 d(s)\|_{\infty} \leq C(N)\|\nabla g(s)\|_{\infty}^2$, i.e., the boundedness of the Riesz transforms $R_jR_k$ on $L^{\infty}(\Rz^N)$, which is false, see Lemma~\ref{le:dievwd} and \cite[Chapter II, \S 4]{Stein1970}. A time--horizon depending only on the Lipschitz norm of the initial condition, uniformly in $\eta$, would moreover contradict the strong ill--posedness of the Euler equations in $\Cc^1$ \cite{BourgainLi2015}. Since the Riesz transforms are bounded on Hölder spaces, the Riccati argument does close in (weighted) Hölder norms. This is carried out for the non--linear scheme in Proposition~\ref{pr:IFASbounds} below.
\end{re}

\subsubsection{Incompressible flow approximation series}

Our second main result establishes the convergence of the non--linear scheme (IFAS) to the unique classical solution of \eqref{eq:NS} on an explicit time interval. In contrast to Theorem~\ref{th:LIFAS}, the a priori bounds on the approximation series cannot be obtained from the regularity of a given solution, but have to be derived from the operators alone. Two obstructions have to be overcome: the dievergence operator maps neither $L^1(\Rz^N, \Rz^N)$ nor $W^{2,\infty}(\Rz^N,\Rz^N)$ into itself (see Lemma~\ref{le:dievwd}), whereas Definition~\ref{de:IFAS} requires $\mu^{\Pf}(x_l)$ to be a divergence--free element of $\Xc^{2,\infty}(\Rz^N) \cap \Wz_1(\Rz^N)^N$ at every partition point, and, as explained in Remark~\ref{re:riccati}, the estimates cannot be closed in the Lipschitz norm. 

Both obstructions are removed by working in weighted Hölder spaces introduced in Section~\ref{sbsec:weightedholder}. Those encode the optimal spatial decay of Proposition~\ref{pr:regclasssol}. 
With the estimates for the flow operators established in Section~\ref{sbsec:flowestims}, we can derive the necessary a priori bounds for the IFAS to mimic the proof of Theorem~\ref{th:LIFAS}. 

To keep notation short, for a partition $\Pf$ of $[0,T]$ with points $(x_l)_{l}$ in increasing order, we write, as in Definition~\ref{de:IFAS},
\begin{equation}\label{eq:nul}
    \nu^{\Pf}_l \coloneqq \biggl[\biggl( e^{\Phi\bigl(\mu^{\Pf}(x_l)\bigr)}_{x_{l+1} - x_l}\circ f^{\Phi\bigl(\mu^{\Pf}(x_l)\bigr)}_{x_{l+1} - x_l}\biggr)\bigl(\mu^{\Pf}(x_l)_i\bigr)\biggr]_{i \in \llbracket 0, N-1 \rrbracket},
\end{equation}
such that $\mu^{\Pf}(x_{l+1}) = s\bigl(\nu^{\Pf}_l\bigr)$. We set $K \coloneqq K_0(1 + K_1)$ where $K_0$ is the maximum of the constants of Lemmas~\ref{le:eflowop} and \ref{le:fflowop} and $K_1$ is the constant of Lemma~\ref{le:sweighted}. We have seen that $K$ depends only on $N$, $\al$ and $\ep$.

\begin{pr}[Time--horizon and a priori bounds for IFAS]\label{pr:IFASbounds}
    Let $N \geq 2$, and $\al, \ep \in (0,1)$, $\eta \geq 0$ and $g_0 \in \mathcal E_{\al,\ep}$, $g_0 \neq 0$. Set 
    \begin{equation}\label{eq:Tstar}
        \begin{aligned}
            &T^*(g_0, 0) \coloneqq \frac{1}{2K\|g_0\|_{\mathcal E}},\\
            &T^*(g_0, \eta) \coloneqq \frac{1}{K\eta^2}\log\Bigl(1 + \frac{\eta^2}{2\|g_0\|_{\mathcal E}}\Bigr)
        \end{aligned}
    \end{equation}
    for $\eta >0$. Then $T^* \geq \frac{1}{K(2\|g_0\|_{\mathcal E} + \eta^2)}$, $T^*$ is continuous and decreasing in $\eta \in [0,\infty)$. 
    
    Moreover, set $\pf_0 \coloneqq \min\bigl(1, \eta^{-2}, (2\|g_0\|_{\mathcal E} + \eta^2)^{-1}\bigr)$ with the convention $0^{-2} \coloneqq \infty$. Then, for every $T > 0$ and every partition $\Pf$ of $[0,T]$ with $\pf(\Pf) \leq \pf_0$, the (IFAS) $\mu^{\Pf}$ of Definition~\ref{de:IFAS} with initial value $g_0$ is well--defined on $[0, \min(T,T^*)]$, and:
    \begin{itemize}
        \item for every partition point $x_l \leq T^*$, $\mu^{\Pf}(x_l) \in \mathcal E$ is divergence--free with mean zero and $\|\mu^{\Pf}(x_l)\|_{\mathcal E} \leq 2\|g_0\|_{\mathcal E} + \eta^2$.
        \item for every $t \in [0, \min(T,T^*)]$, $\mu^{\Pf}(t) \in \mathfrak E$ has mean zero and $\|\mu^{\Pf}(t)\|_{\mathcal E} \leq (1 + K)^2(2\|g_0\|_{\mathcal E} + \eta^2)$.
        \item for every $l$ with $x_{l+1} \leq T^*$, we have $\|\div \nu^{\Pf}_l\|_{\mathcal F} \leq K(1+K) (2\|g_0\|_{\mathcal E} + \eta^2)^2(x_{l+1} - x_l)$, and
        \begin{equation}\label{eq:H3new}
            \bigl\|\nabla\Delta^{-1}\div \nu^{\Pf}_l\bigr\|_{L^1} \leq C(N,\ep)K^2(1+K)(2\|g_0\|_{\mathcal E} + \eta^2)^2(x_{l+1} - x_l).
        \end{equation}
    \end{itemize}
\end{pr}
\begin{proof}
    Set $A \coloneqq 2a_0 + \eta^2$. Let $(x_l)_{l \in \llbracket 0, M \rrbracket}$ be the points of $\Pf$ indexed in increasing order and set for all $l \in \llbracket 0,M \rrbracket$, $\tau_l \coloneqq x_{l+1} - x_l \leq \pf_0$, and $a_l \coloneqq \|\mu^{\Pf}(x_l)\|_{\mathcal E}$ whenever $\mu^{\Pf}(x_l) \in \mathcal E$. In particular $a_0 \coloneqq \|g_0\|_{\Ec}$.

    First, the properties of $T^*$ are trivial for $\eta = 0$. For $\eta > 0$ they follow from $\frac{u}{1+u} \leq \log(1 + u) \leq u$ for $u = \frac{\eta^2}{2a_0}$. Indeed, $T^* \geq \frac{1}{K\eta^2}\cdot\frac{\eta^2/(2a_0)}{1 + \eta^2/(2a_0)} = \frac{1}{K(2a_0 + \eta^2)}$, and $T^* \to \frac{1}{2Ka_0}$ as $\eta \to 0$. The monotonicity in $\eta$ is the concavity of $u \mapsto \log(1 + u)$.

    We show by induction over $l$ that $\mu^{\Pf}(x_l) \in \mathcal E$ with mean zero and $a_l \leq A$ for all $l$ with $x_l \leq T^*$: 
    
    \bigskip

    \textbf{Initialisation:} 

    \noindent
    For $l = 0$ this holds by assumption. 

    \bigskip
    
    \textbf{Induction hypothesis:} 
    
    \noindent
    Suppose it holds for all $k \leq l$ and $x_{l+1} \leq T^*$. 
    
    \bigskip

    \textbf{Induction step:} 

    \noindent
    Set $h \coloneqq \mu^{\Pf}(x_l) \in \mathcal E$, such that $\Ph(h)$ is a well--defined incompressible flow vector with $\ol{\Ph(h)} = -h$. Set $F_l \coloneqq \div \nu^{\Pf}_l$. 
    Since $\eta\sqrt{\tau_l} \leq 1$ and $\tau_la_l \leq \pf_0 A \leq 1$ by assumption, Lemmas~\ref{le:eflowop} and \ref{le:fflowop} apply and yield, for all $r \in [0,\tau_l]$,
    \begin{equation}\label{eq:stepbounds}
        \begin{aligned}
            &\|f^{\Ph(h)}_r(h)\|_{\mathcal E} \leq (1 + K\eta^2\tau_l)a_l,\\ &\bigl\|e^{\Ph(h)}_r\bigl(f^{\Ph(h)}_{\tau_l}(h)\bigr)\bigr\|_{\mathcal E} \leq (1 + K\tau_la_l)(1 + K\eta^2\tau_l)a_l,\\ &\|F_l\|_{\mathcal F} \leq K_0\tau_la_l^2(1 + K\eta^2\tau_l).
        \end{aligned}
    \end{equation}
    Moreover, all these fields have mean zero, with $f^{\Ph(h)}_r(h)$ divergence--free, and $\int x_j F_l \der x = -\int f^{\Ph(h)}_{\tau_l}(h)_j \der x = 0$. Hence, Lemma~\ref{le:sweighted} applies to $\nu^{\Pf}_l$ and yields that, for all $\theta \in [0,1]$, $s_{\theta}(\nu^{\Pf}_l) \in \mathfrak E$, has mean zero and satisfies
    \begin{equation}\label{eq:riccatistep}
        \|s_{\theta}(\nu^{\Pf}_l)\|_{\mathcal E} \leq (1 + K\eta^2\tau_l)a_l\bigl(1 + K_0\tau_la_l + K_1K_0\tau_la_l\bigr) = (1 + K\eta^2\tau_l)(1 + K\tau_la_l)a_l,
    \end{equation}
    and $\mu^{\Pf}(x_{l+1}) = s(\nu^{\Pf}_l) \in \mathcal E$. 
    In particular, $a_{l+1} \leq (1 + K\eta^2\tau_l)(1 + K\tau_la_l)a_l$. 
    
    Now, set $b_k \coloneqq a_ke^{-K\eta^2x_k}$ for all $k \leq l$. Using $(1 + K\eta^2\tau_k)e^{-K\eta^2\tau_k} \leq 1$, we obtain $b_{k+1} \leq b_k(1 + K\tau_ka_k) = b_k\bigl(1 + K\tau_ke^{K\eta^2x_k}b_k\bigr)$, hence,
    \begin{equation}
        \frac{1}{b_{k+1}} \geq \frac{1}{b_k} - \frac{K\tau_ke^{K\eta^2x_k}}{1 + K\tau_ke^{K\eta^2x_k}b_k} \geq \frac{1}{b_k} - K\tau_ke^{K\eta^2x_k}.
    \end{equation}
    Now, summing over $k \leq l$ and comparing the Riemann sum of the increasing function $s \mapsto e^{K\eta^2s}$ with its integral,
    \begin{equation}
        \begin{aligned}
            \frac{1}{b_{l+1}} &\geq \frac{1}{a_0} - K\sum_{k \leq l}\tau_ke^{K\eta^2x_k} \geq \frac{1}{a_0} - K\int_0^{x_{l+1}}e^{K\eta^2s} \der s = \frac{1}{a_0} - \frac{e^{K\eta^2x_{l+1}} - 1}{\eta^2}\\ &\geq \frac{1}{a_0} - \frac{1}{2a_0} = \frac{1}{2a_0},
        \end{aligned}
    \end{equation}
    where we used $x_{l+1} \leq T^*$, i.e., $e^{K\eta^2x_{l+1}} \leq 1 + \frac{\eta^2}{2a_0}$. Therefore, $$a_{l+1} = b_{l+1}e^{K\eta^2x_{l+1}} \leq 2a_0\bigl(1 + \frac{\eta^2}{2a_0}\bigr) = A,$$ which completes the induction and proves the first point.

    For $t \in [x_l, x_{l+1}]$ with $x_{l+1} \leq T^*$, $\mu^{\Pf}(t)$ is given by one of the expressions in \eqref{eq:stepbounds} or \eqref{eq:riccatistep}, and $K\eta^2\tau_l \leq K$, $K\tau_la_l \leq K$ give the second point. 
    If $T^* \in (x_l, x_{l+1})$ for some $l$, the same bounds hold on $[x_l, T^*]$, as they only use $a_l \leq A$. 
    
    The bound on $\|F_l\|_{\mathcal F}$ is \eqref{eq:stepbounds} with $a_l \leq A$, $K_0 \leq K$ and $1 + K\eta^2\tau_l \leq 1 + K$. Finally, \eqref{eq:H3new} follows from $\nabla\Delta^{-1}\div \nu^{\Pf}_l = \nabla\Gamma*F_l$, Lemma~\ref{le:sweighted} and \eqref{eq:Etomoments}: $\|\nabla\Gamma*F_l\|_{L^1} \leq C(N,\ep)|\nabla\Gamma*F_l|_{N+1-\ep} \leq C(N,\ep)K_1\|F_l\|_{\mathcal F}$.
\end{proof}

\begin{re}
    The time $T^*$ in \eqref{eq:Tstar} depends on the initial datum only through its weighted $\Cc^{2,\al}$--norm $\|g_0\|_{\mathcal E}$, and it is bounded from below by $\frac{1}{K(2\|g_0\|_{\mathcal E} + \eta^2)}$, uniformly in the viscosity $\frac{\eta^2}{2N}$. As $\eta \to 0$, $T^*$ tends to $\frac{1}{2K\|g_0\|_{\mathcal E}}$, which is the classical lifespan of Hölder solutions of the Euler equations. This is consistent with Remark~\ref{re:riccati}. The Riccati inequality \eqref{eq:riccatistep} closes in weighted Hölder norms, since the Riesz transforms are bounded on Hölder spaces (see Lemma~\ref{le:sweighted}), while it cannot close in the Lipschitz norm. 
    The term $\eta^2$ in $2\|g_0\|_{\mathcal E} + \eta^2$ and $T^*$ comes from the weights only. The averaged flow dispersion operators do not increase any unweighted $\Cc^{k,\al}$--norm, but spread mass to infinity at rate $\eta^2 t$.
\end{re}

\begin{theorem}\label{th:IFAS}
    Let $N \in \Nz$, $N \geq 2$, $\al, \ep \in (0,1)$, $\eta \geq 0$.   Let $T^* = T^*(g_0,\eta)$ be given by \eqref{eq:Tstar}. Moreover, let 
    \begin{itemize}
        \item $g_0 \in \mathcal E_{\al,\ep}$, $g_0 \neq 0$.
        \item $(\Pf_n)_{n \in \Nz}$ be a sequence of partitions of $[0,T^*]$ with $\lim_{n \to \infty}\pf(\Pf_n) = 0$.
        \item $\mu^{\Pf_n}$ be the (IFAS) of Definition~\ref{de:IFAS} with initial value $g_0$. 
    \end{itemize}
    Then, there exists a unique classical solution $(g,d)$ of \eqref{eq:NS} on $[0,T^*]$ with initial condition $g_0$ satisfiyng $g \in \Cc\bigl([0,T^*], L^2(\Rz^N)\bigr) \cap L^{\infty}\bigl([0,T^*], W^{1,\infty}(\Rz^N)\bigr)$ and $\nabla d \in L^{\infty}\bigl([0,T^*], L^1(\Rz^N) + L^{\infty}(\Rz^N)\bigr)$, and
    \begin{equation}\label{eq:IFASconv}
        \lim_{n \to \infty}\ \sup_{t \in [0,T^*]}\ \Bigl( \sum_{i=0}^{N-1}\bigl\|\mu^{\Pf_n}(t)_i - g(t)_i\bigr\|_{\Wz_1} + \bigl\|\mu^{\Pf_n}(t) - g(t)\bigr\|_{\Cc^2_b(\Rz^N)}\Bigr) = 0.
    \end{equation}
    Moreover, let $K$ be the constant of Proposition~\ref{pr:IFASbounds}. Then, for all $t \in [0,T^*]$, $g(t) \in \mathcal E_{\al,\ep}$ with $\|g(t)\|_{\mathcal E} \leq (1 + K)^2 (2\|g_0\|_{\mathcal E} + \eta^2)$.
    In particular, for $ k \in \llbracket 0, 2 \rrbracket$ and all $(t,x) \in [0,T^*] \times \Rz^N$,
    \begin{equation}\label{eq:decaysol}
        |\nabla^k g(t,x)| \leq (1 + K)^2\bigl(2\|g_0\|_{\mathcal E} + \eta^2\bigr)\bigl(1 + \|x\|^2\bigr)^{-\frac{N + 1 + k - \ep}{2}}. 
    \end{equation}
    For $\eta = 0$, \eqref{eq:NS} is the incompressible Euler equation and the (IFAS) only consist of expectational flow and dievergence steps.
\end{theorem}
\begin{proof}
    Fix $q \coloneqq \frac{1 - \ep}{2}$. Without loss of generality, $\pf(\Pf_n) \leq \pf_0$ for all $n \in \Nz$, such that Proposition~\ref{pr:IFASbounds} applies to every $\Pf_n$.  By Proposition~\ref{pr:IFASbounds} and \eqref{eq:Etomoments}, there exists $B = B(N, \al, \ep, \|g_0\|_{\mathcal E}, \eta)$ such that for all $n\in\Nz$, all $t \in [0,T^*]$, all $i \in \llbracket 0, N-1 \rrbracket$, and all $l$ with $x_{l+1} \leq T^*$,
    \begin{equation}\label{eq:Hnew}
        \begin{aligned}
            &\bigl|\mu^{\Pf_n}(t)_i\bigr|(\Rz^N) + \Mc_q\bigl(\mu^{\Pf_n}(t)_i\bigr) + \|\mu^{\Pf_n}(t)\|_{\Cc^{2,\al}} + \bigl\|\nabla\mu^{\Pf_n}(t)_i\bigr\|_{L^1} \leq B\\ 
            &\bigl\|\nabla\Delta^{-1}\div\nu^{\Pf_n}_l\bigr\|_{L^1} \leq B(x_{l+1} - x_l),
        \end{aligned}
    \end{equation}
    where $\nu^{\Pf_n}_l$ is defined as in \eqref{eq:nul}. Throughout this proof, to simplify notation:
    \begin{itemize}
        \item $\Pf$ denotes one of the $\Pf_n$, with points $(x_l)_{l \in \llbracket 0, M \rrbracket}$ indexed in increasing order,
        \item $\Ph_l \coloneqq \Ph(\mu^{\Pf}(x_l))$ and $\nu_l \coloneqq \nu^{\Pf}_l$,
        \item we set $\tau_l \coloneqq x_{l+1} - x_l$,
        \item $A \coloneqq 2\|g_0\|_{\mathcal E} + \eta^2$ and $A' \coloneqq (1 + K)^2A$,
    \end{itemize}
    Recall $\ol{\Ph_l} = -\mu^{\Pf}(x_l)$ and that $\mu^{\Pf}(x_l) \in \mathcal E$ is divergence--free. We fix $i \in \llbracket 0, N-1 \rrbracket$ and proceed as in the proof of Theorem~\ref{th:LIFAS}.

    \bigskip

    \textbf{Step 1: Equicontinuity in $\Wz_1(\Rz^N)$.}

    \bigskip

    Let $t,s \in [x_l, x_{l+1}]$ lie in the same of the three subintervals of Definition~\ref{de:IFAS}. Then,
    \begin{itemize}
        \item the second point of Lemma~\ref{le:fflowop} and \eqref{eq:Hnew} yield in the first subinterval $\|\mu^{\Pf}(t)_i - \mu^{\Pf}(s)_i\|_{\Wz_1} \leq \frac{9}{2}\eta^2 B|t-s|$.
        \item the second point of Lemma~\ref{le:eflowop} with $\|\ol{\Ph_l}\|_{\infty} \leq B$ yields in the second subinterval $\|\mu^{\Pf}(t)_i - \mu^{\Pf}(s)_i\|_{\Wz_1} \leq 3B^2|t-s|$.
        \item the first point of Lemma~\ref{le:dievop} and \eqref{eq:Hnew} yield in the third subinterval $\|\mu^{\Pf}(t)_i - \mu^{\Pf}(s)_i\|_{\Wz_1} \leq \frac{3|t-s|}{\tau_l}\|\nabla\Delta^{-1}\div\nu_l\|_{L^1} \leq 3B|t-s|$.
    \end{itemize}
    Hence, for all $s,t \in [0,T^*]$,
    \begin{equation}\label{eq:equicontIFAS}
        \|\mu^{\Pf}(t)_i - \mu^{\Pf}(s)_i\|_{\Wz_1} \leq K_2|t-s|, \qquad K_2 \coloneqq \tfrac{9}{2}\eta^2 B + 3B^2 + 3B.
    \end{equation}

    \bigskip

    \textbf{Step 2: The IFAS satisfy the weak formulation in the limit.}

    \bigskip

    Let $\ph \in \Cc^{\infty}_c(\Rz^N, \Rz^N)$ with $\div \ph = 0$ and $t \in [0,T^*]$. Since $s_{\theta}(\nu_l) - \nu_l = -\theta\nabla\Delta^{-1}\div\nu_l$ and $\div\ph = 0$, we have $\sum_j\int_{\Rz^N}\ph_j \der\bigl(s_{\theta}(\nu_l)_j - (\nu_l)_j\bigr) = 0$ for all $\theta \in [0,1]$. Hence, as in \eqref{eq:threesteps},
    \begin{equation}
        \begin{aligned}
            &\sum_{j=0}^{N-1}\int_{\Rz^N} \ph_j \der \bigl(\mu^{\Pf}(x_{l+1})_j - \mu^{\Pf}(x_l)_j\bigr) = \sum_{j=0}^{N-1} \int_{\Rz^N} \ph_j \der \Bigl(f^{\Ph_l}_{\tau_l}\bigl(\mu^{\Pf}(x_l)_j\bigr) - \mu^{\Pf}(x_l)_j\Bigr)\\
            &\qquad + \int_{\Rz^N} \ph_j \der \Bigl( e^{\Ph_l}_{\tau_l}\bigl(f^{\Ph_l}_{\tau_l}(\mu^{\Pf}(x_l)_j)\bigr) - f^{\Ph_l}_{\tau_l}\bigl(\mu^{\Pf}(x_l)_j\bigr)\Bigr),
        \end{aligned}
    \end{equation}
    and the estimates \eqref{eq:fstep} and \eqref{eq:estep} hold with $g(x_l)$ replaced by $\mu^{\Pf}(x_l)$ and $K_2$ there replaced by $B$, so that
    \begin{equation}\label{eq:stepIFAS}
        \begin{aligned}
            &\biggl| \sum_{j=0}^{N-1}\int_{\Rz^N} \ph_j \der \bigl(\mu^{\Pf}(x_{l+1})_j - \mu^{\Pf}(x_l)_j\bigr) - \tau_l \sum_{j = 0}^{N-1}\biggl(\frac{\eta^2}{2N}\int_{\Rz^N} \Delta \ph_j\, \mu^{\Pf}(x_l)_j \der x\\
            &\qquad - \int_{\Rz^N} \nabla \ph_j \cdot \mu^{\Pf}(x_l)\, \mu^{\Pf}(x_l)_j \der x\biggr) \biggr| \leq K_3 \tau_l^{\frac{3}{2}},
        \end{aligned}
    \end{equation}
    with $K_3 \coloneqq NB\|\ph\|_{\Cc^3}(1 + \eta^3)(1 + B)^2$. Now, let $l \in \llbracket 0, M \rrbracket$ and $r \in [x_l, x_{l+1}]$. It remains to proceed like in \eqref{eq:replace} in Step 3 of the proof of Theorem~\ref{th:LIFAS}. For the non--linear term, we use Lemma~\ref{le:interpol}. Indeed, by \eqref{eq:Hnew}, \eqref{eq:equicontIFAS} and \eqref{eq:interpol} applied to $\nu = \mu^{\Pf}(r)_i - \mu^{\Pf}(x_l)_i$, which satisfies $\|\nu\|_{\Wz_1} \leq K_2\pf(\Pf)$ and $\|\nu\|_{\Lip} \leq 2B$, and since $\pf(\Pf) \leq \pf(\Pf)^{\frac{1}{N+2}}$ for $\pf(\Pf) \leq 1$, we have,
    \begin{equation}
        \|\mu^{\Pf}(r) - \mu^{\Pf}(x_l)\|_{\infty} \leq K_4\pf(\Pf)^{\frac{1}{N+2}} ,
    \end{equation}
    where $K_4 \coloneqq N\max\Bigl(\frac{2}{c_N}K_2,\ \Bigl(\frac{2(8B)^{N+1}}{c_N}K_2\Bigr)^{\frac{1}{N+2}}\Bigr)$. Therefore, since $\ph$ has compact support, the non--linear term changes by at most $2B\|\nabla\ph\|_{L^1}K_4\pf(\Pf)^{\frac{1}{N+2}}$ when $x_l$ is replaced by $r$. The diffusion term is handled as in Step 3 of the proof of Theorem~\ref{th:LIFAS} using \eqref{eq:equicontIFAS}. Summing \eqref{eq:stepIFAS} over $l < l^{\Pf}(t)$ and collecting the error terms exactly as there, we obtain a nondecreasing $\Omega : (0,1] \to \Rz_+$ with $\lim_{\de \to 0}\Omega(\de) = 0$, depending only on $\ph, B, \eta, N, T^*$, such that for all $t \in [0,T^*]$,
    \begin{equation}\label{eq:consistencyIFAS}
        \begin{aligned}
            &\biggl| \int_{\Rz^N} \ph \cdot \mu^{\Pf}(t) \der x - \int_{\Rz^N} \ph \cdot g_0 \der x\\
            &\qquad - \int_0^t \int_{\Rz^N} \biggl( \frac{\eta^2}{2N}\Delta \ph \cdot \mu^{\Pf}(r) - \sum_{j = 0}^{N-1} \mu^{\Pf}(r)_j\, \mu^{\Pf}(r)\cdot \nabla \ph_j \biggr) \der x \der r \biggr| \leq \Omega\bigl(\pf(\Pf)\bigr).
        \end{aligned}
    \end{equation}

    \bigskip

    \textbf{Step 3: Compactness and regularity of the limit.}

    \bigskip

    By \eqref{eq:Hnew}, for all $t \in [0,T^*]$ the family $(\mu^{\Pf_n}(t)_i)_n$ has uniformly bounded mass and $q$--moments, hence is relatively compact in $\Wz_1(\Rz^N)$ by Lemma~\ref{le:q-tight} and Corollary~\ref{co:compacttight}. By \eqref{eq:equicontIFAS} and the Theorem of Arzel\`a--Ascoli, after extracting successively for $i \in \llbracket 0, N-1 \rrbracket$, there exist a subsequence $(\Pf_{n_k})_k$ and $\Wz_1$--continuous curves $g_i$ with $g(0) = g_0$ such that $\mu^{\Pf_{n_k}}_i$ converges uniformly to $g_i$ on $[0,T^*]$ in $\Wz_1(\Rz^N)$. 
    
    Let us prove that the convergence is actually in $\Cc^2_b(\Rz^N)$. Fix $t \in [0,T^*]$. 
    For all $k, j \in \Nz$, the function $u_{t,k,j} \coloneqq \mu^{\Pf_{n_k}}(t)_i - \mu^{\Pf_{n_j}}(t)_i$ satisfies $\|u_{t,k,j}\|_{\Lip} \leq 2B$ and $\{\nabla^2 u_{t,k,j}\}_{\Cc^{\al}} \leq 2B$ by \eqref{eq:Hnew}. Now, $\sup_t\|u_{t,k,j}\|_{\Wz_1} \to 0$ as $k,j \to \infty$. 
    Hence \eqref{eq:interpol} gives $\sup_t\|u_{t,k,j}\|_{\infty} \to 0$, and then \eqref{eq:interpolholder} gives $\sup_t\|u_{t,k,j}\|_{\Cc^2_b} \to 0$, i.e., $(\mu^{\Pf_{n_k}}_i)_k$ is a Cauchy sequence in $\Cc([0,T^*], \Cc^2_b(\Rz^N))$.
    Thus, it converges there to some $\tilde g_i \in \Cc([0,T^*], \Cc^2_b(\Rz^N))$. One readily shows that $g_i = \tilde g_i$.
    
    Since the weighted norms of Definition~\ref{de:weightednorms} are lower semicontinuous under pointwise convergence of $\nabla^k u$, $k \leq 2$, Proposition~\ref{pr:IFASbounds} yields $\|g(t)\|_{\mathcal E} \leq A'$ for all $t$. Moreover, $\div g(t) = 0$, since $\div\mu^{\Pf_{n_k}}(x_l) = 0$ at all partition points and $\mu^{\Pf_{n_k}}(t) - \mu^{\Pf_{n_k}}(x_{l^{\Pf_{n_k}}(t)}) \to 0$ in $\Cc^1_b$ by \eqref{eq:equicontIFAS} and Lemma~\ref{le:interpol}. This proves \eqref{eq:decaysol} along the subsequence, and $g \in \Cc([0,T^*], L^2(\Rz^N))$ by dominated convergence.

    \bigskip

    \textbf{Step 4: The limit is a classical solution.}

    \bigskip

    Set $R(r) \coloneqq \frac{\eta^2}{2N}\Delta g(r) + \Lc_{g(r)}g(r)$. Passing to the limit in \eqref{eq:consistencyIFAS} along $(\Pf_{n_k})_k$ gives, for all divergence--free $\ph \in \Cc^{\infty}_c(\Rz^N,\Rz^N)$ and $t \in [0,T^*]$,
    \begin{equation}\label{eq:weakNS}
        \int_{\Rz^N} \ph \cdot \bigl(g(t) - g_0\bigr) \der x = \int_0^t \int_{\Rz^N} \ph \cdot R(r) \der x \der r,
    \end{equation}
    where we integrated by parts using $\div g(r) = 0$. By Step 3, we have $R \in \Cc([0,T^*], \Cc_b(\Rz^N, \Rz^N))$ and $|R(r,x)| \leq C(A')\langle x \rangle^{-N-3+\ep}$, such that for all $t \in [0,T^*]$, $w(t) \coloneqq g(t) - g_0 - \int_0^t R(r) \der r \in L^1 \cap L^{\infty}(\Rz^N, \Rz^N)$ is orthogonal to all divergence--free test fields and therefore a gradient field, i.e., $w(t) = \nabla\psi(t)$, by \cite[Lemma III.1.1]{Galdi2011}. The linearity of $s$ and the identities $s(g(t) - g_0) = g(t) - g_0$, $s(\nabla\psi(t)) = 0$ and $\div\Delta g = 0$ yield
    \begin{equation}\label{eq:mild}
        g(t) - g_0 = \int_0^t \bigl(R(r) + \nabla d(r)\bigr) \der r, 
    \end{equation}
    where $ d(r) \coloneqq -\Gamma * \div\bigl(\Lc_{g(r)}g(r)\bigr)$. Here $\div(\Lc_{g}g) \in \Cc\bigl([0,T^*], \Cc^{1,\al}_b(\Rz^N)\bigr)$ decays like $\langle x \rangle^{-2N-4+2\ep}$ and has mean zero, such that $d(r)$ is well--defined, $d(r,x) \to 0$ as $\|x\| \to \infty$, and $d \in \Cc\bigl([0,T^*], \Cc^2_b(\Rz^N)\bigr)$ with $\nabla d \in \Cc([0,T^*], L^1 \cap L^{\infty}(\Rz^N, \Rz^N))$ by Lemma~\ref{le:sweighted}. Hence the integrand in \eqref{eq:mild} is continuous on $[0,T^*] \times \Rz^N$, such that $g \in \Cc^1([0,T^*] \times \Rz^N)$ with $\partial_tg = \Lc_gg + \frac{\eta^2}{2N}\Delta g + \nabla d$. Therefore, $(g,d)$ is a classical solution of \eqref{eq:NS} on $[0,T^*]$ in the sense of Section~\ref{ch:parstructwas}, belonging to the class in the statement.

    \bigskip

    \textbf{Step 5: Uniqueness and convergence of the whole sequence.}
    
    \bigskip
    
    Uniqueness in the stated class is the classical energy argument, see \cite[Chapter 3]{MajdaBertozzi2002}, \cite{Kato1972}. Thus, every subsequence of $(\mu^{\Pf_n})_n$ has a further subsequence converging, in the sense of \eqref{eq:IFASconv}, to the same limit $g$, which proves \eqref{eq:IFASconv} for the whole sequence, and \eqref{eq:decaysol} for all $t \in [0,T^*]$.
\end{proof}

\begin{re}
    For $g_0 \in \Sc(\Rz^N, \Rz^N)$, the solution of Theorem~\ref{th:IFAS} coincides on $[0,T^*]$ with the solution of Proposition~\ref{pr:regclasssol}, by the uniqueness statement. Thus, Theorem~\ref{th:IFAS} recovers the decay $|x|^{-N-1-k}$ of $\nabla^k g$ up to the loss $\ep$, by a purely Lagrangian argument. On the other hand, Theorem~\ref{th:LIFAS} applies on every compact subinterval of the maximal existence interval of that solution, whereas the convergence of the non--linear scheme is only established up to $T^*$. Whether the (IFAS) converge on the whole existence interval, as the (LIFAS) do, is left open. In view of Remark~\ref{re:riccati}, this would require a priori bounds which do not pass through the Lipschitz norm of $\mu^{\Pf}(x_l)$ alone.
\end{re}

We conclude with two consequences of Theorem~\ref{th:IFAS} which are based on the fact that the time $T^*$ and all bounds in Proposition~\ref{pr:IFASbounds} are explicit and uniform in the viscosity.

\begin{co}[Lifespan, blow--up rate and inviscid limit]\label{co:lifespan}
    Let $N \geq 2$, $\al, \ep \in (0,1)$ and let $K$ be the constant of Proposition~\ref{pr:IFASbounds}.
    \begin{itemize}
        \item Let $\eta > 0$, $g_0 \in \Sc(\Rz^N, \Rz^N)$ be divergence--free, and let $(g,d)$ be the maximal classical solution of \eqref{eq:NS} of Proposition~\ref{pr:regclasssol}, defined on $[0,T_{\max})$. Then $g(t) \in \mathcal E_{\al,\ep}$ for all $t \in [0,T_{\max})$, $T_{\max} \geq T^*(g_0,\eta)$, and, if $T_{\max} < \infty$,
        \begin{equation}\label{eq:blowup}
            \|g(t)\|_{\mathcal E} \geq \frac{1}{2K(T_{\max} - t)} - \frac{\eta^2}{2} \qquad \text{for all } t \in [0,T_{\max}).
        \end{equation}
        \item Let $g_0 \in \mathcal E_{\al,\ep}$, $g_0 \neq 0$ and $T_0 \coloneqq \frac{1}{K(2\|g_0\|_{\mathcal E} + 1)}$. For all $\eta \in [0,1]$, let $(g^{\eta}, d^{\eta})$ be the classical solution of \eqref{eq:NS} on $[0,T_0]$ given by Theorem~\ref{th:IFAS}---such that $(g^0, d^0)$ solves the incompressible Euler equation. Then, for all $\eta \in [0,1]$ and $t \in [0,T_0]$,
        \begin{equation}\label{eq:inviscid}
            \begin{aligned}
                &\|g^{\eta}(t)\|_{\mathcal E} \leq (1 + K)^2(2\|g_0\|_{\mathcal E} + 1),\\
                &\|g^{\eta}(t) - g^0(t)\|_{L^2} \leq C(N,\ep)(1 + K)^2(2\|g_0\|_{\mathcal E} + 1) t e^{(1 + K)^2(2\|g_0\|_{\mathcal E} + 1)t}\eta^2,
            \end{aligned}
        \end{equation}
        and $g^{\eta}$ converges to $g^0$ in $\Cc\bigl([0,T_0], \Cc^1_b(\Rz^N, \Rz^N)\bigr)$ as $\eta \to 0$. 
        
        Moreover, $\|g^{\eta}(t)_i - g^{\eta}(s)_i\|_{\Wz_1} \leq C(N,\al,\ep,\|g_0\|_{\mathcal E})|t - s|$ for all $s,t \in [0,T_0]$, $i \in \llbracket 0, N-1 \rrbracket$ and $\eta \in [0,1]$.
    \end{itemize}
\end{co}
\begin{proof}
    For the first point, Proposition~\ref{pr:regclasssol} yields that $g \in \Cc([0,T_1], \Cc^3_b)$ with $|\nabla^k g(t,x)| \leq C(T_1)(1 + \|x\|)^{-N-1-k}$ for $k \leq 3$ and all $T_1 < T_{\max}$, so that $g(t) \in \mathcal E_{\al,\ep}$. Indeed, the weighted Hölder seminorm of $\nabla^2 g(t)$ is controlled by $\nabla^3 g(t)$ at distances $\leq \frac{1}{2}\min(\langle x \rangle, \langle y \rangle)$ and by the third point of Lemma~\ref{le:weights} otherwise. Fix $t_0 \in [0,T_{\max})$. 
    Theorem~\ref{th:IFAS}, applied with initial condition $g(t_0)$ and time variable $t - t_0$, provides a classical solution $(\tilde g, \tilde d)$ on $[t_0, t_0 + T^*(g(t_0),\eta)]$ with $\tilde g(t_0) = g(t_0)$. It follows $\tilde g = g$ on the common interval, and $(g,d)$ extends to a bounded classical solution on $\bigl[0, t_0 + T^*(g(t_0),\eta)\bigr)$. By maximality, $T_{\max} \geq t_0 + T^*(g(t_0),\eta) \geq t_0 + \frac{1}{K(2\|g(t_0)\|_{\mathcal E} + \eta^2)}$, and the statement follows.

    Set $A_1 \coloneqq (1 + K)^2(2\|g_0\|_{\mathcal E} + 1)$. For the second point, note that $T^*(g_0,\eta) \geq \frac{1}{K(2\|g_0\|_{\mathcal E} + \eta^2)} \geq T_0$ for $\eta \in [0,1]$, such that Theorem~\ref{th:IFAS} applies on $[0,T_0]$ for every $\eta \in [0,1]$, with $2\|g_0\|_{\mathcal E} + \eta^2 \leq 2\|g_0\|_{\mathcal E} + 1$. This gives the first bound in \eqref{eq:inviscid}. 
    
    Let $w \coloneqq g^{\eta} - g^0$ and $p \coloneqq d^{\eta} - d^0$. Then $w \in \Cc([0,T_0], L^2) \cap \Cc^1([0,T_0] \times \Rz^N)$, $\div w = 0$, $w(0) = 0$ and
    \begin{equation}
        \partial_t w = \Lc_{w}g^{\eta} + \Lc_{g^0}w + \frac{\eta^2}{2N}\Delta g^{\eta} + \nabla p.
    \end{equation}
    Now, since  
    \begin{equation}
        \begin{aligned}
            &\int_{\Rz^N} \bigl<w, \Lc_{g^0}w \bigr> \der x = \frac{1}{2}\int_{\Rz^N} \bigl< g^0, \nabla|w|^2 \bigr>\der x = 0,\\
            &\int_{\Rz^N} \bigl<w, \nabla p \bigr>\der x = -\int_{\Rz^N} p \div w \der x =0,
        \end{aligned}
    \end{equation}
    we obtain
    \begin{equation}
        \begin{aligned}
            \frac{\der}{\der t}\|w(t)\|_{L^2} &\leq \|\nabla g^{\eta}(t)\|_{\infty}\|w(t)\|_{L^2} + \frac{\eta^2}{2N}\|\Delta g^{\eta}(t)\|_{L^2}\\
             &\leq A_1\|w(t)\|_{L^2} + \frac{\eta^2}{2N}C(N,\ep)A_1,
        \end{aligned}
    \end{equation}
    where we used $|\Delta g^{\eta}(t,x)| \leq A_1\langle x \rangle^{-N-3+\ep}$, which is square integrable. Now,  Grönwall's inequality yields the second bound in \eqref{eq:inviscid}. 
    Finally, since we have $\|w(t)\|_{W^{2,\infty}} \leq 2A_1$, the inequality $\|u\|_{\infty} \leq C(N)\|u\|_{L^2}^{\frac{2}{N+2}}\|\nabla u\|_{\infty}^{\frac{N}{N+2}}$, proved as \eqref{eq:interpol}, testing with $\|u\|_{L^2}^2 \geq \frac{a^2}{4}|B(x_0,\rh)|$, and \eqref{eq:interpolholder} yield $\sup_{t \in [0,T_0]}\|w(t)\|_{\Cc^1_b} \leq C(N,\al,\ep,a_0)\,\eta^{\frac{2}{N+2}}$, which proves the convergence in $\Cc([0,T_0], \Cc^1_b)$.

    The $\Wz_1$--Lipschitz bound being independent of $\eta$ follows from \eqref{eq:equicontIFAS} by passing to the limit in \eqref{eq:IFASconv}.
\end{proof}

\section{Conclusion}
First, we have shown that the parabolic differential structure and the methods introduced by \citet{Schichl2025} are not limited to probability measures. 
On the normed Wasserstein space $\Wz_1$ of \citet{PiccoliRossiTournus2023}, the averaged flow dispersion operators and the expectational flow operators remain well--defined and source terms---inaccessible in the probability setting, since they do not conserve mass---are encoded by the source operators. 
Along divergence--free tangent vectors, the flow operators act on densities by measure--preserving flows, so that $L^{\infty}$, Hölder, $L^1$, moment and weighted Hölder norms are essentially non--increasing along discrete schemes defined by those operators. 
We expect this to be the main use of the present work: the same operators yield Euler--type schemes for the measure differential equations with sources studied in \cite{PiccoliRossi2014, Piccoli2019, PiccoliRossiTournus2023}, such as transport--diffusion equations with birth and death terms, and the compactness and stability results of Section~\ref{ch:diffprelim} provide the convergence arguments.

The second contribution of this work is the application to the incompressible Navier--Stokes equations. The non--linear scheme (IFAS) is a purely Lagrangian splitting scheme which never differentiates the unknown, and which can be read as the Euler--step analogue of the stochastic Lagrangian representation of \citet{ConstantinIyer2008}, the dispersion step being the expectation over a $2N$--point random walk. 
Theorem~\ref{th:IFAS} shows that, for initial data in the weighted Hölder space $\mathcal E_{\al,\ep}$, which represents the optimal spatial decay of Navier--Stokes flows, the (IFAS) converge in $\Wz_1$ and in $\Cc^2$ to the unique classical solution on the explicit time interval $[0,T^*]$ of Proposition~\ref{pr:IFASbounds}, with $T^* \geq \bigl(K(2\|g_0\|_{\mathcal E} + \eta^2)\bigr)^{-1}$. 
The proof rests on two observations: 
\begin{itemize}
    \item the symmetry of the dispersion step makes the diffusion cost $O(\eta^2 t)$ in weighted norms instead of $O(\eta\sqrt t)$, which is what allows the bounds to be summed along partitions of vanishing mesh (see Lemmas~\ref{le:weights}, \ref{le:eflowop} and \ref{le:fflowop}).
    \item  the Leray projection, unbounded on $L^1$ and on $W^{2,\infty}$, is bounded on a natural subspace of $\mathcal E_{\al,\ep}$ (see Lemma~\ref{le:sweighted}).
\end{itemize}
Since $T^*$ and all bounds are uniform in the viscosity, Corollary~\ref{co:lifespan} recovers, with explicit constants and uniformly in the weighted decay class, the classical lower bound on the lifespan and blow--up rate for Hölder solutions and the inviscid limit with its $O(\eta^2)$ rate \cite{Swann1971, Kato1972, Constantin1986}. 
The novelty lies in the newly introduced methods, the purely Lagrangian approach and in the uniform control of the far field. The linearised scheme (LIFAS) converges to any given classical solution (Theorem~\ref{th:LIFAS}) and yields the $L^{\infty}$, Hölder, $L^1$ and moment estimates of Corollaries~\ref{co:solestim} and \ref{co:L1estim}, in which the pressure gradient appears as the source creating and destroying mass.

Several questions remain open. The convergence of the (IFAS) beyond $T^*$, on the whole existence interval of the classical solution, comes down to closing the Lagrangian a priori bounds globally, which cannot be done through the Lipschitz norm of the approximations alone. From a numerical point of view, the (IFAS) is a projection method of Chorin type with semi--Lagrangian advection and random--walk diffusion \cite{Chorin1968, Chorin1973}. Theorem~\ref{th:IFAS} is a convergence proof for it in the Wasserstein--norm, which is the natural topology for particle approximations, but convergence rates and the spatial discretisation are left open. 


\bibliographystyle{plainnatnew}
\bibliography{MDE}

\end{document}